\documentclass[opre,nonblindrev]{informs3}
\usepackage{amssymb}
\usepackage{amsmath,amsfonts}
\usepackage{subcaption}
\usepackage{natbib}
\usepackage{bbm}
\usepackage[resetlabels]{multibib}
\newcites{app}{Appendix References}
\usepackage{algorithmic}  \usepackage{caption}
\usepackage{algorithm}
\usepackage{enumerate}
\usepackage{amsmath}

\newcommand{\R}{\mathbb{R}}

\newcommand{\hU}{\widehat{U}}

\newcommand{\Z}{\mathbb{Z}}
\newcommand{\E}{\mathbb{E}}

\newcommand{\red}[1]{\textcolor{red}{#1}}

\newcommand{\frakj}{\mathfrak{j}}
\usepackage{natbib}
\bibpunct[, ]{(}{)}{,}{a}{}{,}%
\def\bibfont{\small}%
\newcommand{\Pd}{\mathbb{P}}

\newcommand{\hX}{\widehat{X}}

\newcommand{\hV}{\widehat{V}}

\newcommand{\1}{\mathbbm{1}}

\newcommand{\bbS}{\mathbb{S}}
\newcommand{\Ex}{\mathbb{E}}

\newcommand{\bbZ}{\mathbb{Z}}

\newcommand{\mL}{\mathcal{L}}

\newcommand{\calU}{\mathcal{U}}

\newcommand{\calO}{\mathcal{O}}

\newcommand{\eProof}{\hfill  \qed }
\newcommand{\bProof}{\noindent\textbf{Proof: }}

\newcommand{\bsq}{\vrule height .9ex width .8ex depth -.1ex}

\newcommand{\bea}{\begin{equation*}}
\newcommand{\bas}{\begin{align*}}
\newcommand{\eas}{\end{align*}}
\newcommand{\eea}{\end{equation*}}
\def\approxprop{%
	\def\p{%
		\setbox0=\vbox{\hbox{$\propto$}}%
		\ht0=0.6ex \box0 }%
	\def\s{%
		\vbox{\hbox{$\sim$}}%
	}%
	\mathrel{\raisebox{0.7ex}{%
			\mbox{$\underset{\s}{\p}$}%
	}}%
}

\def\tinf{\rightarrow\infty}

\TheoremsNumberedThrough     
\ECRepeatTheorems

\usepackage[usenames,dvipsnames]{xcolor}

\usepackage[normalem]{ulem}

\begin{document}
	
	\ARTICLEAUTHORS{%
		\AUTHOR{Anton Braverman}
		\AFF{Kellogg School of Management at Northwestern University, \EMAIL{anton.braverman@kellogg.northwestern.edu}}
\AUTHOR{Itai Gurvich}
		\AFF{Cornell school of ORIE and Cornell Tech, \EMAIL{gurvich@cornell.edu}}
		\AUTHOR{Junfei Huang}
		\AFF{Chinese University of Hong Kong, \EMAIL{junfeih@cuhk.edu.hk}}}
	
\RUNAUTHOR{Braverman et. al.}

\TITLE{On the Taylor Expansion of Value Functions}
\RUNTITLE{Taylored DPs}

\TITLE{On the Taylor Expansion of Value Functions}
\RUNTITLE{Taylored DPs}

\ABSTRACT{We introduce a framework for approximate dynamic programming that we apply to discrete time chains on $\mathbb{Z}_+^d$ with countable action sets. Our approach is grounded in the approximation of the (controlled) chain's generator by that of another Markov process. In simple terms, our approach stipulates applying a second-order Taylor expansion to the value function to replace the Bellman equation with one in continuous space and time where the transition matrix is reduced to its first and second moments. In some cases, the resulting equation (which we label {\bf TCP}) can be interpreted as corresponding to a Brownian control problem. When tractable, the TCP serves as a useful modeling tool. More generally, the TCP is a starting point for approximation algorithms. We develop bounds on the optimality gap---the sub-optimality introduced by using the control produced by the ``Taylored'' equation. These bounds can be viewed as a conceptual underpinning, analytical rather than relying on weak convergence arguments, for the good performance of controls derived from Brownian control problems. We prove that, under suitable conditions and for suitably ``large'' initial states, (i) the optimality gap is smaller than a $1-\alpha$ fraction of the optimal value, where $\alpha\in (0,1)$ is the discount factor, and (ii) the gap can be further expressed as the infinite horizon discounted value with a ``lower-order'' per period reward. Computationally, our framework leads to an ``aggregation'' approach with performance guarantees. While the guarantees are grounded in PDE theory, the practical use of this approach requires no knowledge of that theory. }
\HISTORY{\emph{This version:} \today. \emph{File name:} \texttt{\jobname.tex}
}

\maketitle
\vspace*{-1cm}
\section{Introduction} Dynamic Programming (DP) is the primary tool for solving optimization problems where decisions are subject to dynamic changes in the system state. It is used in the study and practice of a variety of applications.

Deriving structural properties is typically a challenge and, computationally, the number of calculations grows exponentially with the size of the state space (the infamous ``curse of dimensionality''). Approximations are useful both as a modeling tool---for the pursuit of structural insight---and for decision-support, i.e., for the derivation of {\em computationally} efficient solutions to practical problems.

As in other classes of optimization problems---combinatorial, stochastic, etc.---approximations are often the only way to gain computational tractability for large problems. This has motivated the development of approximate dynamic programming methods; see, e.g., the books by \cite{powell2007approximate,bertsekas2007approximate}. As a modeling tool, Brownian approximations have made inroads across multiple disciplines, notably in economics, operations management, and electrical engineering. They often capture structural relationships that are inaccessible in the original, ``too'' detailed, dynamic programming problem. Yet, for a variety of reasons, these have not been widely used as a way to reduce computational complexity.

We partly bridge the gap in this paper. What we add is an approximation to dynamic programs that is inspired by {\em perturbation} techniques that were recently developed for the approximation of stationary queues by ``Brownian queues''; see \cite{Gurvich_AAP2014,BravermanDai_AAP2017,huang2016beyond}, and the additional discussion below. The seeds of the idea for extending these methods from performance analysis to optimal control appear in \cite{ata2012optimality} and \cite{huang2016beyond}. This paper seeks to start expanding those ideas---applied to heavy-traffic queues---into an accessible and generalizable framework. 	

The initial step in our approach is intuitively straightforward: we formally replace the value function in the optimality (a.k.a. Bellman) equation with its second-order Taylor expansion to obtain an equation considered over a continuous state space. As an example, consider a discrete time and space Markov chain on $\mathbb{Z}$ collecting a reward $r(x)$ when visiting state $x$ and making transitions following the stochastic matrix $P\equiv P_{x,y}$. Fixing $\alpha\in (0,1)$, the infinite horizon discounted reward
$$V(x)=\Ex_x\left[\sum_{t=0}^{\infty}\alpha^t r(X_t)\right],~x\in\mathbb{Z},$$ satisfies the functional equation
$$V(x)=r(x)+\alpha \sum_{y}P_{x,y}V(y), ~x\in\mathbb{Z},$$ which can be re-written as
$$0=r(x)+\alpha \sum_{y}P_{x,y}(V(y)-V(x))-(1-\alpha)V(x),~x\in\mathbb{Z}.$$

Applying (formally) a second-order Taylor expansion $V(y)\approx V(x)+V'(x)(y-x)+\frac{1}{2}V''(x)(y-x)^2$, we obtain the {\em differential} equation:
$$0=r(x)+ \alpha\mu(x) V'(x)+\alpha \frac{1}{2} \sigma^2(x)V''(x)-(1-\alpha)V(x),~x\in \mathbb{R} ,$$ where $\mu(x):=\Ex_x[X_1-x]=\sum_{y}P_{x,y}(y-x)$ and $\sigma^2(x):=\Ex_x[(X_1-x)^2]=\sum_{y}P_{x,y}(y-x)^2$.

When it exists, the solution $\hV$ to this Taylored equation can be interpreted as corresponding to the infinite-horizon discounted reward of a diffusion process with drift $\alpha \mu(x)$, diffusion coefficient $\alpha\sigma^2(x)$, and exponential discounting $e^{-(1-\alpha)t}$. Such an interpretation, while conceptually useful, is not mathematically necessary. Second-order Tayloring leads naturally to bounds in terms of the the third derivative of $\hV$:
$$|\hV(x)-V(x)|\leq \bar{\frakj}^3\Ex_x\left[\sum_{t=0}^{\infty}\alpha^t|D^3\hV(X_t)|_{X_t\pm \bar{\frakj}}^*\right],$$ where $|D^3\hV(X_t)|_{X_t\pm \bar{\frakj}}^*$ is the maximum of the third derivative in a neighborhood of radius $\bar{\frakj}$ around $X_t$ and $\bar{\frakj}$ is the maximal jump of the Markov chain.

This analysis of performance {\em evaluation} suggests an approach for {\em optimization}. Applying the second-order Taylor expansion to the Bellman equation
$$V(x)=\max_{u \in \mathcal{U}(x)}\left\{r(x, u)+\alpha \sum_{y}P_{x,y}^uV(y)\right\},$$ we obtain a ``Taylored'' Control Problem (TCP); see \S \ref{sec:tayloring}. Formulating the TCP is the first step. The next steps are: (1) to translate the Tayloring-induced error into bounds on optimality gaps, and (2) to build on Tayloring to propose solution algorithms.
For the development of optimality-gap bounds we draw on the theory of Partial Differential Equation (PDEs) to prove a vanishing-discount and an order-optimality result, both under suitable ``smoothness'' conditions on the primitives $\mu$, $\sigma^2$, and $r$. For suitably ``large'' initial conditions we have the following: (1) as $\alpha\uparrow 1$, the optimality gap shrinks, in relative terms, proportionally to $(1-\alpha)$, and (2) the gap can be bounded by the infinite-horizon discounted reward with an immediate-reward function that is of a lower polynomial order. It should not come as a surprise that our approach ``inherits'' some of the challenges and subtleties of PDE theory. This is reflected in the bounds in Theorem \ref{thm:explicit}, where the bound depends on the amount of time that the chain spends in ``corners'' of the state space.

From a computational perspective, because Tayloring collapses the transition matrices into $\mu$ and $\sigma$, multiple chains can induce the same TCP; they are TCP-equivalent. We rely on the TCP to transition from the original chain to another, more tractable, one. The TCP ``couples'' the two chains and supplies bounds on the approximation error.

What we are about to introduce in this paper has intimate connections to, and creates a bridge between, two somewhat disparate streams of the literature.		

\paragraph{\bf Asymptotic optimality in queues and generator comparisons.} Asymptotic optimality arguments typically require the  machinery of weak convergence. Our approach offers (in applicable cases) a simple alternative with explicit bounds. It is motivated by recent developments in queueing theory pertaining to Stein's method.\footnote{This stream of the literature is ``in its infancy'' relative to the well-developed literature on convergence-based asymptotic optimality. A key benefit of asymptotic analysis in controlled queues is so-called ``state-space collapse''---the reduction of problem dimensionality through the convergence of parts of the state space to degenerate points. A framework to incorporate such dimensionality reduction into a Stein-type analysis is still absent.} In performance analysis (i.e., for a given control), Stein's method allows us to bound directly---without resorting to convergence arguments---the (impressive) ``proximity'' between the stationary distribution of a queueing system and its Brownian approximation by comparing their transition probabilities (or, more precisely, their {\em generators})---that of the Markov chain and that of a suitable diffusion process; see \cite{Gurvich_AAP2014} and \cite{BravermanDai_AAP2017}. While the use of the language of generators is mathematically natural, it is simpler and conceptually useful to view this as Taylor expansion applied to equations that characterize stationary performance and/or optimality conditions.

In the papers cited above the problem primitives and the state space are scaled according to a {\em heavy-traffic} scaling; having, e.g., the arrival-rate increase to infinity or utilization increase to $100\%$. In this paper, the parameters {\em are fixed} and no space scaling is used. This simplifies the analysis of the TCP and, in turn, the derivation of optimality-gap bounds; see \S \ref{subsec:bounds}. It also allows us to consider applications that lack a natural notion of scaling. Nevertheless, specializing our results to queueing examples and relating the discount factor to the utilization does shed some light on the nature of our results; see \S \ref{sec:HT}.

Transitioning from performance analysis, as considered in earlier papers, to {\em controlled} chains as we do here, is like considering a family of generators (``indexed'' by the control) instead of a single generator. One can interpret the Taylored equation as the  Hamilton-Jacobi-Bellman (HJB) equation for a suitable Brownian control problem. The relation we seek to uncover is based not on process-limit theory but rather on first principles, namely, the Tayloring of the value function.

{\bf Approximate dynamic programming (ADP).} Approximate value or policy iteration typically starts with the choice of a base from which to construct a candidate value function. The queueing-approximations literature teaches us that, as a heuristic, the value function of a suitable Brownian control problem is a good candidate for a base function; such an approach is taken, for example, in \cite{chen2009approximate}. Our analysis supports this approach: we establish that the TCP solution, even  taken as the sole item in the base, yields an approximation whose performance is related to properties of a closely related differential equation.

Algorithmically, our {\em Taylored Approximate Policy Iteration} (TAPI) algorithm is a modification of policy iteration, in which the policy evaluation portion of iteration-$k$ requires solving a {\em linear} PDE to get an approximate value function $V^{(k)}$, which is subsequently plugged into a policy improvement step (an optimization problem that does not require the solution of a PDE) to produce $u^{(k+1)}$ and so on. The linear PDE can be solved via Finite Difference (FD) or other PDE discretization-based solution methods. The coarser the discrete grid the more efficient is the computation.

An alternative to FD in the implementation of TAPI is to build on the Taylored equation as an intermediate step---a translator---between Bellman equations corresponding to two {\em TCP-equivalent} chains, i.e.,  that induce the same TCP. Given a controlled chain, one possible construction of a TCP-equivalent one is inspired by the transformation put forth in \cite{kushner2013numerical} and \cite{dupuis1998rates} that relates the differential equation to a control problem for a Markov chain, henceforth referred to as the K-D chain---one with a smaller state space and a simpler transition structure. In contrast with the infinitesimal view inherent to the K-D approach (where one takes the discretization to $0$ to approximate continuous state space), we use it with {\em coarse discretization} so that the new Bellman equation can be viewed as an aggregation method---where the state space is reduced to a coarser grid with ``super states''; for existing aggregation ideas see, e.g., \cite[Ch. 6]{bertsekas2007approximate}.

Our Tayloring approach to approximate dynamic programming stands on strong mathematical footing. The gap introduced by using a TCP-equivalent chain can be bounded via the (suitably integrated) third derivative of the PDE solution. From a computational viewpoint, although the algorithm that we propose is not entirely immune to the curse of dimensionality, it pushes computational barriers. This paper introduces the framework, and provides analytical support and initial numerical evidence. A full account of algorithms and computational benefits is left for future work; see \S \ref{sec:conclusions}.

\vspace*{0.5cm}

{\bf Notation:} We use the standard notation $\mathbb{R}_+^d$ for the positive orthant in $\mathbb{R}^d$ and  $\mathbb{R}_{++}^d$ for its interior---the space of strictly positive $d$-dimensional vectors. For a set $\Omega\subseteq \mathbb{R}_+^d$, $\partial{\Omega}$ denotes its boundary. In particular, $\partial \mathbb{R}_+^d=\mathbb{R}_+^d\backslash\mathbb{R}_{++}^d=\bigcup_{i=1}^d\mathcal{B}_i$ where $\mathcal{B}_i:=\{x\geq 0 :x_i=0\}$. The standard Euclidean norm is denoted by $|\cdot|$ and for $x\in\mathbb{R}_+^d$ and $\epsilon>0$, we denote by $x\pm \epsilon$ the {\em set} $\{y\in \mathbb{R}_+^d:|y-x|\leq \epsilon\}$. For a function $f:\mathbb{R}_+^d\to \mathbb{R}$ and a subset $\Omega\subseteq \mathbb{R}_+^d$, we let
$|f|_{\Omega}^*=\sup_{y\in \Omega} |f(y)|$ and, for $\beta\in (0,1]$, $$[f]_{\beta,\Omega}^*=\sup_{y,z\in \Omega}  \frac{|f(y)-f(z)|}{|y-z|^{\beta}}.$$ If $f$ is twice continuously differentiable we write $f_i(x)=\frac{\partial}{\partial x_i}f(x)$ and $f_{ij}(x)=\frac{\partial^2}{\partial x_i\partial x_j}f(x)$. We use $Df(x)$ for the Gradient vector (whose elements are $f_i(x)$) and $D^2f(x)$ for the Hessian matrix (whose elements are $f_{ij}(x)$). We use the standard notation $\mathcal{C}^{2}(\Omega)$ for the family of twice continuously differentiable functions over $\Omega$, and $\mathcal{C}^{2,\beta}(\Omega)$ is the subset of $\mathcal{C}^2(\Omega)$ whose members have a  second derivative that is H{\"o}lder continuous on $\Omega$ with exponent $\beta$; $\beta=1$ corresponds to Lipschitz continuity on $\Omega$.  In this paper when we speak of a {\em solution} to a differential equation we mean that in the classical sense.

Throughout, to simplify notation, we use $\gamma,\Gamma$ to denote Hardy-style constants that may change from one line to the next and that do not depend on the discount factor $\alpha$ or on the state $x$.

\section{Tayloring the Bellman equation\label{sec:tayloring}}

Consider an infinite-horizon discounted Markov Decision Problem on $\mathbb{Z}_+^d$: $$V_{\ast}^{\alpha}(x):=\max_{U}\Ex_{x}^U\left[\sum_{t=0}^{\infty}\alpha^t r(X_t,U(X_t))\right],$$ where $r(x,u)$ is the reward collected at state $x$ under a control $u$. A stationary policy $U=\{U(x), x\in\mathbb{Z}_+\}$ has the property that $U(x)\in\calU(x)$, where $\mathcal{U}(x)$ is the set of actions allowed in state $x$. We assume that $\mathcal{U}(x)$ is {\em discrete} (possibly countably infinite) and it is an intersection of a polyhedron (that can depend on $x$) and a discrete set $\mathbb{D}$ that does not depend on $x$, i.e., for some $m$, $\mathcal{U}(x)=\{u\in\mathbb{R}^m:Au\leq b(x)\} \cap \mathbb{D}$, where $b(x)$ is defined for all $x\in\mathbb{R}_+^d$. We let $\mathbb{U}=\times_{x\in\mathbb{Z}_+^d} \mathcal{U}(x).$

Given $x\in\mathbb{Z}_+^d$ and $u\in \mathcal{U}(x)$, we write $r_{u}(x)=r(x,u)$ and let $P^{u}_{x,y}$ be the probability of transitioning from $x$ to $y$ under an action $u\in\calU(x)$. We write $\Ex_x^{U}[\cdot]$ for the expectation with respect to the law of the $U$-controlled Markov chain $(X_t,t\geq 0)$ with the initial state $x$;  $\Ex_x^u[\cdot]$ is the expectation with respect to the law $P^u_{x,\cdot}$.

Under standard conditions, e.g., \cite[\S  1.4]{bertsekas2007approximate}, $V_{\ast}^{\alpha}(x)$ solves the Bellman equation \begin{align} V(x)&=\max_{u\in\mathcal{U}(x)}\Big\{r_{u}(x)+\alpha P^{u}V(x)\Big\}=
\max_{u\in\mathcal{U}(x)}\Big\{r_{u}(x)+\alpha\Ex_x^u[V(X_1)]\Big\},\label{eq:Bellman}\end{align} where we use the operator notation $P^{u}V(x)=\sum_{y}P^{u}_{x,y}V(y)=\Ex_x^u[V(X_1)]$. Subtracting $V(x)$ on both sides of the Bellman equation we have \be 0= \label{eq:discrete}\max_{u\in\mathcal{U}(x)}\Big\{r_{u}(x)+\alpha(P^{u}V (x)-V(x))-(1-\alpha)V(x)\Big\},~ x\in \mathbb{Z}_+^d.\ee Pretending that $V$ is extendable to $\mathbb{R}_+^d$ and twice-continuously differentiable there, we have
\begin{align*} ((P^{u}-I)V)(x)&=\sum_{y}P^{u}_{x,y}V(y)-V(x)
	\\ & \approx \sum_{y} P^{u}_{x,y}\left(\sum_{i}V_i(x)(y_i-x_i) +\frac{1}{2}\sum_{i,j}V_{ij}(x)(y_i-x_i)(y_j-x_j)\right).
\end{align*} Defining
\begin{align} (\mu_{u})_i(x)&:=\Ex_x^{u}[(X_{1})_i-x_i]=\sum_{y} P^{u}_{x,y}(y_i-x_i),~i=1,\ldots,d,~ \mbox{ and }\\
	(\sigma^2_{u})_{ij}(x)&:=\Ex_x^u[((X_1)_i-x_i)((X_1)_j-x_j)] =\sum_{y}P^{u}_{x,y}(y_i-x_i)(y_j-x_j),~i,j=1,\ldots,d,\label{eq:mudefin}\end{align}  for all $x\in\mathbb{Z}_+^d$ and $u\in\mathcal{U}(x)$ (and extending these to $\mathbb{R}_+^d$; see the discussion after Assumption \ref{asum:1}), we arrive at
\be\label{eq:taylored00}
r_{u}(x)+\alpha(P^{u}V (x)-V(x))-(1-\alpha)V(x)\approx
r_{u}(x)+\alpha\mathcal{L}_{u}V(x)-(1-\alpha)V(x),~ x\in \mathbb{R}_+^d,
\tag{2nd order Taylor}
\ee
where \begin{align*}
	\mathcal{L}_{u}V(x)&=\sum_{i}(\mu_{u})_i(x)V_i(x)+\frac{1}{2}\sum_{i,j}(\sigma_{u}^2)_{ij}(x)V_{ij}(x)\\&=
	\mu_{u}(x)'DV(x) +\frac{1}{2}trace(\sigma_{u}^2(x)'D^2 V(x)).
\end{align*}
This suggests, heuristically at this stage, replacing \eqref{eq:Bellman} with
\begin{align}
0=\max_{u \in \mathcal{U}(x)}\{r_u(x)+\alpha \mathcal{L}_uV (x)-(1-\alpha)V(x)\}, ~ x\in \mathbb{R}_+^d.\label{eq:TCP0}
\end{align}
A solution to \eqref{eq:TCP0} is a pair $(\hU_*^{\alpha}(x),\hV_*^{\alpha}(x))$, where $\hU_*^{\alpha}(x)$ is the maximizer. We can consider $\hV_*^{\alpha}(x)$ as an approximation for the true optimal value $V_*^{\alpha}(x)$, and the restriction of the maximizer $\hU_*^{\alpha}(x)$ to  $\mathbb{Z}_+^d$ gives a feasible control for the original chain, allowing us to refer to the $\hU_*^{\alpha}$-controlled chain.

Implicit in this derivation is an extension of $r_u$, $\mu_u$, and $\sigma^2_u$ from $\mathbb{Z}_+^d$ to $\mathbb{R}_+^d$. We first impose the requirement that the primitives have natural extensions from $\{(x,u):x\in\mathbb{Z}_+^d, u\in\mathcal{U}(x)\}$ to
$\{(x,u):x\in\mathbb{Z}_+^d, u\in\mathbb{D}\}$.

\begin{assumption}[primitives] There exist functions \be f_r(x,u), ~f_{\mu}(x,u),~f_{\sigma}(x,u), ~x\in\mathbb{Z}_{+}^d, u\in\mathbb{D},\label{eq:fdefin}\ee such that $r_u,\mu_u$, and $\sigma^2_u$ are the restrictions of these functions to $x\in\mathbb{Z}_+^d$ and $u\in\mathcal{U}(x)$ and satisfy the following properties: (i) $f_r$ is locally Lipschitz in $\mathbb{Z}_+^d$ (uniformly in $u$), and (ii)
The functions $f_{\mu}$ and $f_{\sigma}$ are globally bounded and Lipschitz uniformly in $u$, i.e., there exists $L>0$ (not depending on $u$) such that
$$|f_{\mu}(\cdot,u)|_{\mathbb{Z}_+^d}^*+[f_{\mu}(\cdot,u)]_{1,\mathbb{Z}_+^d}^*,~|f_{\sigma}(\cdot,u)|_{\mathbb{Z}_+^d}^*+[f_{\sigma}(\cdot,u)]_{1,\mathbb{Z}_+^d}^*\leq L.
$$	 Finally, (iii) $f_{\sigma}(x, u)$ (and, in turn, its restriction $\sigma^2_u(x)$) satisfies the ellipticity condition: there exists $\lambda>0$ (not depending on $u$) such that 	
	\be \tag{elliptic}
\lambda^{-1}|\xi|^2\geq \sum_{ij}\xi_{i}\xi_j (f_{\sigma})_{ij}(x,u)\geq \lambda|\xi|^2, \mbox{ for all } \xi\in\mathbb{R}^d,~x\in\mathbb{Z}_+^d,~u\in\mathbb{D}.\ee
		\label{asum:1}
	\end{assumption}

Under the Lipschitz requirement in Assumption \ref{asum:1},  the McShane-Whitney extension theorem (\cite{mcshane1934extension}) constructs an explicit extension to $\mathbb{R}_{+}^d$ that is itself Lipschitz continuous with the same constant $L$ (or locally Lipschitz in the case of $f_r$). It is sometimes convenient to leave a discontinuity at the boundary (see the OD-boundary condition below and Example \ref{example:SSQ0}). Importantly, the computational algorithm that we use relies on the extension to $\mathbb{R}_+^d$ only (if at all) on the boundary. Continuity properties of these extensions do matter for our analytical results. {\em Henceforth, $f_r(\cdot,u)$, $f_{\mu}(\cdot,u)$ and $f_{\sigma}(\cdot,u)$ are the extensions to $\mathbb{R}_+^d$. }

Finally, because every discrete state space can be embedded in $\mathbb{Z}_+$, it is fair to ask what are the requirements that we impose on the original problem. As in Assumption \ref{asum:1}, these requirements are and will be stated as constraints on $\mu$ and $\sigma^2$. Our optimality-gap bounds require, for example, that the optimally controlled chain has bounded jumps; see Theorems \ref{thm:diffusion} and \ref{thm:explicit}. The bound's magnitude, in turn, depends on the maximal jump-size as it depends on the Lipschitz constant $L$ in Assumption \ref{asum:1}. The embedding of a two-dimensional chain into one dimension might induce $\mu_u$ and/or a maximal jump size that are significantly larger than in the original two-dimensional model.

Conceptually speaking, our approach is relevant to settings where (1) there is a natural meaningful metric on the state space so that $\mu$ can be interpreted as the average step size starting at $x$, (2) one can speak of large and small initial states, and (3) boundaries have physical meaning. Thus, for example, inventory and queueing problems are natural candidates for this approach, but a Markov chain where the states are colors or letters might not be.

\subsection{Boundary conditions} Equation \eqref{eq:TCP0}, while well defined, poses a challenge insofar as we want to apply existing PDE theory as exposed, e.g., in \cite{gilbarg2015elliptic} and \cite{lieberman2013oblique}. The theory covers mostly first-order conditions on the boundary, i.e., those where either $DV$ or $V$ appear but not $D^2V$. We will consider two such conditions: (i) first-order Tayloring, and (ii) an oblique-derivative condition that supports second-order Tayloring on the boundary.\footnote{The discussion of boundary conditions is unnecessary if the state space is $\mathbb{Z}^d$---as in the inventory example in \S \ref{sec:examples}. But even in these, computation requires truncating the state space, making boundary conditions relevant.}

{\bf First-order Tayloring (FOT) boundary:} Applying first-order Tayloring in boundary states, i.e., replacing $V(y)-V(x)\leftarrow  DV(x)'(y-x)$ for $x\in \partial \mathbb{R}_+^d$, leads to
\begin{align} 0&=\max_{u \in \mathcal{U}(x)}\left\{
	r_{u}(x)+\alpha\mathcal{L}_{u}V(x)-(1-\alpha)V(x)\right\},~ x\in \mathbb{R}_{++}^d,\label{eq:FOT1}\\
	0& =\max_{u\in\mathcal{U}(x)}\left\{ r_{u}(x)+\alpha\mu_{u}(x)'DV(x)-(1-\alpha)V(x)\right\},~~ x\in \partial\mathbb{R}_+^d.\label{eq:FOT2}
\end{align}

We say that the FOT boundary condition is {\em control-independent} if $\mu_u(x)\equiv \mu(x)$ for all $x\in\partial \mathbb{R}_+^d$. In that case the maximizer on the boundary $\hU_*(x), x\in\partial \mathbb{R}_+^d$ does not depend on the value of $\hV_*$ and $D\hV_*$ there.

{\bf Oblique-derivative (OD) boundary:} Under certain assumptions on the behavior of $\mu$ near the boundary, certain first-order boundary conditions imply that \eqref{eq:TCP0} also holds (as a second-order equation) on the boundary. Informally, suppose that there exists a vector $\eta(x)$ such that for $y$ close to a boundary point $x\in\partial \mathbb{R}_+^d$ and all $u\in\mathbb{D}$
$$f_{\mu}(y,u)-f_{\mu}(x,u)\approxprop \eta_i(x),$$ i.e., that the boundary change in the drift is approximately proportional to $\eta$. Then, $(\hU_*,\hV_*)$, with $\hV_*\in \mathcal{C}^2(\mathbb{R}_+^d)$, that solves the {\em OD-boundary TCP}:
\begin{align} &  0=\max_{u \in\mathcal{U}(x)}\{ r_{u}(x)+\alpha\mathcal{L}_{u}V(x)-(1-\alpha)V(x)\},~ x\in \mathbb{R}_{++}^d, \label{eq:OD1} \\
&	0 = \eta(x)'DV(x),~x\in \partial\mathbb{R}_+^d,
\label{eq:OD2} \end{align}
also solves the 2nd-order TCP \eqref{eq:TCP0}. See Lemma \ref{lem:oblique} in the appendix for the formal statement.

An advantage of the TCP with OD boundary condition is its interpretability as the HJB of a control problem for a reflected diffusion; see, e.g., \cite{borkar2004ergodic}. FOT-boundary, in contrast, imposes fewer structural requirements. While queueing settings provide an intuitive way to identify $\eta$---see Example \ref{example:SSQ} and \S \ref{sec:exmultipool}---FOT is more direct and requires less context-specific expertise. It does come, however, at the cost of weaker bounds; see Remark \ref{rem:firstorder}.

\begin{example}[A discrete-time single-server queue]\label{example:SSQ0} Consider a controlled random walk on $\mathbb{Z}_+$ where, for $x\geq 1$, $P_{x,x-1}^u=u, P^u_{x,x+1}=1-u$ and $P^u_{0,1}\equiv 1$. We take $\mathcal{U}(x)=\mathbb{D}=[0,1]\cap \mathbb{Q}$ ($\mathbb{Q}$ denotes the rational numbers) for all $x\in\mathbb{Z}_+$. Then, $\mu_u(x)=1-2u=:f_{\mu}(x,u)$ for $x\geq 1$ and $f_{\mu}(0,u)=1$. Also,  $\sigma_u^2(x)\equiv 1=:f_{\sigma}(x,u)$.
	We use a reward function that penalizes for large states (holding cost) and for speedy service (effort cost) $r_u(x)=-x^4-\frac{c_s}{1-u}$ where $c_s>0$. 	
	
	We use the discontinuous extension for $f_{\mu}(x,u)$ that has $f_{\mu}(x,u)=1-2u$ for all $x>0$ and $f_{\mu}(0,u)=1$, so that $f_{\mu}(0+,u)-f_{\mu}(0,u)=-2u\propto -1$ and the OD-boundary condition is $V'(0)=0$. This condition---familiar from performance equations for reflected Brownian motion \cite[\S 6.3]{Har:2013}---finds a natural justification in Lemma \ref{lem:oblique}: If the control $\hU_*$ is continuous at $0$, then per Lemma \ref{lem:oblique}, a solution to
	\begin{align*} 0&=\max_{u \in \mathcal{U}(x)}\left\{r_u(x)+\alpha(1-2u)V'(x)+\frac{\alpha}{2}V''(x)-(1-\alpha)V(x)\right\},~ x>0,\\
		0&=V'(0),\end{align*} satisfies \eqref{eq:TCP0} at $x=0$.
In this example, the FOT-boundary condition reduces to the (control-independent) equation $$0=\max_{u\in\mathbb{D}}\{r_u(0)+\alpha V'(0)-(1-\alpha)V(0)\}=-c_s+\alpha V'(0)-(1-\alpha)V(0).$$ \hfill \bsq
	
\end{example}

\noindent {\bf State-space truncation and boundary conditions:} For computational algorithms the state space must be truncated. We will impose the truncation of $\mathbb{Z}_+^d$ to the square $\mathbb{S}_M=\{x\in \mathbb{Z}_+^d: \max_{i}x_i\leq M\}$. The boundary conditions for the TCP will depend on the way we truncate the state space in the original chain. It is natural to define the transition probabilities for $x\in\mathbb{S}_M$ by
$$\tilde{P}^u_{x,y}=\left\{\begin{array}{ll} 0 & \mbox{ for } y\notin \mathbb{S}_M,\\ \frac{P^u_{x,y}}{\sum_{z\in\mathbb{S}_M}P^u_{x,z}}, & \mbox{ otherwise. }\end{array}\right.$$ In the random walk of Example \ref{example:SSQ0}, this simply means $\tilde{P}_{M,M+1}^u=0$ and $\tilde{P}^u_{M,M-1}=1$, which leads naturally to the OD boundary condition $V'(M)=0$.

\subsection{The Initial Tayloring Bound\label{sec:initial}}

In what follows, for a fixed stationary policy $U$, and a function $f:\mathbb{Z}_+^d\to\mathbb{R}$, we write
$$V_{U}^{\alpha}[f](x)=\Ex_x^{U}\left[\sum_{t=0}^{\infty}\alpha^t f(X_t)\right].$$ We drop the argument $f$ when the immediate reward function is $r_u(x)$ and clear from the context. Thus, for example, $V^{\alpha}_{\hU_*}(x)$ is the value under the policy $\hU_*$ with the reward function $r_u(x)$.

Given a stationary policy $U$, we define $\mathfrak{j}_{U}$ to be the smallest integer (allowing for infinity) such that, for all $x,y\in\mathbb{Z}_+^d$ with $|y-x|>\frakj_U$,    $\Pd_{x,y}^{U(x)}=0$. We say that the chain has uniformly bounded jumps if $$\bar{\frakj}:=\sup_{U\in \mathbb{U}}\frakj_U<\infty.$$

Note that the controls $\hU_*$ and $U_*$ are likely to depend on $\alpha$, but for notational convenience we do not make this dependence explicit.

\begin{theorem}[initial bound with 2nd-order Tayloring at the boundary] Fix $\alpha\in(0,1)$ and suppose that there exists a solution $(\hU_*,\hV_{\ast})$ to \eqref{eq:TCP0} with $\hV_{\ast} \in \mathcal{C}^{2,\beta}(\mathbb{R}^d_{+})$ for some $\beta\in(0,1)$. Suppose further that  $\frakj_{\hU_*},\frakj_{U_*}<\infty$ and that $|\hV_{\ast}(x)|\leq \Gamma(1+|x|^m)$ for some $m$ and $\Gamma$ (that can depend on $\alpha$). Then, for $x\in\mathbb{Z}_+^d$,
	\begin{align}\nonumber
		&\left(|\hV_{\ast}(x)-V_{\ast}^{\alpha}(x)|\vee \left|V_{\hU_*}^{\alpha}(x)-V_{\ast}^{\alpha}(x)\right|\right) \\
		& \qquad \leq \frakj_{\hU_*}^{2+\beta}\vee \frakj_{U_*}^{2+\beta} \left(\Ex_x^{\hU_*}\left[\sum_{t=0}^{\infty} \alpha^t[D^2\hV_{\ast}]_{\beta,X_t\pm \frakj_{\hU_*}}^*\right]+\Ex_x^{U_*}\left[\sum_{t=0}^{\infty} \alpha^t [D^2\hV_{\ast}]_{\beta,X_t\pm \frakj_{U_*}}^*\right]\right).\label{eq:initial}
	\end{align}	\label{thm:diffusion}
\end{theorem} \vspace{0.2cm}

\begin{remark}[performance approximation] {\em We make the obvious observation that Theorem \ref{thm:diffusion} applies to the performance analysis of a given control. Fixing a control $U$ is the same as taking control sets $\mathcal{U}(x)$ that contain the single action $U(x)$. Equation \eqref{eq:initial} reduces to
		\begin{align}
	&|\hV(x)-V_{U}^{\alpha}(x)| \leq \frakj_{U}^{2+\beta} \Ex_x^{U}\left[\sum_{t=0}^{\infty} \alpha^t[D^2\hV]_{\beta,X_t\pm \frakj_{U}}^*\right].\label{eq:initialperf}
	\end{align}  In this case the TCP is a linear PDE. }	\hfill \bsq
	\end{remark}

\begin{remark}[unbounded jumps] \label{rem:unboundedjumps} {\em The bound can be easily adjusted to unbounded jumps with suitable finite moments. In this case the right-hand side of \eqref{eq:initial} takes the form
	\begin{equation*}
\begin{split}  \Ex_x^{\hU_*}\left[\sum_{t=0}^{\infty}\alpha^t \Ex_{X_t}[|\Delta_{X_t}|^{2+\beta}[D^2\hV]_{\beta,X_t\pm |\Delta_{X_t}|}^*]\right]+\Ex_x^{U_*}\left[\sum_{t=0}^{\infty}\alpha^t \Ex_{X_t}[|\Delta_{X_t}|^{2+\beta}[D^2\hV]_{\beta,X_t\pm |\Delta_{X_t}|}^*]\right],
	\end{split}
	\end{equation*} where $\Delta_{X_t}=X_{t+1}-X_t$; see the proof of Theorem \ref{thm:diffusion}.}\hfill \bsq
\end{remark}

\begin{remark}[FOT boundary]\label{rem:firstorder} {\em With first order Tayloring on the boundary, \eqref{eq:initial} is replaced with \begin{align*}
		&\left(|\hV_{\ast}(x)-V_{\ast}^{\alpha}(x)|\vee \left|V_{\hU_\ast}^{\alpha}(x)-V_{\ast}^{\alpha}(x)\right|\right) \\
		& \qquad \leq \frakj_{\hU_*}^2\vee \frakj_{U_*}^2 \left(\Ex_x^{\hU_*}\left[\sum_{t=0}^{\infty} \alpha^t|\mathfrak{e}[\hV_{\ast},\frakj_{\hU_*}]|_{\beta, X_t\pm \frakj_{\hU_*}}^*\right]+\Ex_x^{U_*}\left[\sum_{t=0}^{\infty} \alpha^t |\mathfrak{e}[\hV_{\ast},\frakj_{U_*}]|_{\beta, X_t\pm \frakj_{U_*}}^*\right]\right),
	\end{align*} where
	$$|\mathfrak{e}[f,z]|_{\beta,\Omega}^*=z^{\beta} [D^2f]_{\beta,\Omega}^*+\sum_{i}\1\{x\in\mathcal{B}_i\}\sum_{j\neq i}|f_{ij}|_{\Omega}^*.$$ Relative to \eqref{eq:initial}, the second derivative on the boundary factors into the optimality gap.\footnote{ In addition, the H\"{o}lder bounds we state for the second derivative are weaker in the case of FOT boundary condition; see Lemmas \ref{lem:PDEbounds} and \ref{lem:PDEboundsFOT}.  } }  \hfill \bsq
\end{remark}

Theorem \ref{thm:diffusion} is a starting point. It assumes the existence of a smooth solution and it leaves unspecified the magnitudes of the second-derivative's H{\"o}lder exponents. We will address both concerns in subsequent sections.

\paragraph{\bf Towards computability: TCP-equivalent chains.} The primitives of the MDP are the reward function(s), $r_u$, the transition matrices $P^u$---from which we build $\mu_u$ and $\sigma^2_u$---and the discount factor $\alpha\in (0,1)$.

There are multiple MDPs (or primitives) that induce the same TCP. Specifically, consider an MDP with the same state and action spaces. Let $\{\widetilde{P}^u\}$ be a family of transition matrices and $\tilde{\alpha}(x)\in (0,1)$ be a (possibly state-dependent) discount factor that jointly satisfy the constraints
$$\sum_{y}\widetilde{P}_{x,y}^{u}(y_i-x_i)=\frac{\alpha(1-\widetilde{\alpha}(x))}{\widetilde{\alpha}(x)(1-\alpha)}(\mu_u)_i(x),~ \sum_{y}\widetilde{P}_{x,y}^{u}(y_i-x_i)(y_j-x_j)=\frac{\alpha(1-\widetilde{\alpha}(x))}{\widetilde{\alpha}(x)(1-\alpha)}(\sigma^2_u)_{ij}(x).$$
and take the reward function  $\tilde{r}_u(x)=\frac{1-\tilde{\alpha}(x)}{1-\alpha}r_{u}(x)$.

These ``tilde'' primitives then induce {\em the same} TCP as  the original primitives. The two chains are {\em TCP-equivalent}. Let $\tilde{U}_*$ be the optimal policy for this new optimal control problem. It then follows that
\begin{align}\label{eq:TCPequivBound} |V_{*}^{\alpha}(x)-\widetilde{V}_{\ast}^{\alpha}(x)|\leq & \Gamma \left(\Ex_x^{\hU_*}\left[\sum_{t=0}^{\infty} \alpha^t[D^2\hV_{\ast}]_{\beta,X_t\pm \frakj_{\hU_*}}^*\right]+\Ex_x^{U_*}\left[\sum_{t=0}^{\infty} \alpha^t [D^2\hV_{\ast}]_{\beta,X_t\pm \frakj_{U_*}}^*\right]\right.\\ & \left.~~+ \Ex_x^{\hU_*}\left[\sum_{t=0}^{\infty} \bar{\alpha}^t[D^2\hV_{\ast}]_{\beta,\widetilde{X}_t\pm \frakj_{\hU_*}}^*\right]+\Ex_x^{\tilde{U}_*}\left[\sum_{t=0}^{\infty} \bar{\alpha}^t [D^2\hV_{\ast}]_{\beta,\widetilde{X}_t\pm \frakj_{\tilde{U}_*}}^*\right]\right),\nonumber
\end{align} where $\frakj_{\tilde{U}_*}$ is the maximal jump of the chain $\widetilde{X}$ under the policy $\tilde{U}_*$ and $\bar{\alpha}_t=\sup_{x\in\mathbb{Z}_+^d} \tilde{\alpha}(x)$, and $\Gamma$ is an appropriate constant that depends on $\frakj_{\hU^*},$ $\frakj_{U^*}$, and $\beta$.

Among all TCP-equivalent chains it is reasonable to look for one that introduces significant computational benefits. There are substantial degrees of freedom in making this choice. The state-space of the new chain could, for example, be a {\em strict} subset of that of $X$. The K-D chain that we use in \S \ref{sec:TAPI} has this property.

This transition from one Markov chain to a different, but TCP-equivalent one does not require solving any continuous state-and-time control problem. The TCP merely serves as the basis for optimality-gap guarantees. What we pursue next is making these guarantees more explicit.

\subsection{Explicit bounds\label{subsec:bounds}}

The following example captures in a simple setting the essential ingredients of the forthcoming optimality-gap bounds.

\begin{example}[The discrete queue revisited] \label{example:SSQ}
	In the setting of Example \ref{example:SSQ0}, let us {\em fix the control} to $U(x)\equiv 1/2$; see Remark \ref{rem:ssqbound} and Example \ref{example:pdescaling} for the full control version. The OD-boundary TCP is given by
	 \begin{align*}
		0&=-x^4-\frac{c_s}{1-U(x)} +\alpha (1-2U(x))V'(x)+\alpha\frac{1}{2}V''(x)-(1-\alpha)V(x),~ x>0,\\
		0 & = V'(0),
	\end{align*} and admits the unique solution\footnote{The FOT solution is equal to $\hV_U(x)=g^{\alpha}(x)$, 
		where $|D^3g^{\alpha}(x)|\leq \Gamma (1-\alpha)$ and $|D^2g^{\alpha}(0)|\leq \Gamma \sqrt{1-\alpha}$, which does not change the conclusion of this discussion.}

	$$\hV_{U}(x)=-\frac{x^4}{1-\alpha }-\frac{6 \alpha  x^2}{(1-\alpha )^2}-\frac{6 \alpha ^2}{(1-\alpha )^3}-\frac{2c_s}{1-\alpha },~x\geq 0,$$
so that $$[D^2\hV_U]_{1,[0,x]}^*\leq |D^3\hV_U|_{[0,x]}^*\leq \frac{24 x}{1-\alpha},~ x\geq 0.$$
Because the maximal jump is 1 ($\bar{\frakj}=1$), Theorem \ref{thm:diffusion} and equation \eqref{eq:initialperf} imply that
	\begin{align*}
		\left|\hV_{U}(x)-V_{U}^{\alpha}(x)\right|\leq
		\Ex_x^U\left[\sum_{t=0}^{\infty}\alpha^t
		|D^3\hV_U|_{X_t\pm 1}^*\right]\leq 24\Ex_x^{u}\left[\sum_{t=0}^{\infty} \alpha^t\frac{X_t+1}{1-\alpha}\right], ~ x\in \mathbb{Z}_+,
	\end{align*}
	where, recall, $V_U^{\alpha}(x)$ is the infinite horizon discounted reward with the immediate reward $r_u$ and under the policy $U$. For all $x\geq 0$, $\frac{x}{1-\alpha}\leq (1-\alpha)x^4+\frac{1}{(1-\alpha)^2}$, so that \begin{align*} \Ex_x^{U}\left[\sum_{t=0}^{\infty} \alpha^t \frac{X_t+1}{1-\alpha}\right]&\leq
		\Ex_x^{U}\left[\sum_{t=0}^{\infty} \alpha^t \left((1-\alpha)X_t^4 +\frac{1}{(1-\alpha)^2}\right)\right]+\frac{1}{(1-\alpha)^2}\\&\leq (1-\alpha)|V_{U}^{\alpha}(x)|+\frac{2}{(1-\alpha)^3}.\end{align*} We claim that $|V_{U}^{\alpha}(x)|\geq \frac{\gamma}{(1-\alpha)^4}$ for all $x\geq \frac{1}{1-\alpha}$, so that
	\begin{align*}
		\left|\hV_{U}(x)-V_{U}^{\alpha}(x)\right|\leq \Gamma(1-\alpha)|V_{U}^{\alpha}(x)|,\end{align*} for all such $x$;
	see Corollary \ref{cor:alphaopt} and its proof.
	Furthermore, because $\frac{x}{1-\alpha}\leq x^3+\frac{1}{(1-\alpha)^{\frac{3}{2}}}$ for all $x\geq 0$, we have
	\begin{align*} \Ex_x^{u}\left[\sum_{t=0}^{\infty} \alpha^t \frac{X_t+1}{1-\alpha}\right]&\leq
		\Ex_x^{u}\left[\sum_{t=0}^{\infty} \alpha^t \left(X_t^3 +\frac{1}{(1-\alpha)^{\frac{3}{2}}}\right)\right]+\frac{1}{(1-\alpha)^2}\\&\leq V_{U}^{\alpha}[f_3](x)+\frac{2}{(1-\alpha)^{\frac{5}{2}}},\end{align*} where $V_{U}^{\alpha}[f_3](x)$ is the value under the control $U$ with the ``lower-order'' cost function $f_3(x)=x^3$ replacing $x^4+\frac{c_s}{1-u}$.
	We claim that $V_{U}^{\alpha}[f_3](x)\geq \frac{1}{(1-\alpha)^{5/2}}$ for all $x\geq \frac{1}{(1-\alpha)^{5/8}}$, leading to the {\em order-optimality} result
	$$\left|\hV_{U}(x)-V_{U}^{\alpha}(x)\right|\leq \Ex_x^{U}\left[\sum_{t=0}^{\infty} \alpha^t \frac{X_t+1}{1-\alpha}\right]\leq \Gamma V_{U}^{\alpha}[f_3](x),$$ for all such $x$; see Corollary \ref{cor:orderopt}.
	Considering ``large'' initial states is important.
	For small initial values of $x$, the error may be as large as the value itself.
	
	The arguments in this example are not the tightest, but they illustrate the generalizable arguments in \S \ref{subsec:bounds}. \hfill \bsq \end{example}

\noindent Some preliminary construction and definitions will be needed for the statement of our bounds.

\begin{figure}[tb]
	\begin{subfigure}{0.52\textwidth}
		\includegraphics[scale=0.4]{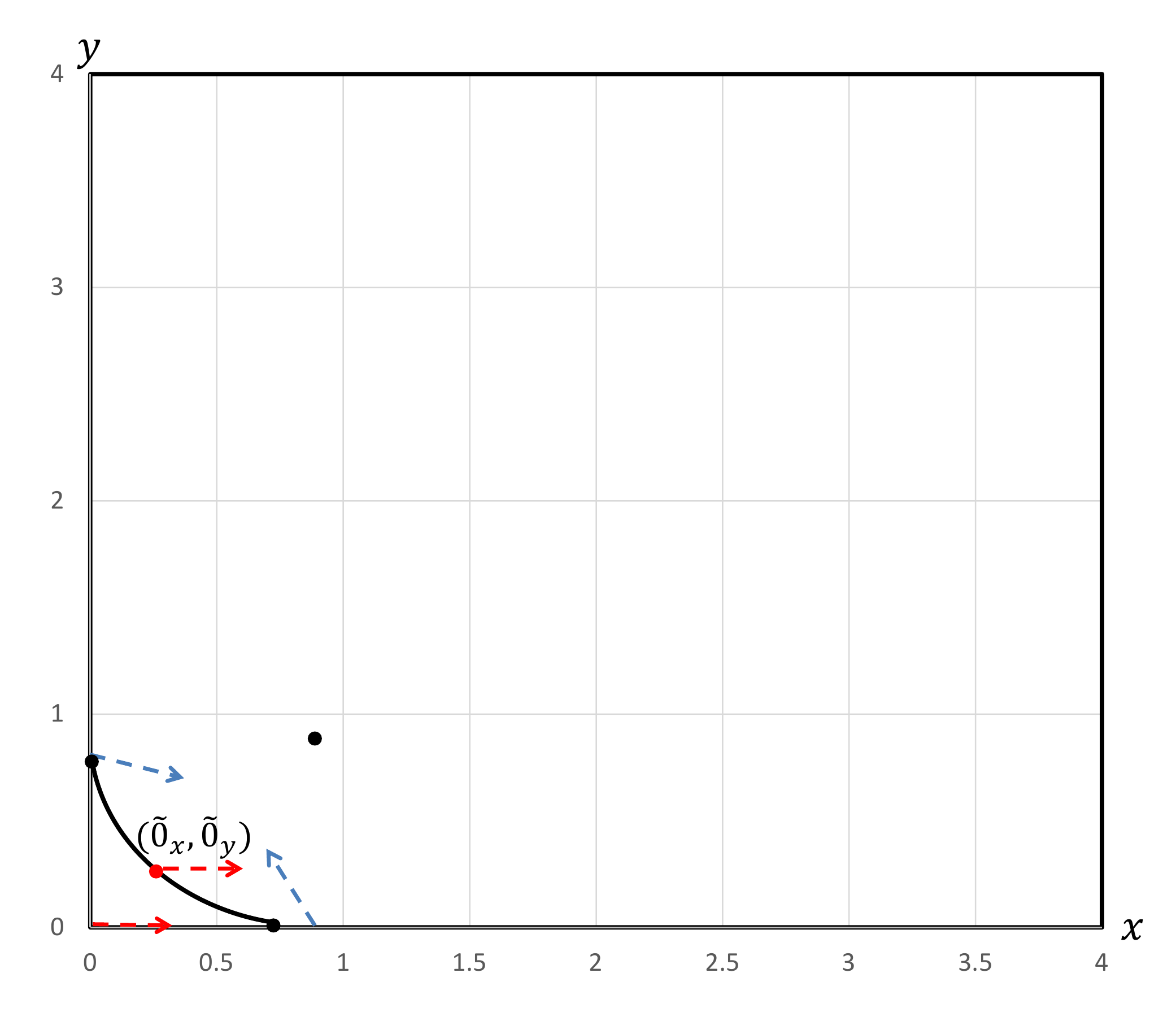}\end{subfigure}
	\begin{subfigure}{0.48\textwidth} \includegraphics[scale=0.39]{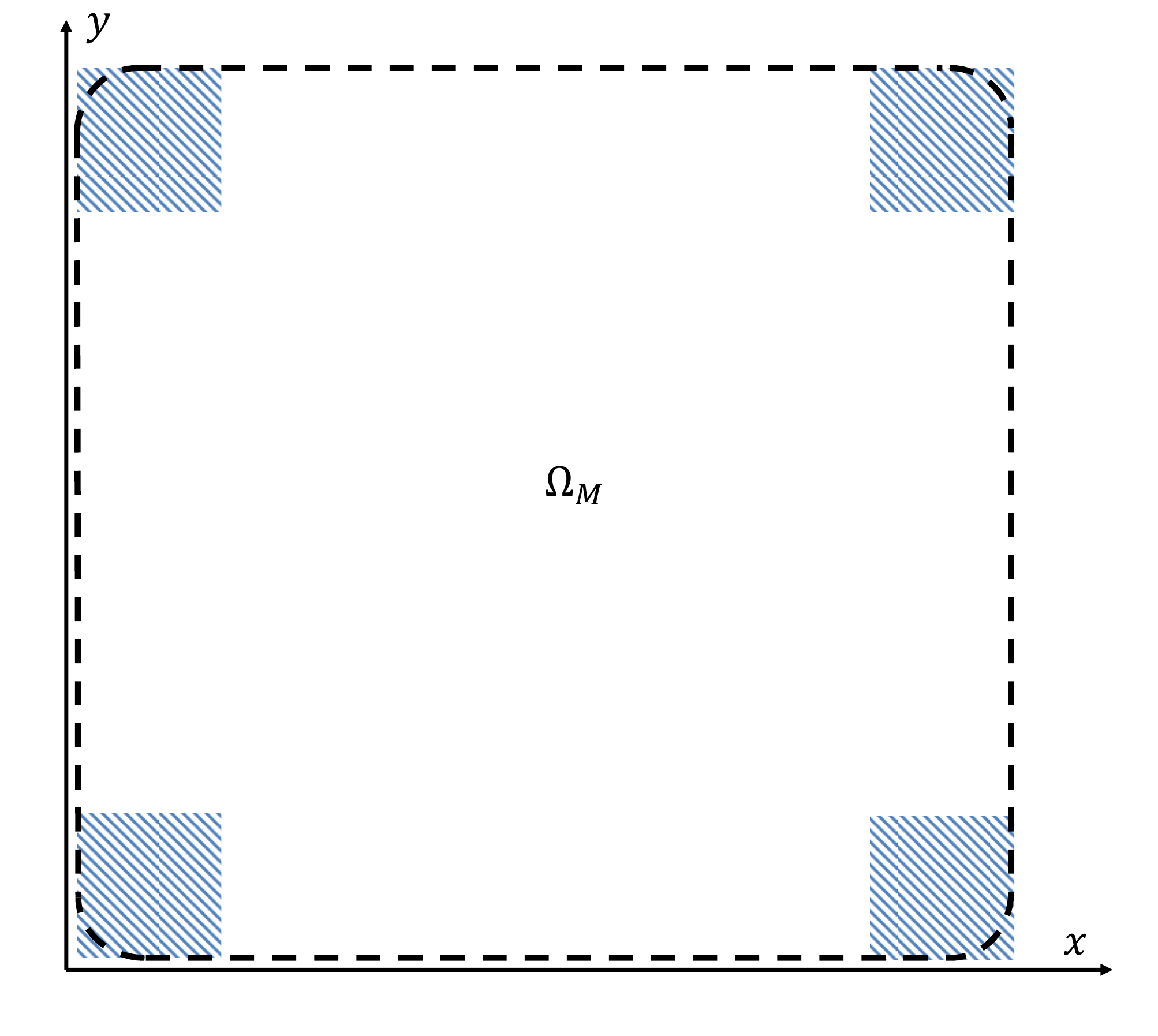} \end{subfigure}
	\caption{A two-dimensional example of truncating the state space and ``curving'' the boundary. The right-hand side displays the smoothed and truncated domain. The squares are the points outside of $\Omega_-^{\ell}$. \label{fig:curved} }	\end{figure}

\paragraph{Smoothing the State Space:} PDEs do not, in general, admit classical solutions in domains with corners; see \cite{dupuis1990oblique}. Fortunately, in our framework we have some freedom in smoothing the domain without compromising the bounds.\footnote{By contrast, in the queueing-approximations literature the state is scaled and smoothing the scaled state space can compromise the optimality gaps.} Consider a two-dimensional controlled chain on the ``square'' state space $\{x\in \mathbb{Z}_+^2: x_1,x_2 \leq M\}$ in Figure \ref{fig:curved}. We can replace the point $0$ with a point $\tilde{0}$---through which we can ``pass'' a smooth boundary while preserving the transition probabilities and the reward function; see Figure \ref{fig:curved}. This {\em does not change} the value function $V_{\ast}^{\alpha}$ or the optimal control $u^*$ of the original chain. It will change the extensions of $r_u(x)$ as well as the values (and extension) of $\mu_u$ and $\sigma^2_u$ but, notice, these change only at states that connect to the moved corner $\tilde{0}$. In this way, the discreteness of the state space supports our analysis; see also \S \ref{sec:conclusions}.

The boundary of the truncated and smoothed state space is illustrated for $d=2$ on the right-hand side of Figure \ref{fig:curved}. We refer to $\Omega_M$ as the open subset of $\mathbb{R}_+^d$ defined by this boundary. This is the domain on which we consider the PDE.

\paragraph{Corners:}  We will distinguish between ``corners,'' such as the point $\tilde{0}$ discussed above, and one-dimensional boundaries. Given a strictly positive function $\varrho(x)$, define the set
$$\Omega_-^{\varrho(\cdot)}:=\{x\in \overline{\Omega_M}: d(x,\mathcal{B}_i)\leq \varrho(x), \mbox{ for at most one $i$}\}.$$

The set $\Omega_-^{\varrho(\cdot)}$ includes states that, while possibly close to an axis, are far from points where two or more axes meet. In the two-dimensional case, taking a constant $\varrho(x)=\ell$ is the same as carving out squares from the corners of $\Omega_M$. Our bounds distinguish between points in $\Omega_-^{\varrho(x)}$ and those points outside of this set. In the one-dimensional case $(d=1)$, $\Omega_-^{\rho(\cdot)}=\overline{\Omega_M}=[0,M]$.

\paragraph{Actions space:} We consider here action sets that do not depend on the state $x$. We can accommodate $\mathcal{U}(x)=\{u:Au\leq b(x)\}\bigcap \mathbb{D}$ by introducing a penalty of the form $\Theta[Au-b(x)]^+$ for suitably large $\Theta$. Notice that the updated reward function that incorporates the penalty---$\bar{r}_u(x)\leftarrow r_u(x)-\Theta[Au-b(x)]^+$---has $$[\bar{r}_u]_{1,\Omega}^*\leq [r_u]_{1,\Omega}^*+\Theta [b]_{1,\Omega}^*.$$

The following is an indirect corollary of the PDE literature. All statements below focus on the case of the OD-boundary condition. Analogues for the FOT boundary condition appear in the appendix; key is the fact that with FOT, \eqref{eq:lipbound} is further restricted to the smaller set $\Omega_{--}^{\varrho(\cdot)}:=\{x\in \Omega_M: d(x,\mathcal{B}_i)> \varrho(x), \mbox{ for all $i$}\}\subset \Omega_-^{\varrho(\cdot)}.$

\begin{lemma}	\label{lem:PDEbounds}
	Let $\Omega_M$ be an open bounded domain as above with boundary in $\mathcal{C}^{2,1}$, and suppose that $\mathcal{U}(x)\equiv \mathcal{U}=\mathbb{D}$ and that Assumption \ref{asum:1} holds. Suppose further that $\eta\in\mathcal{C}^{1,1}(\partial\Omega_M)$ with
	$|\eta|_{\partial \Omega_M}^*+|D\eta|_{\partial \Omega_M}^*+[D\eta]_{1,\partial \Omega_M}^*\leq L$ and that there exists  $\nu_0>0$ such that $\eta(x)\cdot \theta(x)\geq \nu_0|\eta(x)|$ where $\theta(x)$ is the outward normal at $x\in\partial \Omega_M$. Finally, suppose that $\sup_{u\in \mathcal{U}}|r_u|_{\overline{\Omega_M}}^*+\sup_{u\in \mathcal{U}}[r_u]_{1,\overline{\Omega_M}}^*<\infty$.
	
	Then, the TCP with the OD boundary condition has a unique solution $\hV_*\in \mathcal{C}^{2,\beta}(\Omega_M)$ for some $\beta\in (0,1)$ (that does not depend on $\alpha$).  Moreover, given a function $\varrho:\mathbb{Z}_+^d\to \mathbb{R}_{++}$, we have the bound	 	
	\be \label{eq:lipbound}
	[D^2 \hV_*]_{\beta,x\pm \frac{\varrho(x)}{2}}^*\leq \Gamma\left( \frac{|\hV_*|_{\mathcal{B}^{\varrho(x)}(x)}^*}{\varrho^{2+\beta}(x)}+\max_{u \in \mathcal{U}}[r_{u}]_{\beta,\mathcal{B}^{\varrho(x)}(x)}^*\right), ~ ~ x\in  \Omega_{-}^{\varrho(\cdot)}. \ee We also have the global bound $$[D^2 \hV_*]_{\beta,x\pm \frac{\varrho(x)}{2}}^*\leq \Gamma\Theta_M, ~ x\in \overline{\Omega_M},$$ where
$\Gamma\equiv \Gamma(d,\lambda,L,\nu_0,\partial \Omega_M)$ where $\Theta_M:=|\hV|_{\overline{\Omega_M}}^*+\max_{u \in \mathbb{D}}[r_u]_{\beta,\overline{\Omega_M}}^*$.  \end{lemma}

Several comments are due: (i) the bound depends on $\alpha$ only through $|\hV_*|_{\mathcal{B}^{\varrho(x)}(x)}^*$ and possible dependence of the {\em controlled} reward function $r_{\hU_*}$ on $\alpha$; (ii) in the one-dimensional case, $\Omega_{-}^{\varrho(\cdot)}=\overline{\Omega_M}$ so that the bound in \eqref{eq:lipbound} holds in fact in all of $\overline{\Omega_M}$.

For the following, define, $f_k(x)=|x|^k$ for $k\in \mathbb{N}$ and $x\in\mathbb{R}_+^d$ and let $$\mathcal{T}_{CO}=\{t\geq 0: X_t\notin \Omega_{-}^{\varrho(\cdot)}\},$$
which is the set of times the chain spends close to the ``corners'' of $\Omega_M$. This set is, by definition, empty in the one-dimensional case.

\begin{theorem}[explicit bounds] Let $(\hU_*,\hV_*)$ be a solution to the TCP as in Lemma \ref{lem:PDEbounds}. In addition to the requirements in Lemma \ref{lem:PDEbounds}, suppose that $\bar{\frakj}<\infty$ and that for all $x\in\mathbb{R}_+^d$,
\be\max_{u \in \mathcal{U}}[r_{u}]_{\beta,0\pm |x|}^*\leq \Gamma(1+|x|^{k-\beta}), \mbox{ and } |\hV_*|_{0\pm |x|}^*\leq \Gamma\left(\frac{1}{(1-\alpha)^{m}}+\frac{|x|^k}{1-\alpha}\right).\label{eq:requirements2}\ee
Then,
	$$[D^2\hV_*]_{\beta,x\pm \bar{\frakj}}^*\leq \Gamma\left((1+|x|)^{k-\beta}+\frac{1}{(1-\alpha)^\frac{m(k-\beta)}{k+2}}
	+\frac{1}{(1-\alpha)^{\frac{1}{2}(k-\beta)}}\right),$$
for all $x\in\Omega_{-}^{\varrho(\cdot)}$.
In turn, for a stationary policy $U\in\mathbb{U}$ and all $x\in\overline{\Omega_M}$,
	\begin{align*} \Ex_x^{U}\left[\sum_{t\notin \mathcal{T}_{CO}}\alpha^t[D^2\hV_*]_{\beta,X_t\pm \bar{\frakj}}^*\right]\leq&\Gamma\left(V_U^{\alpha}[f_{k-\beta}](x)+\frac{1}{(1-\alpha)^{\frac{m(k-\beta)}{k+2}+1}}+\frac{1}{(1-\alpha)^{\frac{1}{2}(k+2-\beta)}}\right),\end{align*}
	so that
	\begin{align} |V_{\hU_*}^{\alpha}(x)-V_*^{\alpha}(x)|\leq &  \Gamma\left(\max\{V_{U_*}^{\alpha}[f_{k-\beta}](x),V_{\hU_*}^{\alpha}[f_{k-\beta}](x)\}+\frac{1}{(1-\alpha)^{\frac{m(k-\beta)}{k+2}+1}}+ \frac{1}{(1-\alpha)^{\frac{1}{2}(k+2-\beta)}}\right.\nonumber \\ &\qquad + \left. \Theta_M
		\Ex_x^{\hU_*}\left[\sum_{t\in \mathcal{T}_{CO}}\alpha^t\right]+\Theta_M
		\Ex_x^{U_*}\left[\sum_{t\in \mathcal{T}_{CO}}\alpha^t\right]\right),\label{eq:expbound} 	
		\end{align} where $\Theta_M \leq \Gamma \left(\frac{1}{(1-\alpha)^m}+\frac{|M|^k}{1-\alpha}\right). $
	 \label{thm:explicit}
\end{theorem} Above the term on the second line is the contribution of the corners. The magnitude of the boundary effect depends on how much time the chain spends near ``corners.''  We revisit this point in \S \ref{sec:HT}.

\begin{remark}[a priori requirements on $\hV_*$ and the value of $m$] We can gain some insight into the value of $m$ in the requirement \eqref{eq:requirements2} by considering the one-dimensional case. Suppose that $(\hU_*,\hV_*)$ satisfy	$$0=r_{\hU_*}(x)+\alpha\mathcal{L}_{\hU_*}\hV_*(x)-(1-\alpha)\hV_*(x),~x>0,$$ and $\hV_*'(0)=0$. Suppose, further, that (i) ~ $\mu_{\hU_*}(x)\leq 0$ for all $x>0$ and that, (ii) uniformly in $x$, $0<\underline{\sigma}\leq  \sigma_{\hU_*}(x)\leq \bar{\sigma}<\infty$. It follows from basic arguments (see Lemma \ref{lem:1dimbound} in the appendix) that, for $k>1$,
	$$\hV_*(x)\leq \Gamma\left(\frac{x^k}{1-\alpha}+\frac{1}{(1-\alpha)^{\frac{k+2}{2}}}\right),$$ so that we can take $m=\frac{k+2}{2}$.
	
	If there exists $\kappa>0$ such that $\mu_{\hU_*}(x)\leq -\kappa$ for all $x>0$, then $m=1$, i.e., $\hV_*(x)\leq \Gamma\left(\frac{x^k}{1-\alpha}+\frac{1}{1-\alpha}\right)$. If  $\mu_{\hU_*}(x)$ is not necessarily negative but it is bounded--- $|\mu_{\hU_*}(x)|\leq \kappa$ (not dependent on $\alpha$ or $x$)---then we can take $m=k+1$; see the proof of Lemma \ref{lem:1dimbound} for further comments.
\hfill \bsq \end{remark}

If $r_u(x)\geq \gamma f_k(x)$ for some $k$ and all $u$, then it trivially holds that $|V_{U}^{\alpha}(x)|\geq \gamma V^{\alpha}_{U}[f_k](x)$. When this ``super-polynomial'' property holds we have the following corollary.
\begin{corollary}[vanishing discount optimality] Let $(\hU_*,\hV_*)$ be a solution to the TCP and suppose the assumptions of Theorem \ref{thm:explicit} hold with  $m\leq k+1$. Then, for every stationary policy $U$,
	\begin{align*} \Ex_x^{U}\left[\sum_{t\notin \tau_{OC}}\alpha^t[D^2\hV_*]_{\beta,X_t\pm \bar{\frakj}}^*\right]\leq  \Gamma (1-\alpha)^{\beta}V_{U}^{\alpha}[f_k](x),	
	\end{align*}
for all $x:|x|\geq \frac{1}{1-\alpha}$.
Consequently, if
\begin{equation}
|V_{\ast}^{\alpha}(x)|=|V_{U_*}^{\alpha}(x)|\geq \gamma V_{U_*}^{\alpha}[f_k](x)\mbox{ and  } |V_{\hU_*}^{\alpha}(x)|\geq \gamma V_{\hU_*}^{\alpha}[f_k](x),~x\in\mathbb{Z}_+^d,
\end{equation}
for some $\gamma>0$ that does not depend on $\alpha$, then
\begin{align*} |V_{\hU_*}^{\alpha}(x)-V_*^{\alpha}(x)|\leq & \Gamma (1-\alpha)^{\beta}\max\{|V_{\ast}^{\alpha}(x)|,|V_{\hU_*}^{\alpha}(x)|\}\\ &+\Gamma\Theta_M \left(\Ex_x^{\hU_*}\left[\sum_{t\in \mathcal{T}_{CO}}\alpha^t\right]+
	\Ex_x^{U_*}\left[\sum_{t\in \mathcal{T}_{CO}}\alpha^t\right]\right).\end{align*} 	\label{cor:alphaopt}	
\end{corollary}

Of course, if $M< \frac{1}{1-\alpha}$ the statement above is moot, because then $|x|\leq 1/(1-\alpha)$ for all $x\in \overline{\Omega_M}$.

\begin{corollary}[order optimality] Suppose the assumptions of Theorem \ref{thm:explicit} hold and let $\zeta(m,k)=\frac{1}{k+1-\beta}\max\{\frac{m(k-\beta)}{k+2}+1,\frac{1}{2}(k+2-\beta)\}$. Then, for all $x$ such that $|x|\geq \frac{1}{(1-\alpha)^{\zeta(m,k)}}$,
\begin{align*} |V_{\hU_*}^{\alpha}(x)-V_*^{\alpha}(x)|\leq & \Gamma\max\{V_{U_*}^{\alpha}[f_{k-\beta}](x),V_{\hU_*}^{\alpha}[f_{k-\beta}](x)\}\\
	&+\Gamma\Theta_M \left(\Ex_x^{\hU_*}\left[\sum_{t\in \mathcal{T}_{CO}}\alpha^t\right]+
	\Ex_x^{U_*}\left[\sum_{t\in \mathcal{T}_{CO}}\alpha^t\right]\right).\end{align*} 			
	\label{cor:orderopt}
\end{corollary}

In the case that $D^2\hV$ is Lipschitz continuous ($\beta=1$), the first term of the bound in 	Corollary \ref{cor:alphaopt} is proportional to a  $(1-\alpha)$ fraction of the value while that in Corollary \ref{cor:orderopt} corresponds to order optimality with an error proportional to the value with $f_{k-1}$.

\begin{remark}[Example \ref{example:SSQ} revisited] In Example \ref{example:SSQ}, where we were considering a fixed control, we have $m=3=k-1$ and $\beta=1$, so that $\frac{1}{2}(k+2-\beta)=\frac{m(k-\beta)}{k+2}+1=2.5$ and the terms depending on $(1-\alpha)$ in the first line of \eqref{eq:expbound} correspond to $1/(1-\alpha)^{2.5}$, as we obtained in Example \ref{example:SSQ} through explicit derivations. Moving from control to optimization we have--- because $r(x,u)\leq -x^4$---that the requirements of both corollaries \ref{cor:alphaopt} and \ref{cor:orderopt} are satisfied and their conclusions apply. Also, there are no ``corners'' in this one-dimensional example, so the second line of \eqref{eq:expbound} is $0$.
 \label{rem:ssqbound} \hfill \bsq
\end{remark}

				\section{Computation: Tayloring-based approximate dynamic programming}
				
				The two staple algorithms for solving DPs are the value and policy iteration algorithms. The curse of dimensionality renders both incapable of solving large-scale problems and motivates the development of approximation algorithms. We present here an algorithm that we call Taylored Approximate Policy Iteration (TAPI), which offers a reduction in computational effort compared to solving the original DP; in \S \ref{sec:exmultipool}, for example, we solve a problem where TAPI takes less than 10 minutes for some instances where the exact solution takes more than 15 hours.  The algorithm comes with verifiable conditions for convergence, and optimality-gap bounds that are grounded in our analysis of Tayloring. 				
				
				Whereas the state space of the DP is unbounded, truncating is inevitable for computation. We use $\mathbb{S}$ for the truncated state space, $\widetilde{\mathbb{S}}$ for its continuation, and $\partial\widetilde{\mathbb{S}}$ for the boundary of $\widetilde{\mathbb{S}}$. Mostly, we restrict attention to the case in which the original state space is $\Z^d_+$ and introduce the truncated state space $\mathbb{S} = [0,M]^d \cap \Z^d_+$ and its continuation $\widetilde{\mathbb{S}} = [0,M]^{d}$.

						Because $\mathbb{S}$ is finite, the Bellman equation
				\begin{align*}
					V_{\ast}(x) = \max_{u \in \mathcal{U}(x)} \Big\{ r_{u}(x) + \alpha \E_{x}^u \big[V_{*}(X_1) \big] \Big\}, \quad x \in \mathbb{S},
				\end{align*} has a unique solution and the PI algorithm is guaranteed to converge to this solution in finitely many iterations.
				\begin{algorithm}[hh!]
					\caption{(Standard) Policy Iteration (PI)} \label{alg:pi}
					\begin{enumerate}
						\item Start with some initial stationary policy $U^{(0)}\in\mathbb{U}$.
						\item For $k = 0,1,\ldots$
						\begin{enumerate}
							\item  {\em Policy evaluation:} solve for the infinite horizon discounted performance under $U^{(k)}$, i.e., find $V^{(k)}(x)$ that satisfies
							\begin{align}
								V^{(k)}(x) = r(x,U^{(k)}(x)) + \alpha \E_{x}^{U^{(k)}(x)} \big[V^{(k)}(X_1)\big], \quad x \in \bbS. \label{eq:PIperformance}
							\end{align}
							\item {\em Policy improvement:} Find
							\begin{align}
								U^{(k+1)}(x) = \arg \max_{u \in \mathcal{U}(x)} \Big\{ r(x,u) + \alpha \E_x^u [V^{(k)}(X_1)\big] \Big\}, \quad x \in \bbS. \label{eq:greedyexact}
							\end{align}
						\end{enumerate}
					\end{enumerate}
				\end{algorithm}
				
				The computational bottlenecks of PI are well understood:
				
				\textbf{Value-function storage:} we require $\calO(|\mathbb{S}|)$ space to store $V^{(k)}(x)$, and $|\mathbb{S}|$ may grow exponentially with the dimension of the problem.
				
				\textbf{Transition-matrix storage:} In the $k^{th}$ step of the PI algorithm, we invert $P^{(k)} - I$, where $P^{(k)}$ is the transition probability matrix associated with policy $U^{(k)}(x)$. Depending on density (or sparsity) of this matrix, storing it may require as much as  $\calO(|\mathbb{S}|^2)$. The value-iteration algorithm does not require such storage.
				
				\textbf{Optimization complexity:} Each iteration includes a policy improvement step: finding (greedy) optimal actions relative to the value-function approximation $V^{(k)}(x)$. For each state $x \in \mathbb{S}$ and action $u\in\mathcal{U}(x)$, computing the expectation $\E_x^u [V^{(k)}(X_1)\big]$ may require as many as $\calO(|\mathbb{S}|)$ function evaluations, depending on the transition structure of the Markov chain. Furthermore, the optimization may require exhaustive search over the discrete action space $\mathcal{U}(x)$. The total cost of the policy improvement step can therefore be up to $\calO( | \mathcal{U}|_{max}|\mathbb{S}|^2)$, where $  | \mathcal{U}|_{max}= \sup_{x \in\mathbb{S}}| \mathcal{U}(x)|$ is an upper bound on the number of feasible actions. Our example in \S \ref{sec:exmultipool} is one where this worst-case cost is realized.


				\subsection{Taylored Approximate Policy Iteration (TAPI)} \label{sec:TAPI}
				
				The basic idea in approximate policy iteration (API) is to produce an approximation of $V^{(k)}(x)$ for the value function at iteration $k$, and then use it in the policy improvement step. Linear architecture is a popular approximation scheme that uses an element of the space $\{\Phi r | r \in \R^{F}\}$, where  $\Phi$ is an $|\mathbb{S}| \times F$ feature matrix. The columns of $\Phi$ are so-called feature vectors that capture pre-selected properties of each state $x \in \mathbb{S}$. The features of a particular state can be generated on the fly, so that there is no need to store the entire matrix $\Phi$; it suffices to store $r$ to represent all elements of $\{\Phi r | r \in \R^{F}\}$. This produces computational benefits when $F$ is significantly smaller than $|\mathbb{S}|$. The optimality gaps of API depend on the ``richness" of the feature vectors, which are typically chosen based on structural insight into the problem at hand; see \cite{Bert2011} for a survey of API methods.
				
				Tayloring offers a generalizable way to approximate $V^{(k)}(x)$ that requires little ad-hoc intuition. In this scheme, the intermediate solution $V^{(k)}(x)$ is replaced in the policy evaluation step by the solution to the associated partial differential equation. In addition to approximating the policy evaluation step, we can also approximate the policy improvement step. The details of TAPI are presented in Algorithm~\ref{alg:tcppi}.
				
				\begin{algorithm}
					\caption{Taylored Approximation Policy Iteration (TAPI)  \label{alg:tcppi}}
					\begin{enumerate}
						\item Start with initial stationary policy $U^{(0)}(x)$.
						\item For $k = 0,1,\ldots$
						\begin{enumerate}
							\item {\em Approximate policy evaluation:} approximate, using the Taylored equation, the infinite horizon discounted performance under $U^{(k)}$, i.e., find $V^{(k)}(x)$ that satisfies
							\begin{align}
								& r(x,U^{(k)}(x)) +  \alpha \mathcal{L}_{U^{(k)}(x)}V^{(k)}(x) - (1-\alpha)V^{(k)}(x) = 0, \quad x \in \widetilde{\mathbb{S}},\label{eq:tapi}\\\nonumber
								& \eta(x)'DV^{(k)}(x)=0, \quad x\in \partial \widetilde{\mathbb{S}}.
							\end{align}
							\item {\em Approximate policy improvement:} Let $U^{(k+1)}(x)$ be the greedy policy associated with $V^{(k)}(x)$ in the Taylored equation:
							\begin{align}
								U^{(k+1)}(x) = \arg \max_{u \in \mathcal{U}(x)} \Big\{ r(x,u) +  \alpha \mathcal{L}_{u}V^{(k)}(x)- (1-\alpha)V^{(k)}(x) \Big\}, \quad x \in \widetilde{\mathbb{S}}. \label{eq:tcpgreedy}
							\end{align}
						\end{enumerate}
					\end{enumerate}
				\end{algorithm}

				In the $k^{th}$ step, assuming the PDE \eqref{eq:tapi} has a solution, this solution can be numerically approximated by any of a number of solution methods for {\em linear} PDEs. The most standard of these is the finite difference (FD) method; see, e.g., \cite{larsson2008partial}. FD returns a solution defined on a suitably spaced grid. In our experiments this grid is a subset of the state space $\bbS$. The efficiency gains of TAPI cover all three of the previously identified computation bottlenecks of PI:
				
				{\bf Value-function storage:} Any method to solve \eqref{eq:tapi} involves either a discretized grid or some other state space partitioning scheme (as in the finite element method). As the discretization gets finer, the approximation converges, under suitable conditions, to the true solution of the PDE. Choosing a {\em coarser} grid reduces the cost of storing the value function estimates.
				
				{\bf Transition-matrix storage:} In \eqref{eq:tapi} the transition probability matrix ``collapses'' into the lower dimensional $\mu_u(x)$ and $\sigma^2_u(x)$. In contrast to the standard PI algorithm, we are not inverting the full matrix $P^{(k)}$ here.
				
				\textbf{Optimization complexity:} The computational benefit of the approximate policy improvement in TAPI comes from the fact that $\mathcal{L}_{u}V^{(k)}(x)$ depends on $u$ only through $\mu_u(x)$ and $\sigma^2_u(x)$, while $V^{(k)}(x)$ and its derivatives do not depend on $u$. For finite action and state spaces the quantities $\mu_u(x)$ and $\sigma^2_u(x)$ can be pre-computed once in advance (or computed on the fly and kept in memory). Contrast this with PI, where the term $\E_{x}^u[V^{(k)}(X_1)]$ has to be re-computed for each $u$ and $x$ at {\em each} iteration $k$ and computing the expectation requires going over all the ``neighbors'' of $x$. The computational cost of the approximate policy improvement is, consequently, $\calO( | \mathcal{U}|_{max}|\mathbb{S}|)$ per iteration compared to $\calO( | \mathcal{U}|_{max}|\mathbb{S}|^2)$ for exact policy improvement. This cost may be further reduced if, given $x$, one has tractable expressions for the dependence of $\mu_u(x)$ and $\sigma^2_u(x)$ on $u$; see the examples in the next section. 			
				
				Given the discrete nature of the controls, the exact policy improvement step in \eqref{eq:greedyexact} can be computationally expensive. Our example in \S \ref{sec:exmultipool} is one where it is difficult to avoid exhaustive search. In
				\cite{MoalKumaVanr2008}, the authors show how to leverage an ``affine-expectations'' assumption to approximate the solution of this problem by that of a linear program.  The approximate policy-improvement step in TAPI is a generalizable way to simplify this step.

				
				\subsubsection{Convergence and Error Bounds}
				As stated, the PDE in \eqref{eq:tapi} might not be mathematically meaningful; $U^{(k)}(x)$ could be such that a solution does not exist to the PDE in the policy evaluation step. One implementation we propose---developed for the specific PDEs arising from diffusion control problems---is put forth in \cite{kushner2013numerical}.  Roughly speaking, applying certain finite difference schemes to \eqref{eq:tapi} leads back to discrete (time and space) MDPs.  The goal in \cite{kushner2013numerical} is to solve the continuous control problem by taking the discretization to be increasingly finer. We take the ``opposite approach,'' and choose a coarse grid to reduce the computational effort relative to the original Bellman equation. The construction of \cite{kushner2013numerical} generates one concrete TCP-equivalent chain. The approximation error that it introduces is related to the third derivative of the PDE multiplied by a number that captures the coarseness of the grid; recall the discussion closing \S \ref{sec:initial} and see \eqref{eq:KDbound}.
				
				TAPI then reduces to PI for a (newly constructed) chain on a finite state space. The convergence to the optimal policy then follows from standard results for policy iteration.
				
				The following one-dimensional example illustrates the K-D construction.


				\begin{example}[K-D construction in one-dimension] Consider the one-dimensional TCP on the truncated state space $[0,M]$:
					\begin{align}
						0 = \max_{u \in \mathcal{U}(x)}\Big\{r(x,u) +\alpha \Big(\mu_{u}(x)V'(x)+\frac{1}{2}\sigma^2_{u}(x)V''(x)\Big) - (1-\alpha)V(x)\Big\}, \quad x \in (0,M), \label{eq:tcponedim}
					\end{align} with the boundary condition $$V'(0)=V'(M)=0.$$
					Fix $h>0$ and let $\mathbb{S}_h=\{0,h,2h,\ldots,M\}$ be the discretized space. Let us make the assumption that $M$ is divisible by $h$.
					
					The K-D chain construction is most intuitive under the ``small-drift'' assumption
					\begin{align}
						\sigma^2_{u}(x)\geq |\mu_{u}(x)|h, \quad x \in  \mathbb{S}_h, u \in \mathcal{U}(x). \label{eq:smdrift}
					\end{align}
					
					For each $x\in \mathbb{S}_h\backslash\{0,M\}$, let us replace $V'(x)$ and $V''(x)$ with the appropriate ``central'' differences $$V'(x)\leftarrow \frac{V(x+h)-V(x-h)}{2h}, \mbox{ and } V''(x)\leftarrow \frac{V(x+h)-2V(x)+V(x-h)}{h^2},$$
					to get
					\begin{align}
						&(1-\alpha)V(x) \notag \\
						=&\max_{u\in\mathcal{U}(x)}\left\{r(x,u) +\alpha \Big(\mu_{u}(x)\frac{V(x+h)- V(x-h)}{2h} +\frac{1}{2}\sigma^2_{u}(x)\frac{V(x+h)-2V(x)+V(x-h)}{h^2}\Big)\right\} \nonumber \\\label{eq:onedimKD}
						=&\max_{u\in\mathcal{U}(x)}\left\{r(x,u) +\alpha\Big(\frac{\mu_{u}(x)h+\sigma_{u}^2(x)}{2h^2}V(x+h) +\frac{-\mu_{u}(x)h+\sigma^2_{u}(x)}{2h^2}V(x-h) - \frac{\sigma^2_{u}(x)}{h^2} V(x)\Big)\right\}
					\end{align}
					Let $\Sigma(x) = \sup_{u \in \mathcal{U}(x)} \sigma^2_{u}(x) > 0$. Multiplying both sides of \eqref{eq:onedimKD} by $h^2/\alpha \Sigma(x)$, we arrive at
					\begin{align}\notag
						V(x) =&\max_{u\in\mathcal{U}(x)}\left\{ \frac{\alpha_h(x) h^2 r(x,u)}{\alpha \Sigma(x)}  + \alpha_h(x)\left(\frac{\mu_{u}(x)h+\sigma_{u}^2(x)}{2\Sigma(x)}V(x+h) +\frac{-\mu_{u}(x)h+\sigma^2_{u}(x)}{2\Sigma(x)}V(x-h)+\right. \right.\\ & \qquad\left. \left.\Big(1 - \frac{\sigma^2_{u}(x)}{\Sigma(x)}\Big) V(x)\right)\right\}, \label{eq:interior}
					\end{align}
					where $$\alpha_h(x):=\left(1 + \frac{h^2}{\Sigma(x)}\left(\frac{1}{\alpha} - 1 \right)\right)^{-1}.$$

					Let, for $x\in \mathbb{S}_h\backslash\{0,M\}$,
					\begin{align} P^{u,h}_{x,x+h}=	\frac{\mu_{u}(x)h+\sigma_{u}^2(x)}{2\Sigma(x)},
						~  P^{u,h}_{x,x-h}=
						\frac{-\mu_{u}(x)h+\sigma_{u}^2(x)}{2\Sigma(x)},~ P^{u,h}_{x,x}=\ 1- P^{u,h}_{x,x+h}- P^{u,h}_{x,x-h},\end{align}
					and $\tilde{r}_h(x,u)=\alpha_h(x) h^2 r(x,u)/(\alpha \Sigma(x))$. Notice that these are well-defined probabilities because of assumption \eqref{eq:smdrift}. Also notice that $\tilde{r}_h(x,u)=\frac{1-\alpha_h(x)}{1-\alpha}r(x, u)$. We arrive at the equation
					\begin{align}\notag
						V(x) =&\max_{u\in\mathcal{U}(x)}\left\{ \tilde{r}_h(x,u)  + \alpha_h(x) \left(P_{x,x+h}^{u,h}V(x+h) +P_{x,x-h}^{u,h}V(x-h)+P_{x,x}^{u,h} V(x)\right)\right\}.
					\end{align}
					
					This equation, in the interior, is recognizable as a Bellman equation for a new Markov chain with state space $\mathbb{S}_h$, transition probabilities and reward function as specified, and {\em state dependent} discount factor $\alpha_h(x)$.
					
					For the boundary points $x=0$ and $x=M$ we cannot use central differences, because the points $-h$ and $M+h$ are not available. We can use instead the forward difference $V'(0)\leftarrow \frac{V(h)-V(0)}{h},$ at $x=0$ and the backward difference $V'(M)\leftarrow \frac{V(M)-V(M-h)}{h},$ which leads then to the added equation
					$V(h)=V(0)$ and $V(M-h)=V(M)$. No discount factor is associated with these boundary states. A thorough treatment of reflecting boundaries appears in \cite[Chapter 5.7]{kushner2013numerical}. This construction can be easily modified to have first-order Tayloring on the boundary instead of the oblique-derivative condition.
					\hfill \bsq \label{example:KDonedim}
				\end{example}
				
				The K-D chain is but one concrete construction of a TCP-equivalent chain. A nice property of this construction is the ``sparsity of neighbors''---that from each state one can only transition to at most $2^d$ neighbors. In fact, it follows from \cite{kushner2013numerical} that in the setting of the oblique-derivative boundary condition there {\em always} exists a TCP-equivalent construction on the coarser grid. While this construction need not be as simple as in the one-dimensional case above, it always maintains the desirable properties of sparsity of ``neighbors.''
				
				In the multi-dimensional case, the state space of the K-D chain is $$\mathbb{S}_h= \times_{i=1}^d\Big\{ [0,M_i]\cap \{h\bbZ_+\cup \{M_i\}\}\Big\}$$ where $h\bbZ_+=\{0,h,2h,3h,\ldots\}$. We denote by $X^h=\{X_t^h,~t=1,2,\ldots\}$ the (controlled) Markov chain on the state space $\mathbb{S}_h$  arising from the K-D construction. As in Example \ref{example:KDonedim}, we let $V_*^{h}(x)$ be the solution to the Bellman equation for the K-D chain and denote by $U_*^h$ the optimal stationary policy for this chain.
				
				Then, under the requirements on $\hV_{\ast}$ in Theorem \ref{thm:diffusion}, equation \eqref{eq:TCPequivBound} implies the bound
				\be  |\hV_{\ast}(x)-V_{\ast}^{h}(x)| \leq h^{2+\beta}\left(\Ex_x^{\hU_*}\left[\sum_{t=0}^{\infty}\bar{\alpha}_h^t [D^2\hV_{\ast}]_{\beta,X_t^h\pm h}^*\right]+\Ex_x^{U_*^h}\left[\sum_{t=0}^{\infty} \bar{\alpha}_h^t [D^2\hV_{\ast}]_{\beta,X_t^h\pm h}^*\right]\right),~x \in \mathbb{S}_h.\label{eq:KDbound}\ee
				
				This bound is similar in spirit to and inspired by \cite{dupuis1998rates}. A challenge here is that $V_{\ast}^{h}(x)$ and $U_*^h$ are only defined on the coarse grid $\mathbb{S}_h$. To borrow a term from the ADP literature, we must now ``disaggregate'' these to generate a policy for the original chain.
				
			Let us assume that we have an extension of $V_{\ast}^h(x)$, denoted by $\widetilde{V}_{\ast}^h(x)$ that is defined for all $x\in \mathbb{S}$. Let $U^h$ be the policy obtained (for the original chain) by one-step improvement starting off the value function $\widetilde{V}_{\ast}^h(x)$, i.e., $U^h$ is the greedy policy relative to $\widetilde{V}_{\ast}^h$, i.e.,
$$U^h(x)=\argmax_{u\in \mathcal{U}(x)}\{ r_{u}(x)+\alpha\Ex^u[\widetilde{V}_{\ast}^h(X_1)]\}.$$
 Then we have---see \cite[Proposition 1.3.7]{bertsekas2007approximate}---that
$$\sup_{x\in\mathbb{S}}|V_{U^h}(x)-V_{\ast}(x)|\leq
\frac{\max_{x \in\mathbb{S}}D(x)}{1-\alpha}.$$
where
\begin{align*}
D(x):=&\ |V_*(x)-\widetilde{V}_{\ast}^{h}(x)| \leq |V_*(x)-\hV_{\ast}(x)|+|\widetilde{V}_*^{h}(x)-\hV_{\ast}(x)|.
\end{align*}

Notice that the first term on the right can be bounded using Theorems \ref{thm:diffusion} and \ref{thm:explicit}. The second term depends on the way in which we extend $V_*^h$ to $\widetilde{V}_*^h$. We can, for example, use a piecewise constant extension where we let $\tilde{V}_*^h(x) =V_*^h(A(x))$, with $A: \widetilde{\mathbb{S}} \to \mathbb{S}_h$ being, for example, the nearest neighbor map $A(x) = \inf_{y \in \mathbb{S}_h} |x-y|$. Then
$$|\widetilde{V}_*^{h}(x)-\hV_{\ast}(x)|=
|V_*^{h}(A(x))-\hV_{\ast}(x)|\leq |V_*^{h}(A(x))-\hV_{\ast}(A(x))|+|\hV_{\ast}(A(x))-\hV_{*}(x)|,$$ where the first term can be bounded using \eqref{eq:KDbound} and the second is bounded by $h|D\hV_{*}|_{x\pm h}$. If, instead of the piecewise constant extension, we define, for each $x$, $\widetilde{V}^h(x)$ to be a quadratic interpolation of suitably chosen neighbors, then the error will in fact be proportional to $|D^3\hV_{*}|_{x\pm 2h}$.

				\begin{remark}[exact policy improvement]\label{rem:exactimprovement} Algorithmically, one can replace the approximate policy improvement step in \eqref{eq:tcpgreedy} with an {\em exact} policy improvement step where we let $U^{(k+1)}(x)$ be the greedy policy associated with $V^{(k)}(x)$, i.e.
					\begin{align*}
						U^{(k+1)}(x) = \arg \max_{u \in \mathcal{U}(x)} \Big\{r(x,u) + \alpha \E_x^u[V^{(k)}(X_1)]\Big\}, \quad x \in \mathbb{S}.
					\end{align*}
					
					Because the policy improvement is done exactly---using the transitions and state space of the Markov chain---we must extend $V^{(k)}(x)$ to the state space $\bbS$ (say, by interpolation). In our examples, we find that while this has no convergence guarantees, it may result in a smaller optimality gap. \hfill \bsq \end{remark}
				
				\begin{remark}[relationship with aggregation algorithms] \label{rem:aggr}
The construction of a TCP-equivalent chain on a coarser grid can be viewed as an aggregation procedure. Bounds on the optimality gaps introduced by aggregation methods are often stated in terms of oscillations of the true value function over the coarser grid; see, e.g., \cite[Section 6.7]{bertsekasneuro}. Thus, the bound depends on the same quantity we are trying to avoid computing. 					
					
					Our construction of the approximate coarse chain is grounded in the TCP, which also provides a grounding for optimality-gap analysis. In the cases where our bounds apply, they are stated in terms of the approximation $\hV_*$ rather than by the value function itself.	\hfill \bsq  \end{remark}

				\section{Examples} \label{sec:examples}
				
				The three examples we study are intended to illustrate the performance of the proposed algorithm. The first two examples have a one-dimensional state space and are hence computationally cheap, even for an exact solution. We use them because the one-dimensional case simplifies visualization and helps underscore some observations.
				
				\subsection{Service-rate control}

				This example is a variant of Example \ref{example:SSQ}. Suppose a holding $c(x)=x^2$ when there are $x$ customers in the system and a control cost $f(u)=1/(1-u)$. The cost minimization problem is equivalent to a reward maximization problem with the negative reward $-x^2-1/(1-u)$. The control set consists of the rational numbers (denoted by $\mathbb{Q}$) in $[0,1]$. The Taylored equation
				\begin{align*} 0&=\min_{u\in [0,1]\bigcap\mathbb{Q}} \left\{x^2+\frac{1}{1-u}+\alpha\left((1-2u)V'(x)+\frac{1}{2}V''(x)\right)-(1-\alpha)V(x)\right\},\\0&= V'(0).\end{align*}
				
				Per Lemma \ref{lem:PDEbounds}, this equation has a solution $\hV_*\in \mathcal{C}^{2,\beta}$ for some $\beta\in (0,1)$ and it satisfies the bounds in Theorem \ref{thm:explicit}.
				
				We use TAPI based on the K-D chain to obtain the optimal control $\hU_*^h$ for this chain. We build a control for the original chain by extending to $\mathbb{Z}_+$ in a piecewise constant manner: the control at point $mh$ is kept constant for all points $mh,\ldots,(m+1)h-1$. We denote this control by $U^h$. We try $h\in \{1,2\}$. In our computations we allow the control to be any number in $[0,1]$ which allows, in this example, to write the control as an explicit function of the value in neighboring states.
				
				\begin{figure}[hh!]
					\begin{center}
						
						\includegraphics[scale=0.29]{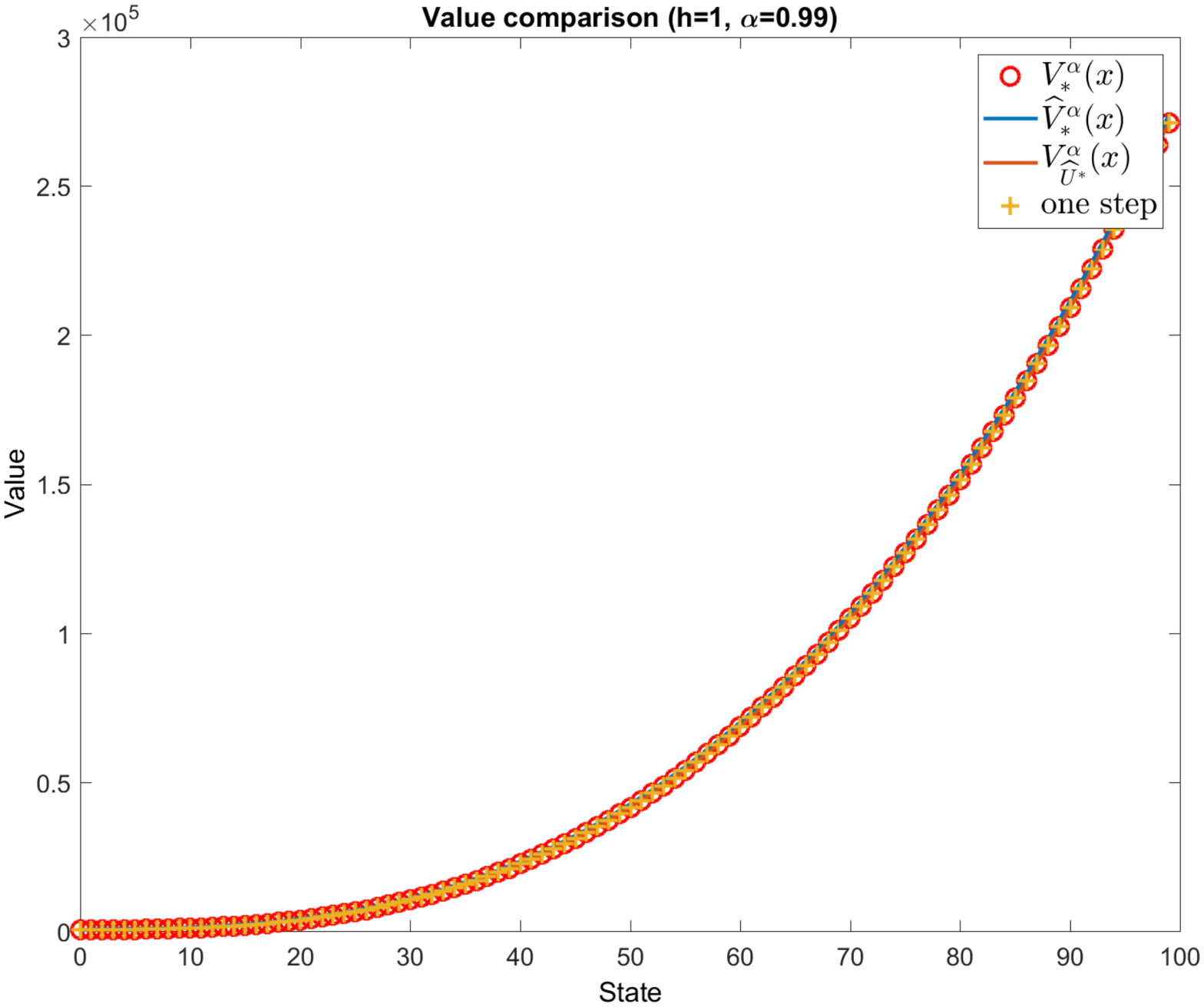}		\includegraphics[scale=0.29]{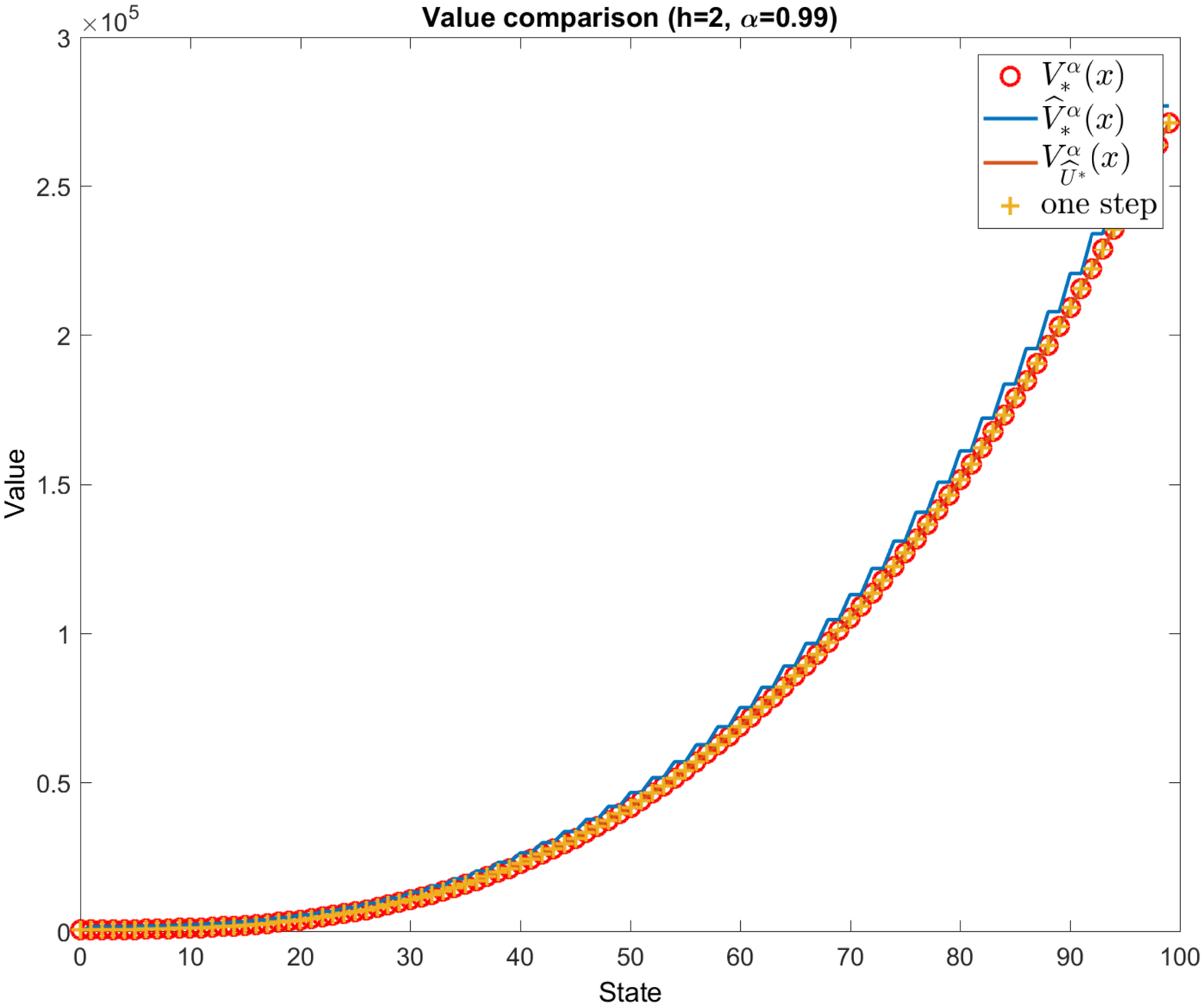}	\\
						\includegraphics[scale=0.29]{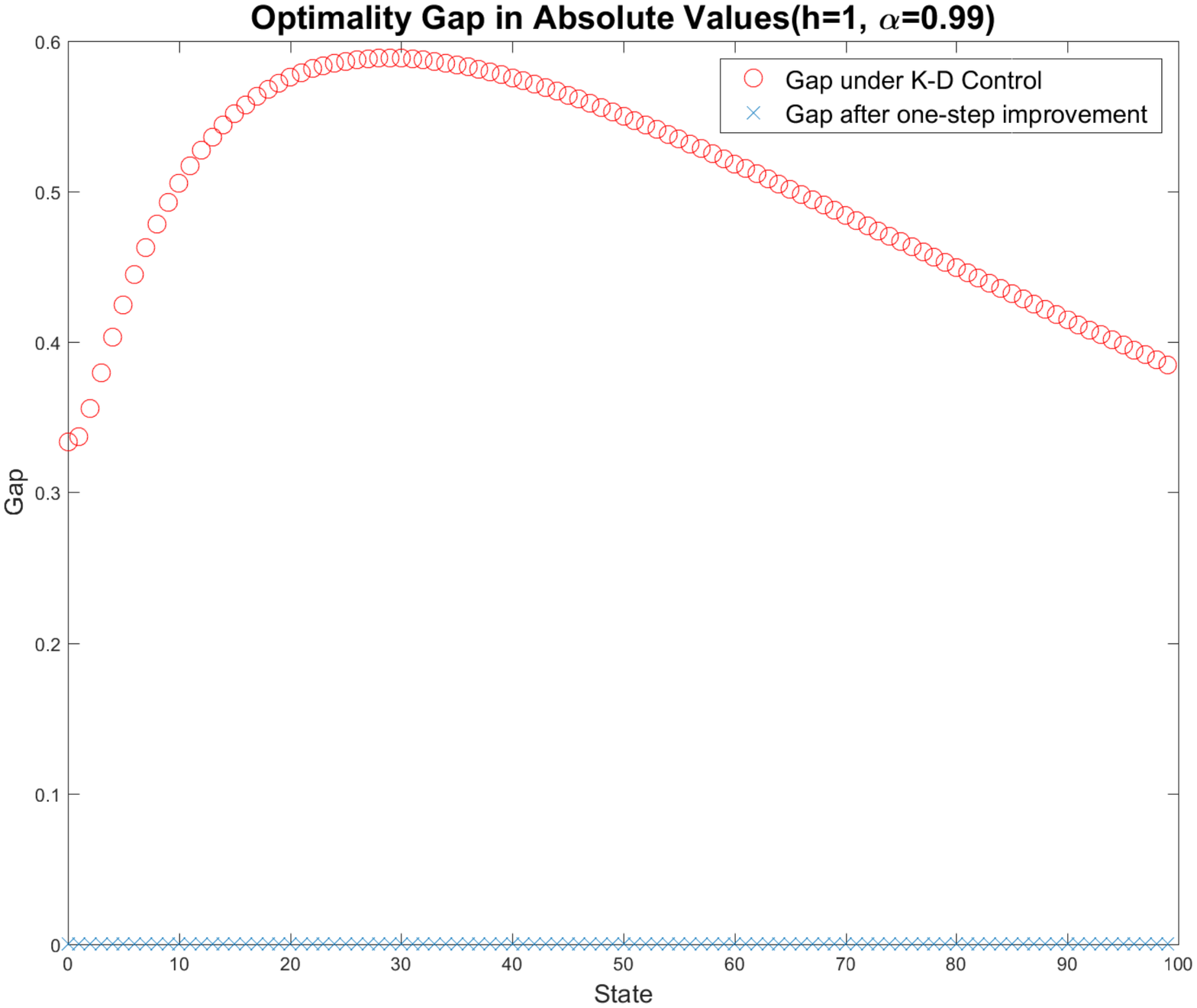}		\includegraphics[scale=0.29]{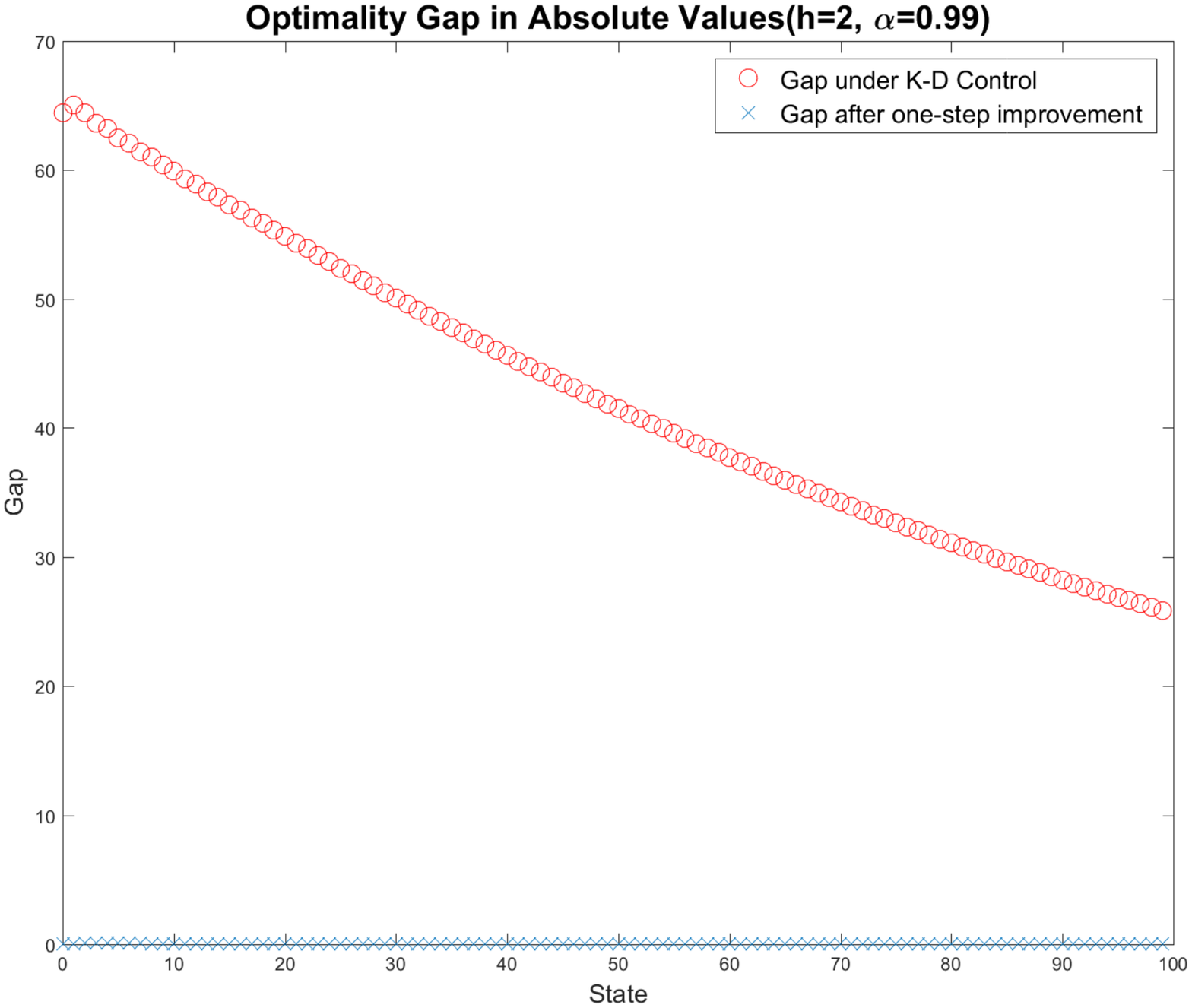}	
					\end{center} \caption{Algorithm performance for the service-rate control problem: (TOP) The optimal value compared against the approximate values. The lines are indistinguishable at the scale of the optimal value. (BOTTOM) The absolute optimality gap: the difference between the cost under the proposed policy and the optimal cost. The blue series is the result after one-step improvement. It is indistinguishable from $0$. \label{fig:barSSQgap}} 	
				\end{figure}

\begin{figure}[hh!]
					\begin{center}
						\includegraphics[scale=0.29]{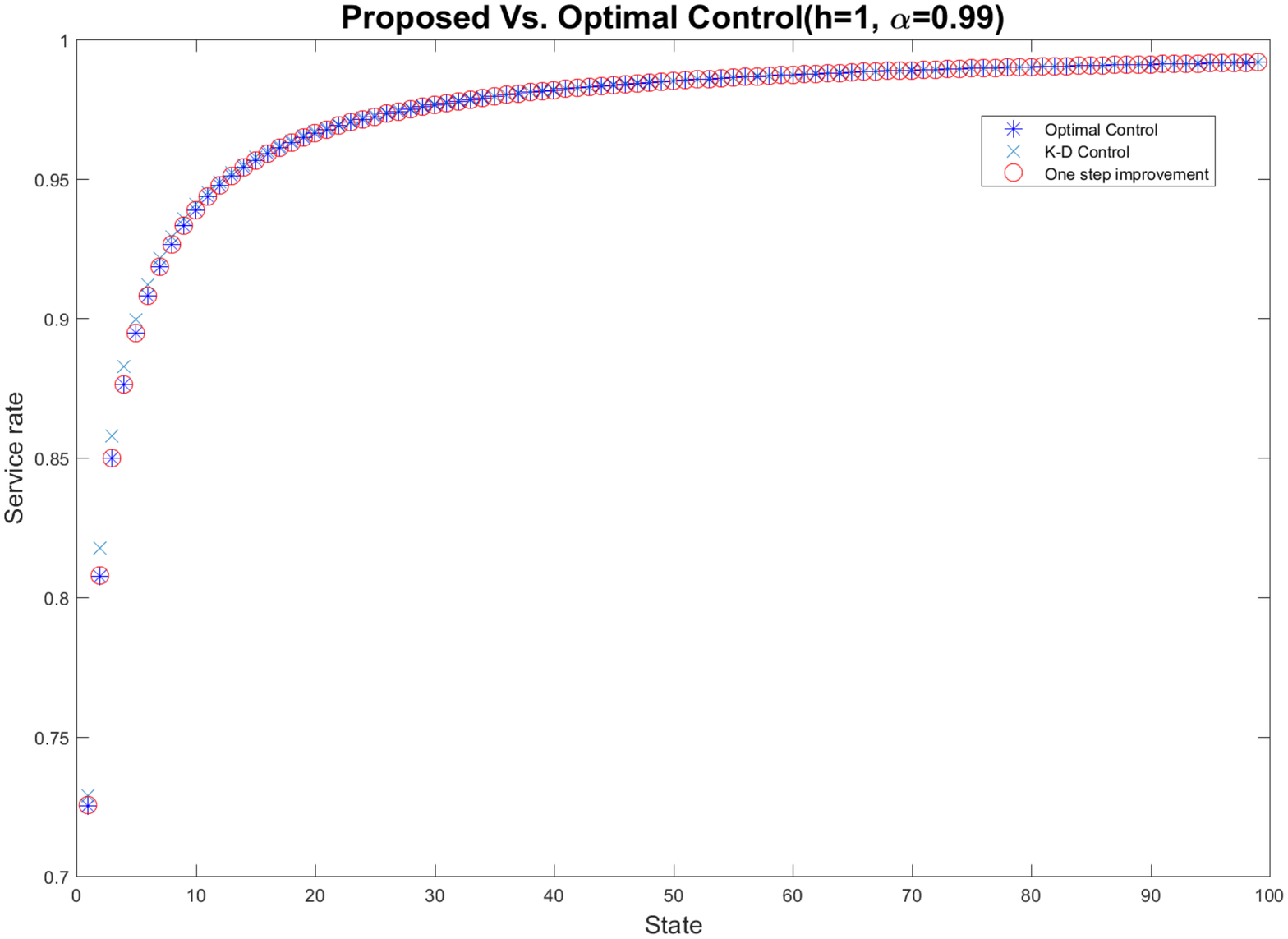}		\includegraphics[scale=0.29]{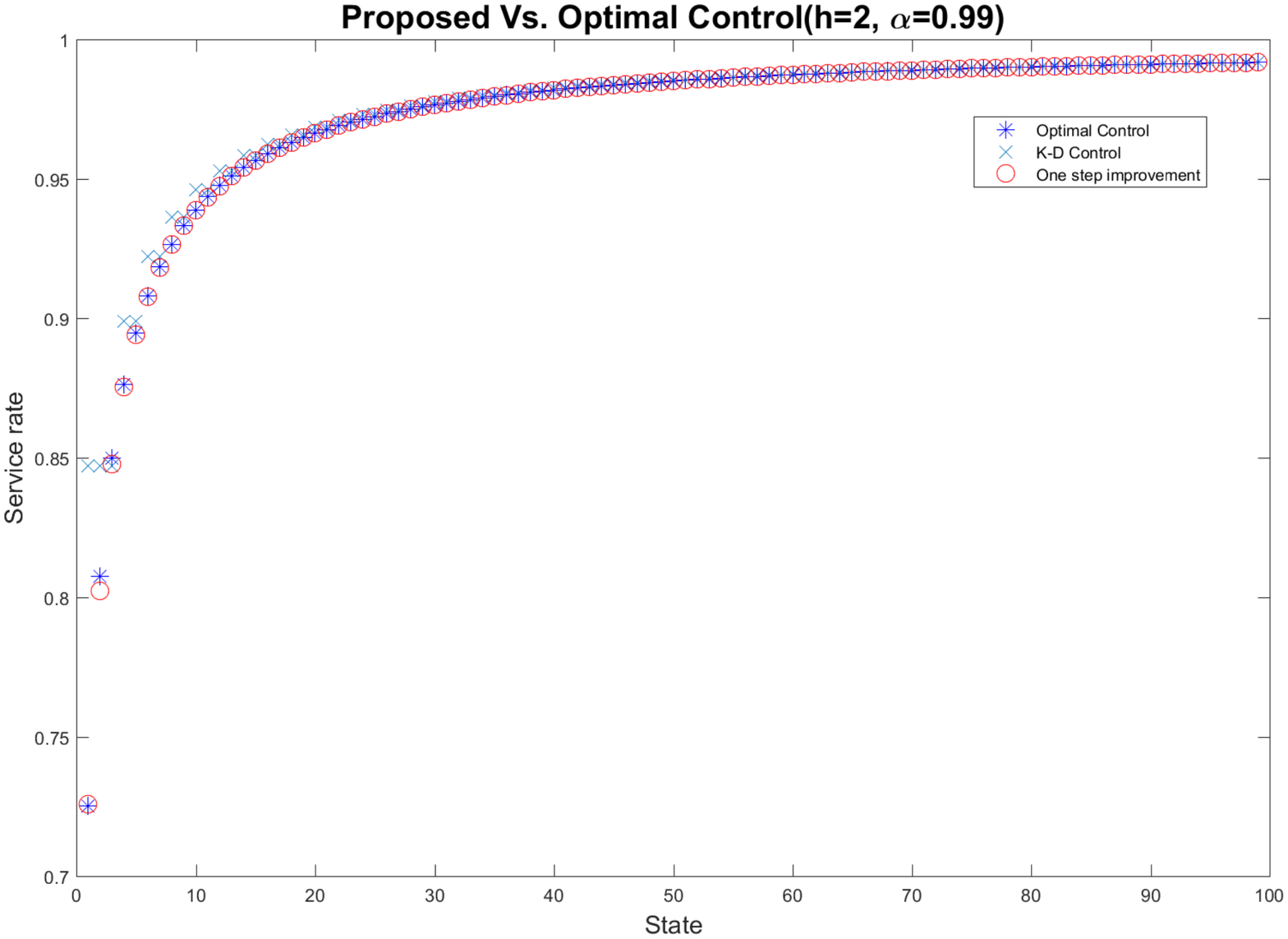}	
					\end{center}
					\caption{Control comparison for the service-rate control problem: (LEFT) $h=1$, (RIGHT) $h=2$. The control obtained from one-step improvement is almost identical (at all states) to the optimal control. 	\label{fig:barSSQcontrol}}
					
				\end{figure}

				In Figure \ref{fig:barSSQgap} we plot the absolute (rather than relative) optimality gap $V_{U^h}(x)-V_{*}(x)$ for $\alpha=0.99$ and discretization $h$ equal to 1 and 2. It is important that even in the case of $h=2$, an optimality gap of 30 (at $x=100$) is negligible relative to the optimal value at that state, which is of the order of $3*10^5$. More impressive is the performance after one-step policy improvement. The greedy policy achieves an optimality gap  that is indistinguishable from $0$. This result is explained by Figure \ref{fig:barSSQcontrol}, where we report the comparison of {\em actions}. The plot also includes the control after one-step policy improvement starting at the K-D chain interpolated value.
								
				Per Theorem \ref{thm:diffusion}, the error should be of the order of the integrated third derivative. If this derivative is uniformly bounded by $\Gamma$, the optimality gap must be smaller than (or equal to) $\Gamma/(1-\alpha)$. The central differences proxy with $h=1$ for the third derivative of $\hV_{\ast}(x)$ is given by $$\hV_{\ast}'''(x)\approx \frac{1}{2}\left(V_{\ast}^{h}(x+2)- 2V_{\ast}^{h}(x+1) + 2V_{\ast}^{h}(x-1) - V_{\ast}^{h}(x-2)\right).$$
							The peak of this proxy for $h=1$ and $\alpha=0.99$ is 1.8, generating in our theorem an error bound of $\frac{1}{1-\alpha}*1.8=180$, which is 0.0006 of the value of more than 270,000 for the state $x=\frac{1}{1-\alpha}=100$.
				
				\subsection{An inventory problem}
				
				In the inventory problem we study next, the optimal policy is a so-called order-up-to level policy, implying large jumps in some states. This problem seems to pose a challenge to our bounds, which depend on the size of the maximal jump.
				
				Period $t$ demand $D_t$ is drawn from a Poisson distribution with mean $\Ex[D_t]=\lambda$. Demand is independent across periods. There is a backlog cost of $b$, a per-period cost $H$ for holding a unit in inventory, and a per-unit order cost of $c$. The lead-time is $0$.

				The state is the inventory position. The per-period cost is given by
				$$r_{u}(x)=cu+ H\Ex[(x+u-D)^+]+b\Ex[(D-(x+u))^{+}],$$ where the action $u$ is the amount ordered. Orders are placed (and received) and then demand is realized and backorder and holding costs are incurred. Transitions have the form
				$$X_t\to (X_t+U_t-D_t),$$ where $U_t$ is the order quantity in period $t$. The Bellman equation is given by $$V(x)=\min_{u\in\mathbb{Z}_+}\{r_{u}(x)+\alpha \Ex[V(x+u-D)]\}.$$
				
				The drift and diffusion coefficient are given by
				\begin{align*} \mu_u\equiv \mu_{u}(x)&=\Ex[X_1-x]=\Ex[x+u-D-x]=u-\lambda, \mbox{ and }\\ \sigma^2_u\equiv \sigma^2_{u}(x)&=\Ex[(X_1-x)^2]=\Ex[(x+u-D-x)^2]=(u-\lambda)^2+\lambda,\end{align*} so that the TCP is
				$$0=\min_{u\in\mathbb{Z}_+}\left\{r_{u}(x)+\alpha\left( \mu_{u}V'(x)+\frac{1}{2}\sigma_{u}V''(x)\right)-(1-\alpha)V(x)\right\}, ~x\in\mathbb{R}.$$
				Notice that $\sigma^2_u\geq \lambda$ for all $x$ and $u\in\mathbb{Z}$ so that strict ellipticity holds. Because the state space includes all the integers, there are no boundary conditions here. These are artificially introduced in our numerical computation to make the state space finite. Specifically, we truncate the space at state $M$ and $-M$ where $\mu_u(-M)=u$, $\sigma_u^2(-M)=u^2$ and $\mu_u(M)=-\lambda$, $\sigma_u^2(M)=\lambda+\lambda^2$. We then introduce the boundary condition $V'(M)=V'(-M)=0$. Notice that for constant $u$, the drift and diffusion coefficient are constant in $x$ and in particular Lipschitz. The cost function is, as well, Lipschitz in $x$ uniformly in $u$. The existence of a solution follows from Lemma \ref{lem:PDEbounds}. 
				
\begin{figure} \hspace*{-0.5cm}
					\includegraphics[scale=0.3]{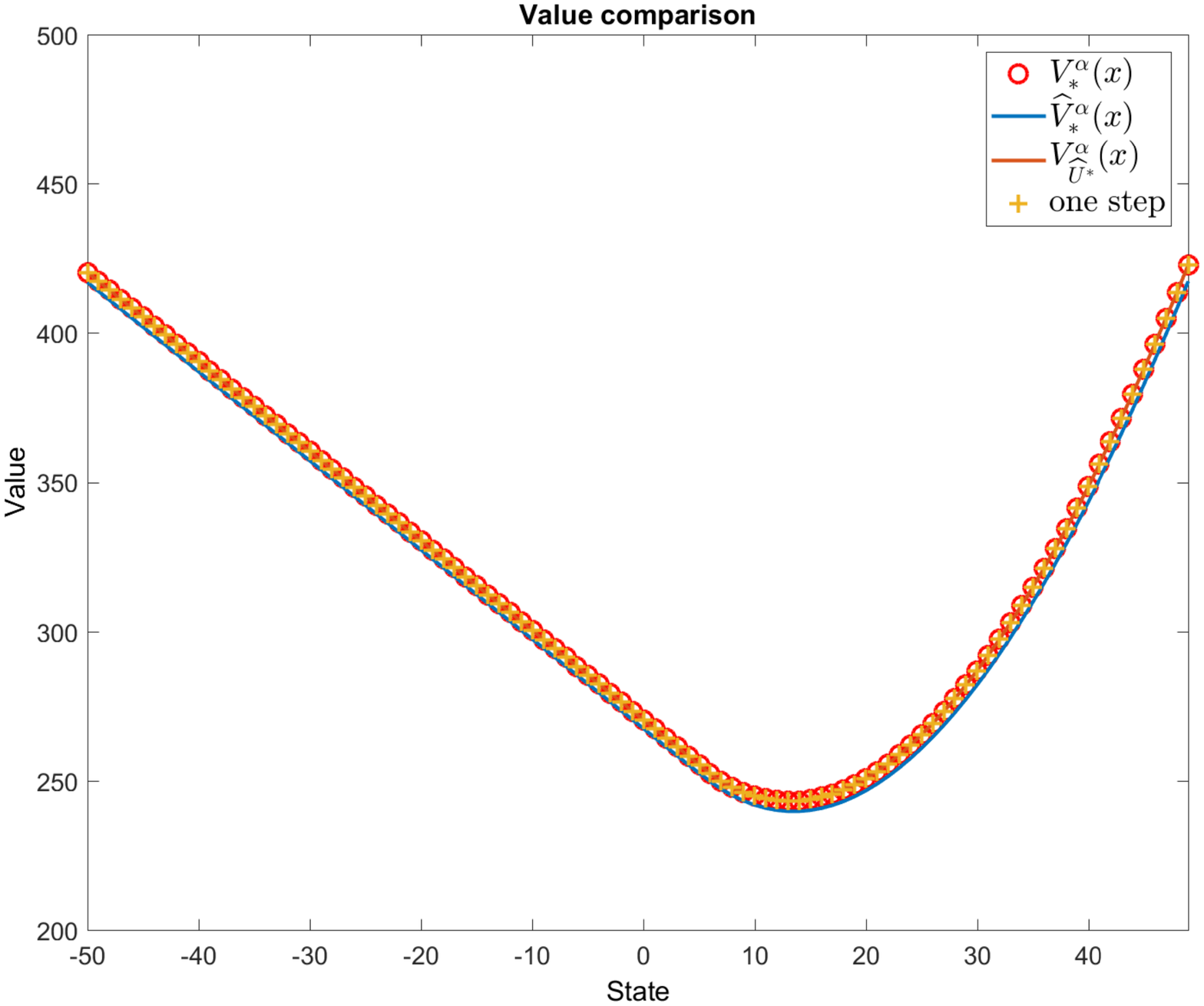}	\includegraphics[scale=0.3]{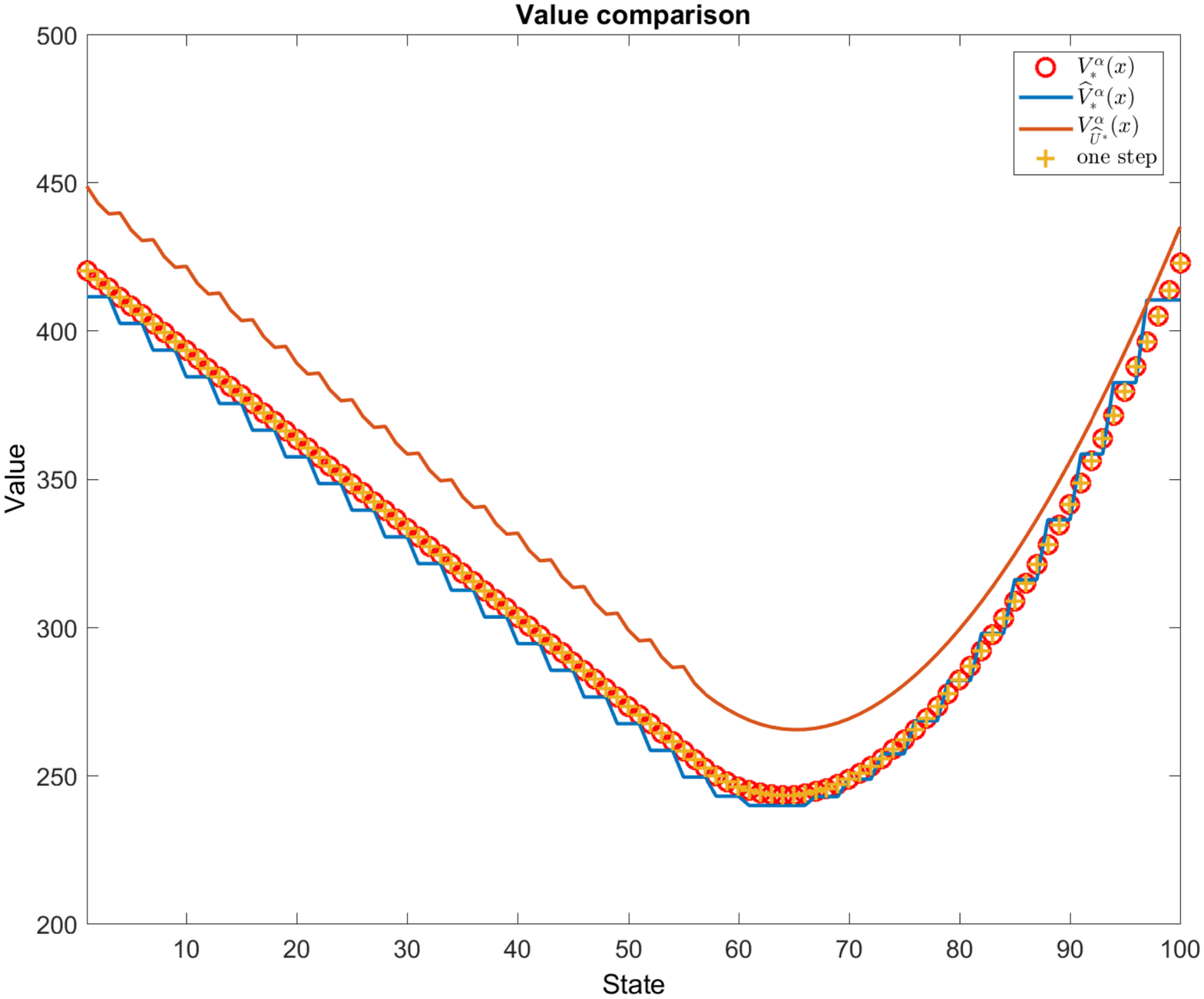}
					\caption{Algorithm performance for the inventory problem: (LEFT) $h=1$, (RIGHT) $h=3$. The performance of the K-D {\em control} is not as good, but one-step improvement relative to the K-D value function $V_{\ast}^{h,\alpha}$ generates an extremely accurate value for all states. \label{fig:inventory}}\end{figure}

				We use the K-D chain for value of coarseness $h\in \{1,3\}$. For each $h$ we use (with some abuse of notation)  $\hV_{\ast}^{\alpha}(x)$ for the value from the K-D approximation, which is a proxy for the TCP value. For $h=3$ we extended it to the integers in a piecewise constant manner. We also take the control $\hU_h^*$ and interpolate it to the whole state space in a piecewise constant manner. $V_{\hU_*}(x)$ is the value in the original chain when using this control. Finally, the value after one-step improvement is the infinite-horizon discounted reward under $U^h$---the greedy control relative to $V_{\ast}^{h,\alpha}(x)$. Figure \ref{fig:inventory} displays the computational results.
				
				\subsection{A routing problem \label{sec:exmultipool}}
				This example is based on the inpatient-management queueing model studied in \cite{pengyi}. The task is to optimally route patients from dedicated queues to internal hospital wards so as to minimize the aggregate cost of holding and routing.
				
				The dynamics of the queues are modeled via a discrete-time queueing model with $J$ server pools (the internal wards), where pool $j$ has $N_j$ servers (beds), a dedicated inflow of customers, and an infinite-size dedicated buffer. This buffer will be truncated for the numerical experiments. Customers from the $j^{th}$ inflow are referred to as type-$j$ customers. These customers can be served by their dedicated pool $j$ but also by other pools.
				
				We let $X_j(t)$ be the number of customers of type $j$ in the system at time $t = 0, 1, 2, \ldots$, and let $X(t) \in \Z_+^{J}$ be the vector whose components are $X_j(t)$. A customer in the system can either be in service, or waiting in a buffer to be served.
				
				The (controlled) chain's evolution is as follows: At the start of time-period $t$, customers waiting in buffer $j$ enter service in pool $j$ until the buffer is emptied or all idle servers are taken. If any customers remain in buffer $j$, we proceed to the \emph{overflow} decision. This is the overflow control. At a cost of $B_{ij}$ per customer, we can choose to assign a customer waiting in buffer $i$ to immediately enter service in pool $j \neq i$ if that pool has an idle server available. We can also decide not to overflow any customers. We let $U_{ij}(t) = U_{ij}(X(t))$ be the number of customers overflown from buffer $i$ to pool $j$ in time period $t$. After overflows are executed, a holding cost $H_i$ per customer waiting in buffer $i$ is incurred. Next, departures are resolved: a type $j$ customer in service completes service and leaves the system with probability $p_j$ (service time is geometric with mean $1/p_j$). Otherwise, the customer remains in service until the next period.
				
				After departures are resolved, new arrivals occur: the number of type-$j$ customers to arrive per period is Poisson distributed with mean $\lambda_j$. Arrivals are independent across types and across time periods. An incoming customer will either occupy an idle server in her dedicated pool or, if there are no such servers, enter the buffer and wait for service.
				
				Under a stationary control $u$, $X(t)$ satisfies the dynamics
				\begin{align}
					X_i(t) = X_i^P(t-1)+ A_i(t-1) - D_i(X_i^P(t-1)), \label{eq:dynamicsQ}
				\end{align}
				where
				\begin{align*}
					X_i^P(t-1) =  X_i(t-1) + \sum_{j \neq i} U_{ji}(X(t-1)) - \sum_{j \neq i} U_{ij}(X(t-1))
				\end{align*}
				is the post-action state in period $t-1$,
				\begin{align*}
					D_i(x) \sim \text{Binomial}\big(x\wedge N_i, p_i\big)
				\end{align*}
				is the number of departures from pool $i$, and
				\begin{align*}
					A_i(t) \sim \text{Poisson}(\lambda_i)
				\end{align*}
				is the number of new unblocked type-$i$ arrivals, independent across periods.
				
				The state space is $\Z_+^{J}$ and the action space is \begin{align*}
					\mathcal{U}(x) = \Big\{u \in \Z_+^{J\times J}\  \Big|\  \sum_{j\neq i} u_{ji} \leq (N_i - x_i)^{+},  \text{ and } \sum_{j \neq i} u_{ij} \leq (x_i - N_i)^{+} \Big\}.
				\end{align*}
				The first constraint guarantees that the overflow to pool $i$ does not exceed the number of idle servers, and the second constraint guarantees that the number of overflowed type-$i$ customers does not exceed the number of customers waiting in buffer $i$. The action space satisfies the structure $\mathcal{U}(x)=\mathbb{D}\bigcap \{u: Au\leq b(x)\}$. Finally, the per-period cost includes the cost of overflow and linear cost of holding and is given by \begin{align}
					r_u(x) =
					\sum_{i} \sum_{j \neq i} B_{ij} u_{ij} + \sum_{i} H_i \times \big(x_i - \sum_{j \neq i} u_{ij} - N_i\big)^{+}. \label{eq:cost}
				\end{align} The goal is to make overflow decisions that minimize the expected infinite-horizon discounted cost.
				
				The Bellman equation for this dynamic program is computationally challenging. A modest system with $J = 3$, $N_i\equiv  40$, and $M = 60$ has above one million states. Moreover, depending on the policy $U(x)$, a state can have many ``neighboring'' states, making the transition probability matrix dense and expensive to store. Finally, because actions are discrete, the only option is exhaustive search over the very large action space: we have to decide how many customers to overflow from {\em each} buffer into {\em each} pool.
				
				We next construct the TCP's ingredients. For $x \in\mathbb{S}$ and $u \in \mathcal{U}(x)$, let $x_i^{P}(u) = x_i + \sum_{j \neq i} (u_{ji} - u_{ij})$. Then,
				\begin{align*}
					\big(\mu_{u}(x) \big)_{i} =&\ \sum_{j \neq i} u_{ji} - \sum_{j \neq i} u_{ij} + \E[A_i(t)] -  \E[ D_i(x_i^{P}(u))] \\
					=&\ \sum_{j \neq i} (u_{ji} - u_{ij}) + \lambda_i - p_i \Big( \big(\lceil x_i\rceil  + \sum_{j \neq i} (u_{ji} - u_{ij}) \big) \wedge N_i\Big)=:(f_{\mu}(u,x))_i,
				\end{align*}
				and
				\begin{align*}
					\big(\sigma^2_{u}(x)\big)_{ij} &=\ \E \bigg[\Big(\sum_{k \neq i} (u_{ki} - u_{ik}) + A_i(t) -  D_i(x_i^{P}(u))\Big) \Big(\sum_{k \neq j} (u_{kj} - u_{jk}) + A_j(t) -  D_j(x_j^{P}(u))\Big) \bigg]\\
					&=:(f_{\sigma}(u,x))_{ij}.
				\end{align*}
				
				We use the oblique derivative condition
				$$\eta(x)'D\hV(x)=0,~x\in \partial\mathbb{R}_+^d$$ where $\eta_i(x)=p_i$ if $x_i=0$ and is $0$ otherwise.
				
				This is grounded in intuition about the ``pushback'' at $0$, but is also mathematically supported by choosing the suitable extension of $\mu_u$ and $\sigma^2_u$. The function $f_{\mu}(u,x)$ is well defined for all $x\in\mathbb{Z}_+^d$ and $u\in \mathbb{D}$ and can be continuously extended to $x\in\mathbb{R}_{++}^d$ in multiple ways. We choose to extend $\lceil x_i \rceil$ so that it is constant (and equal to $1$) for all $x_i\in (0,1]$. We extend $f_{\sigma}(u,x)$ so that it is continuous on all of (not just the interior of) $\mathbb{R}_+^d$. For $x$ with $x_i<N_i$, $\mathcal{U}(x)$ contains only $u$ with  $u_{ij}=0$ for all $j\neq i$.
				If $\hU_*$ is piecewise constant and continuous at the boundary\footnote{$\hU_*(x_n)\rightarrow \hU_*(x)$ for all $x\in \partial\mathbb{R}_+^d$ and sequences $\{x_n\}$ with $x_n\in \mathbb{R}_{++}^d$ and $x_n\rightarrow x$.} then  $(\eta(x))_i=(\mu_{\hU_*}(x_{-i},0+))_i-(\mu_{\hU_*}(x_{-i},0))_i=p_i$.
				
				We also tried the more direct FOT boundary conditions for this example, obtaining similar performance to what is reported in Table \ref{table:threepool1}-\ref{table:threepool3} below.
				
				In our computational experiments we truncate the state space by using finite buffers and truncating arrivals in an intuitive way. We use exhaustive search over $u\in\mathcal{U}(x)$ in the policy improvement search, rather than relaxing the integrality constraints; see, e.g., \cite{MoalKumaVanr2008}. We do so because we wish to capture the error induced by the Taylor expansion without confounding it by approximations to the action space. Still, the computational savings of TAPI were significant: a single iteration of TAPI took a few minutes compared to about 3 hours for a PI iteration, leading to a reduction of total running time from over 15 hours to below 10 minutes.
				
				Tables~\ref{table:threepool1}-\ref{table:threepool3} present the results of applying TAPI to multiple three-dimensional $(J=3)$ instances of the model. The table shows the difference between the value function under the proposed policies, $V_{U^h_*}^{\alpha}(x)$, and the actual optimal value $V_{\ast}^{\alpha}(x)$. The maximal relative error is computed by $\max_{x \in \mathbb{S}}\frac{|V_{U^h_*}(x)-V_{\ast}(x)|}{V_{\ast}(x)}$, where $U^h_*$ is the policy suggested by the approximation algorithm. The mean relative error column reports  $\frac{1}{|\mathbb{S}|} \sum_{x \in \mathbb{S}}\frac{|V_{U^h_*}(x)-V_{\ast}(x)|}{V_{\ast}(x)}$.
				
				\begin{table}[hh!]\caption{$N_1=N_2=N_3=10$, $M = 14$,  $(p_1, p_2,p_3) = (0.8,0.8,0.8)$,  $(H_1, H_2,H_3) = (1,2,3)$, $(B_{12}, B_{13}) = (1,1)$, $(B_{21}, B_{23}) = (4,1)$, and $(B_{31}, B_{32}) = (2,1)$.\label{table:threepool1}} \begin{center}
						\begin{tabular}{c|c|c|c|c|c}
							& & \multicolumn{2}{c}{$\lambda_i = 0.7 N_i p_i$} & \multicolumn{2}{c}{$\lambda_i = 0.8 N_i p_i$}\\
							\hline
							$\alpha$	& h	& Maximum rel. error &
							Mean rel.\ error & Maximum rel. error &
							Mean rel.\ error \\
							\hline
							0.99&	2& 0.023 & 0.001 &	0.019	&0.002\\
							&	4& 0.016 & 0.0004 &	0.017	&0.0005 \\
							&	8& 0.019 & 0.0004 &	0.014	&0.005 \\\hline
							0.999&	2& 0.004 & 0.001 &	0.004	&0.002  \\
							&	4& 0.003 & 0.0004 &	0.003	&0.0005 \\
							&	8& 0.003 & 0.0004 &	0.007	&0.006 \\\hline
					\end{tabular}\end{center}
				\end{table}
				
				\begin{table}[h]
					\caption{$N_1=N_2=N_3=10$, $M = 14$,  $(p_1, p_2,p_3) = (0.4, 0.6, 0.1)$ , $(H_1, H_2,H_3) = (10,2,6)$, $(B_{12}, B_{13}) = (5, 2)$, $(B_{21}, B_{23}) = (3,7)$, and $(B_{31}, B_{32}) = (7, 9)$.\label{table:threepool2}}
					\begin{center}
						\begin{tabular}{c|c|c|c|c|c}
							& & \multicolumn{2}{c}{$\lambda_i = 0.7 N_i p_i$} & \multicolumn{2}{c}{$\lambda_i = 0.8 N_i p_i$}\\
							\hline
							$\alpha$	& h	& Maximum rel. error &
							Mean rel.\ error & Maximum rel. error &
							Mean rel.\ error \\
							\hline
							0.99&	2&	0.206	&0.011 & 0.096 & 0.005 \\
							&	4&	0.184	&0.012 & 0.120 & 0.005 \\
							&	8&	0.110	&0.032 & 0.055 & 0.011 \\\hline
							0.999&	2&	0.037	&0.007 & 0.028 & 0.002 \\
							&	4&	0.036	&0.010 & 0.014 & 0.003 \\
							&	8&	0.051	&0.039 & 0.016 & 0.011 \\\hline
						\end{tabular}
					\end{center}
				\end{table}
				
				\begin{table}[h]
					\caption{ $N_1=N_2=N_3=10$, $M = 14$,  $(p_1, p_2,p_3) = (0.2, 0.7, 0.5)$, $(H_1, H_2,H_3) = (1, 1, 4)$, $(B_{12}, B_{13}) = (5, 2)$, $(B_{21}, B_{23}) = (7, 1)$, and $(B_{31}, B_{32}) = (7, 9)$.
						\label{table:threepool3}	}\begin{center}
						\begin{tabular}{c|c|c|c|c|c}
							& & \multicolumn{2}{c}{$\lambda_i = 0.7 N_i p_i$} & \multicolumn{2}{c}{$\lambda_i = 0.8 N_i p_i$}\\
							\hline
							$\alpha$	& h	& Maximum rel. error &
							Mean rel.\ error & Maximum rel. error &
							Mean rel.\ error \\
							\hline
							0.99&	2&	0.077	&0.021 & 0.053 & 0.009 \\
							&	4&	0.060	&0.009 & 0.039 & 0.005 \\
							&	8&	0.069	&0.042 & 0.025 & 0.014 \\\hline
							0.999&	2&	0.034	&0.024 & 0.016 & 0.009 \\
							&	4&	0.019	&0.010 & 0.009 & 0.004 \\
							&	8&	0.059	&0.055 & 0.017 & 0.016 \\\hline
						\end{tabular}
					\end{center}
				\end{table}
				
\begin{figure}
					\hspace{-1cm}
					\begin{minipage}{0.55\textwidth}
						\includegraphics[scale=0.35]{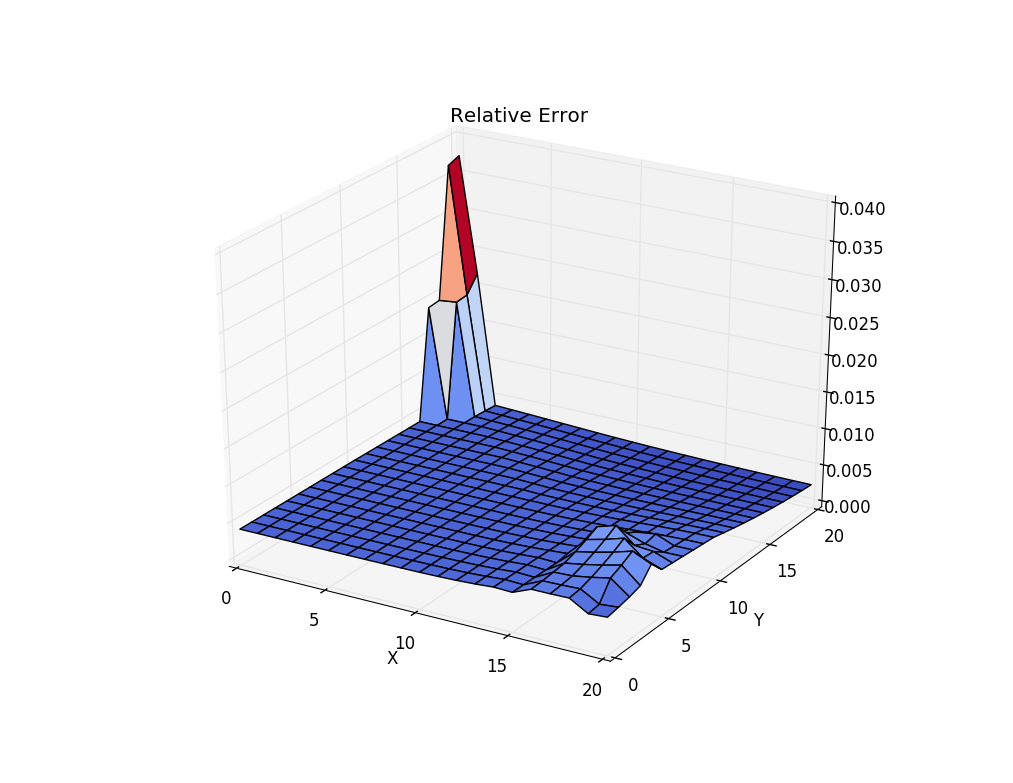}\end{minipage}\begin{minipage}{0.35\textwidth} 		
						\includegraphics[scale=0.35]{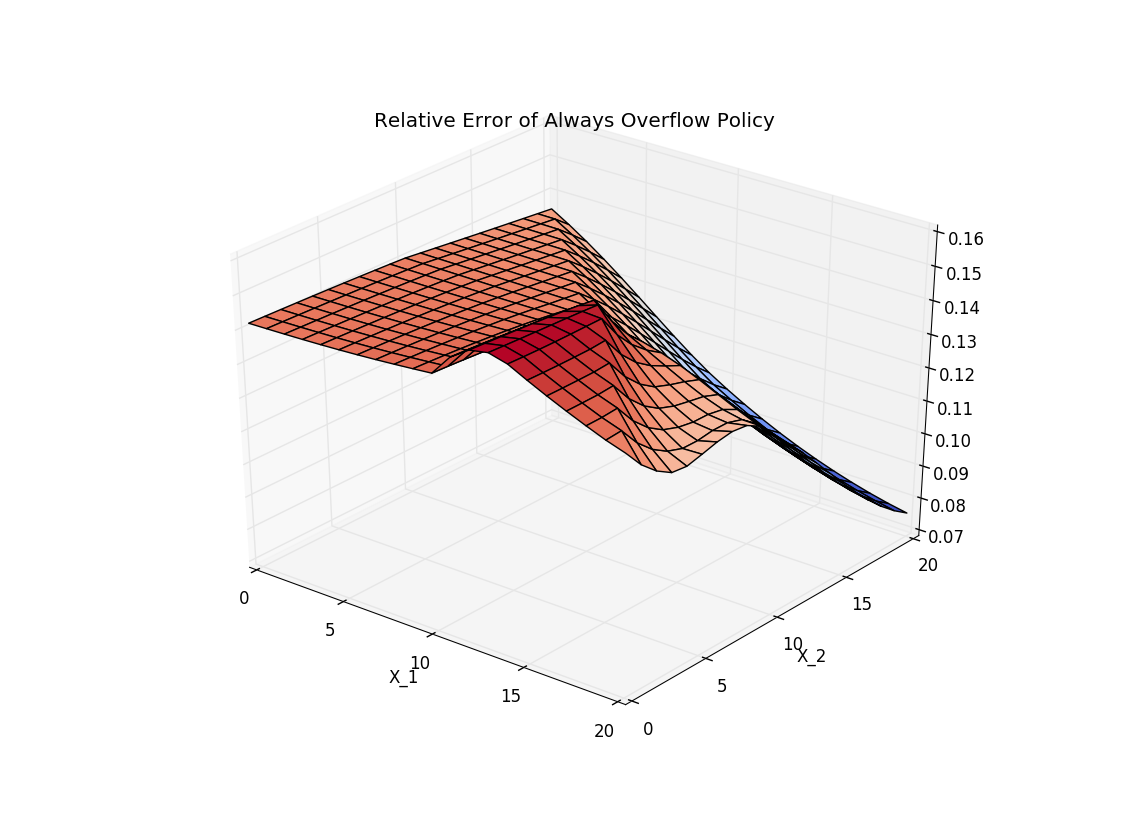}
					\end{minipage}
					\caption{Concentration of error: $J=2$, $N_1=N_2=10$, $M = 10$, $(p_1, p_2) = (0.56,0.56)$, $(H_1, H_2) = (1,4)$, $(B_{12}, B_{21}) = (5,1)$ and $\alpha=0.99$, $h =2$ and $\lambda_i = 0.8 N_i p_i$. On the left we plot the relative error $\frac{|V_{U^h_*}(x)-V_{\ast}(x)|}{V_{\ast}(x)}$ over the entire domain $\mathbb{Z}_+^2$. The $3.7\%$ maximal relative error reported in Table~\ref{tab:multipool2} is caused by only a small portion of the state space. For contrast, the plot on the right shows the distance to $V_{*}(x)$ of the value function under the (suboptimal) policy that overflows as many customers as possible. \label{fig:badpol}}
				\end{figure}

				While the maximal error can be fairly large, the mean relative error is rather small. In Figure \ref{fig:badpol} we use a two-dimensional example to visualize this fact. Even in cases where the {\em maximal} error is as large as $3.7\%$, such errors are confined to a very small portion of the state space (close to a boundary) and are much smaller in most of the state space. This is captured in Figure \ref{fig:badpol}, where we plot the relative error in a two-dimensional ($J=2$) case.
				
				Finally, in reference to Remark \ref{rem:exactimprovement}, we compare the performance of TAPI with a heuristic in which the policy-improvement is executed exactly rather than approximately. Introducing exact improvement, while having no convergence guarantees, can result in better performance; see the ``+exact improv.'' column in Table \ref{tab:multipool2}. This performance, is however, matched by using TAPI as is (with approximation improvement) and adding at the very end a {\em single} exact policy improvement step; see the column ``one step'' in the same table.
				
				\begin{table}[!hh] \TABLE{Relative error: TAPI vs. TAPI with exact-improvement step\label{tab:multipool2}}{
						\begin{tabular}{c|c|c|c|c}
							\multicolumn{4}{c}{$\lambda_i = 0.8N_i p_i$}\\
							$\alpha$	& h	& TAPI &
							$+$ exact improv& one step  \\
							\hline
							0.99&	1&	0.0376 & 0.0086	& 0.0095\\
							&	2 &0.0373&	0.0081 & 0.0088	\\
							&	4&0.0346&	0.0067 & 0.0079	\\\hline
							0.999&	1&0.0093 &	0.0033	&  0.0051 \\
							&	2&0.0082&	0.0031	& 0.0045 \\
							&	4&0.0048&	0.0023 & 0.0032	\\\hline
						\end{tabular} \quad
						\begin{tabular}{c|c|c|c|c}
							\multicolumn{4}{c}{$\lambda_i = N_i p_i$}\\
							$\alpha$	& h	& TAPI &
							$+$ exact improv & one step  \\
							\hline
							0.99&	1&	0.0103	&0.0089 & 0.0083\\
							&	2&0.0107&	0.0069	& 0.0100 \\
							&	4&0.0013 &	0.0014	& 0.0014 \\\hline
							0.999&	1&0.0013 &	0.0026 & 0.0030 \\
							&	2&0.0012&	0.0026 & 0.0043	\\
							&	4&0.0001&	0.0002 & 0.0002	\\\hline
					\end{tabular}}
					{\small Parameters: $J=2$, $N_1=N_2=10$, $M = 10$, $(p_1, p_2) = (0.56,0.56)$, $(H_1, H_2) = (1,4)$, and $(B_{12}, B_{21}) = (5,1)$. }
				\end{table}

				\begin{remark}[smoothing the state space] We use this example, with $J=2$, to illustrate the domain-smoothing described in \S \ref{subsec:bounds}. We replace the corner point $(0,0)$ in the Markov chain's state space with the point $\tilde{0}=(\epsilon,\epsilon)$ so that we can ``pass'' a smooth boundary through this point. Let us denote the drift and diffusion coefficient by $\tilde{\mu}_u$ and $\tilde{\sigma}^2_u$. Because $P^u_{\tilde{0},y}=P_{0,y}^u$, we have  $$(\tilde{\mu}_u)_i(\tilde{0})=\sum_{y} P_{0,y}^u(y_i-\epsilon)=
					\sum_{y} P_{0,y}^u (y_i-0)-\epsilon=\mu_i(0)-\epsilon=\lambda_i-\epsilon,$$
					and $$(\tilde{\sigma}^2_u)_{ij}=\sum_{y} P_{0,y}^u(y_i-\epsilon)(y_j-\epsilon)=\Ex[(A_i(t)-\epsilon)(A_j(t)-\epsilon)].$$
For all other states $(\tilde{\mu}_u)_i(x)=(\mu_u)_i(x)+P_{x,0}^u\epsilon$ and $(\tilde{\sigma}^2_u)_{ij}(x)=(\sigma^2_u)_{ij}(x)+P_{x,0}^u\epsilon(\epsilon-x_i-x_j)$. Notice that $\epsilon$ appears only in the new ``corner''; there are no states of the form $(\epsilon,x_2)$ or $(x_1,\epsilon)$ for $x_1,x_2\geq 1$. Finally, the vector $\eta$ is given by $(p_1,0)$ when $x_1=0$ and $x_2\geq 1$ and by $(0, p_2)$ when $x_2=0$ and $x_1\geq 1$. It is $\eta(x)=(p_1,p_2)$ at $x=\tilde{0}$. We then extended to the curved boundary in a smooth way. All of these have obvious extensions to more than two dimensions.
					\label{rem:smoothing} \hfill \bsq
				\end{remark}

\section{Tayloring and Queues in Heavy-traffic\label{sec:HT}}

In this short section we illustrate, informally and via the simplest of examples, the connection between Tayloring and heavy-traffic approximations.

Consider the discrete-time queue with holding cost $r(x)=x$, $P_{x,x+1}=\lambda$, and $P_{x,x-1}=\mu:=1-\lambda>0.5$. These parameters are fixed, so the question is one of performance approximation, but it is sufficient for illustration; Example \ref{example:pdescaling} in the appendix adds control.

 Let us consider two avenues to approximate the infinite-horizon discounted reward.

\paragraph{Process convergence.} Let $\mu\downarrow \frac{1}{2}$ so that $\rho=\frac{\lambda}{\mu}=\frac{1-\mu}{\mu}\uparrow 1$---the queue is in {\em heavy-traffic}; see, e.g., \cite[Chapter 9]{whitt2002stochastic}. Then, it is a standard heavy-traffic result that extending time by $(1-\rho)^{-2}$ and shrinking space by $(1-\rho)$,
$$(1-\rho)X_{\lceil (1-\rho)^{-2}t\rceil} \approx \hat{X}(t),$$  where $\hat{X}(t)$ is a so-called reflected Brownian motion with drift $-\frac{1}{2}$ and diffusion coefficient $\sigma^2\equiv 1$. This result is formalized by weak convergence arguments. For purposes of this discussion it suffices to note that, with discount factor equal to $1$, the infinite-horizon discounted reward of the diffusion $\hV(x)=\Ex_x[\int_0^{\infty}e^{-s}r(\hX(s))ds]$ solves the ODE
$$0=x-\frac{1}{2}V'(x)+\frac{1}{2}V''(x)-V(x),~ V'(0)=0.$$ It is then possible to show that if one takes the discount factor to be $$\alpha_{\rho}=1-(1-\rho)^2,$$ then as $\rho\uparrow 1$
$$\hV(x)\approx (1-\rho)^3\Ex_x\left[\sum_{t=0}^{\infty}  \alpha_{\rho}^t X_t\right].$$ In words, $\hV$ approximates the value function up to an error that is small relative to $(1-\rho)^{-3}$.

\paragraph{Tayloring.}
The TCP, for given $\lambda$ and $\mu$, is given by
$$0=x+\alpha((\lambda-\mu)V'(x)+\frac{1}{2}V''(x))-(1-\alpha)V(x),$$ and $V'(0)=0$. With $\lambda<\mu$, this ODE has the solution $$\hV^{\alpha}(x)=-\alpha\frac{\mu-\lambda}{(1-\alpha)^2}+\frac{x}{1-\alpha}+c_1e^{\gamma_- x},$$ where
$\gamma_{-}=-\sqrt{(\mu-\lambda)^2+2\frac{1}{\alpha}(1-\alpha)}+(\mu-\lambda)<0$, and  $c_1=-\frac{1}{(1-\alpha) \gamma_-}$. The parameter $\gamma_-$ depends on $\alpha$ and satisfies $\gamma_-\approx - \frac{1-\alpha}{\mu-\lambda}$ for $\alpha\approx 1$ so that $\hV^{\alpha}(x)\approx \frac{\mu(1-\rho)}{1-\alpha}+\frac{x}{1-\alpha},$ as $\alpha\uparrow 1$. We also have
$$|D^3\hV(x)|=\frac{\gamma_-^2}{1-\alpha}e^{\gamma_- x}\leq \frac{\gamma_-^2}{1-\alpha},$$
so that
$$\Ex_x\left[\sum_{t=0}^{\infty}\alpha^t |D^3\hV^{\alpha}|_{X_t\pm 1}^*\right]\leq \frac{1}{1-\alpha}\sup_{x\geq 0}|D^3\hV^{\alpha}(x)|\leq \frac{\gamma_-^2}{(1-\alpha)^2}\approx \frac{1}{(\mu-\lambda)^2}\mbox{ as }\alpha\uparrow 1.$$
If $\mu$ and $\lambda$ are chosen so that $\rho$ is away from 1, then the error remains bounded as $\alpha\uparrow 1$ while both the approximate value $\hV^{\alpha}$ and the true value grows like $1/(1-\alpha)$ as  $\alpha\uparrow 1$.

Let $\mu\downarrow \frac{1}{2}$ and $\lambda=1-\mu$, to put the queue in heavy-traffic as before. Taking, as above, $\alpha_{\rho}=1-(1-\rho)^2,$ we have $\rho\uparrow 1$,
$$|V^{\alpha_{\rho}}(x)-\hV^{\alpha_{\rho}}(x)|\leq \frac{\Gamma}{(1-\rho)^2}.$$

Because $V^{\alpha_{\rho}}(x)\geq \frac{1}{(1-\rho)^3}$ for $x\geq \frac{1}{1-\rho}$, it follows that  $$|V^{\alpha_{\rho}}(x)-\hV^{\alpha_{\rho}}(x)|\leq \Gamma(1-\rho)V^{\alpha(\rho)}(x), \mbox{ for all } x\geq \frac{1}{1-\rho}.$$ This establishes the asymptotic correctness of a Brownian approximation by means of Tayloring rather than by those of weak convergence. In contrast to Brownian approximations, {\em Tayloring is a purely analytical device. }

Finally, this is an opportunity to revisit the contribution of ``corners'' to the bound in Theorem \ref{thm:explicit}. In a queueing network where all stations operate at utilization of $1-\sqrt{1-\alpha}$ as above, the fraction of time spent near corners will be of the order of (or smaller than) $(1-\alpha)$, which will multiply $\Theta_M$ in the bound.

\vspace{0.5cm}

\section{Final comments\label{sec:conclusions}}

In this paper, we have introduced Tayloring as a rigorous framework for value function approximation. Applied to a controlled chain in discrete time and space, we derive bounds grounded in PDE theory and propose a solution algorithm with performance guarantees. This paper is a first, and by no means last, step. Much remains open in terms of the scope---continuous time and space, finite and long-run average problems---and various algorithmic aspects. Below is a short informal discussion of these directions.\vspace*{0.2cm}

\paragraph{Continuous-time.} Consider the continuous-time version of the $M/M/1$ queue---Poisson arrivals with rate $\lambda$ and a single server with service times that are exponential with (controlled) parameter $u(x)$. Given holding and service rate cost $r(x,u)$, consider the problem $$\min_U\Ex_x\left[\int_0^{\infty}e^{-(1-\alpha)t}r(X(t),U(t))dt\right].$$ The Bellman equation is given by
$$0=\min_{u\geq 0}\left\{r_u(x)+\lambda (V(x+1)-V(x))+u\1\{x>0\}(V(x-1)-V(x))-(1-\alpha) V(x)\right\}.$$
Second-order Tayloring leads to the TCP
$$0=\min_{u\geq 0}\left\{ r_u(x)+(\lambda-u)V'(x)+\frac{1}{2}(\lambda+u)V''(x)-(1-\alpha) V(x)\right\}, ~x>0,$$ and we use the boundary condition $V'(0)=0$.
Generally, for a continuous time chain on $\mathbb{Z}_+^d$ with transition-rate matrix $q_u(x,y)$ at state $x$, we have
$$(\mu_{u})_i(x)=\sum_{y}q_u(x,y)(y_i-x_i),
\mbox{ and } (\sigma_{u})_{ij}(x)=\sum_{y}q_u(x,y)(y_i-x_i)(y_j-x_j).$$
It seems fair to conjecture that similar results to Theorems \ref{thm:diffusion} and \ref{thm:explicit} can be derived in this context in reasonable generality, noting that in contrast to the transition-probability matrix $P^u(x,y)$, the transition-rate matrix is, in a variety of practical settings, unbounded. \vspace*{0.2cm}

\paragraph{Continuous state space.} For continuous state space processes, we conjecture that Theorems \ref{thm:diffusion} and \ref{thm:explicit} still work with modifications to $\mu_u$ and $\sigma_u$: for all $x\in\mathbb{R}_+^d$,
\begin{align} (\mu_{u})_i(x)&=\Ex_x^{u}[(X_{1})_i-x_i]=\int_{\mathbb{R}_+} P^{u}(x,dy)(y_i-x_i),~i=1,\ldots,d,\\
	(\sigma^2_{u})_{ij}(x)&=\Ex_x^u[((X_1)_i-x_i)((X_1)_j-x_j)] =\int_{\mathbb{R}^2_+}P^{u}(x, dy)(y_i-x_i)(y_j-x_j),~i,j=1,\ldots,d.\label{eq:mudefin}
\end{align}

The continuous state space may seem to simplify things, insofar as there is no need to extend $\mu_u$ and $\sigma_u$ or to smooth the state space. This simplification, however, implies also losing the freedom, for example, to extend the $\mu$ and $\sigma^2$ in a Lipschitz continuous way or smooth the corners of the state space. It is these freedom that facilitated, in this paper, the application of the PDE theory of classical solutions in smooth domains.\vspace*{0.2cm}

\paragraph{Finite-horizon problems.} Consider a discrete time and space dynamic program on a finite horizon of length $T$. Let $V_t(x)$ be the value with $t$ steps to go and starting at state $x$. The Bellman equation is given by
\begin{align*} V_t(x)  = \max_{u}\left\{  r_{u}(x)+\sum_{y}P_{x,y}^{u,t}V_{t-1}(y)\right\},\end{align*} which we can rewrite as
\begin{align*} 0  =& \max_u \left\{ r_{u}(x)+\sum_{y}P_{x,y}^{u,t} (V_{t}(x)-V_{t}(y))+V_{t-1}(x)-V_t(x)\right.\\
	& \left. -\sum_{y}P_{x,y}^{u,t} [V_{t-1}(x)-V_{t-1} (y)+V_t(x)-V_t(y)]\right\}.
\end{align*}

Let $$\mu_{u,t}(x)=\sum_{y}P_{x,y}^{u,t}(y-x),\mbox{ and } \sigma^2_{u,t}(x)=\sum_{y}P_{x,y}^{u,t}(y-x)^2.$$ Taking a second-order expansion in $x$ and a first-order expansion in $t$, we arrive at the equation
\begin{align*} 0 =\max_{u}\left\{ r_{u}(x)+\mu_{u,t}(x)\frac{\partial}{\partial x}V_t(x)+\frac{1}{2}\sigma_{u,t}^2(x)\frac{\partial^2}{\partial x^2}V_t(x)-\frac{\partial}{\partial t}V_t(x)\right\}.\end{align*}
We drop any consideration of boundary condition from this informal outline. One suspects that the errors will depend in the finite horizon setting on the second derivative in $t$ and the third derivative in the state $x$, in addition to the cross-derivative in $x$ and $t$ that arises from the term $\sum_{y}P_{x,y}^{u,t} [V_{t-1}(x)-V_{t-1} (y)+V_t(x)-V_t(y)]$. While the connection between the original chain and the Taylored equation seems a straightforward extension of what we have done in this paper, it is not clear that the K-D approach can be used similarly, as it is built for infinite horizon problems. This remains an interesting direction for future exploration.\vspace*{0.2cm}

\paragraph{Computation with non-uniform grids.} In our computational examples we do not rationalize the choice of the coarseness level $h$ and, once $h$ is chosen, we use it uniformly in the state space. Our bounds, however, might suggest a direction for improvement.

Finite difference methods for PDEs use finer grids in regions where large gradients are expected and coarser grids where the function is relatively ``flat.'' This brings computational efficiency at little cost to accuracy. A similar logic applies to TAPI as is nicely captured in the bound \eqref{eq:KDbound}, which we rewrite here as
\begin{align}
	\Ex_x^{\hU_*}\left[\sum_{t=0}^{\infty}\bar{\alpha}_h^t h^{2+\beta}[D^2\hV_{\ast}]_{\beta,X_t^h\pm h}^*\right]+\Ex_x^{U^h_\ast}\left[\sum_{t=0}^{\infty} \bar{\alpha}_h^t h^{2+\beta}[D^2\hV_{\ast}]_{\beta,X_t^h\pm h}^*\right].\label{eq:fderror}
\end{align}

The error is dependent on the interaction of the step size $h$ with the supremum of $[D^2\hV_{\ast}]$ over a  neighborhood of $x$. It makes sense to use large ``boxes'' where $[D^2\hV_{\ast}]$ is relatively flat and vice versa. Ad-hoc knowledge of the problem can help here. In the inventory problem of \S \ref{sec:examples}, for example, $\hV$ is approximately linear far from the origin, suggesting that one could use a large $h$ in that part of the state space. We can use such understanding to build a (computationally beneficial) TCP-equivalent chain with non-uniform jump sizes over the state space.

We leave full exploration of these computational aspects for future work.
\bibliographystyle{ormsv080}
\bibliography{central}
 \section*{Appendix}

Below we state a lemma that formalizes the connection between the oblique-derivative TCP and the second-order TCP in \eqref{eq:TCP0}. Recall that $f_{\mu}(x,u),f_{r}(x,u),f_{\sigma}(x,u)$ are assumed to be defined for all $x\in \mathbb{R}_+^d$ and $u\in\mathbb{D}$.

\begin{lemma} Suppose that Assumption \ref{asum:1} holds, that $f_r(x,u), f_{\sigma}(x,u)$ are continuous in $x\in\mathbb{R}_+^d$ for each $u$ and that $\hV_*\in \mathcal{C}^2(\mathbb{R}_+^d)$ solves the OD-boundary TCP
	\begin{align} &  0=\max_{u \in\mathcal{U}(x)}\{ r_{u}(x)+\alpha\mathcal{L}_{u}V(x)-(1-\alpha)V(x)\},~ x\in \mathbb{R}_{++}^d, \label{eq:OD1} \\
	&	0 = \eta(x)'DV(x),~x\in \partial\mathbb{R}_+^d,
	\label{eq:OD2} \end{align}
	where, for $x\in\mathcal{B}_i$, $\eta(x)$ is proportional to $\iota(x,u)$ given by the limit (assumed to exist) $$  \iota_i(x,u)=\lim_{x_i\downarrow 0}\left(f_{\mu}((x_{-i},x_i),u)-f_{\mu}((x_{-i},0),u)\right);$$ $\eta_i(x)=0$ for $x\notin\mathcal{B}_i$.
	Suppose, in addition, that the maximizer $\hU_*$ in \eqref{eq:OD1} is continuous at the boundary, i.e., $\hU_*(x^n)\rightarrow \hU_*(x)$ for {\em any} sequence $\{x^n\}$ such that $x^n\in\mathbb{R}_{++}^d$ and $x^n\rightarrow x\in\partial \mathbb{R}_+^d$. Then, $(\hU_*,\hV_*)$ also solves the 2nd-order TCP \eqref{eq:TCP0}.\label{lem:oblique} \end{lemma}

\bProof Take $x\in\mathcal{B}$ and a sequence $\{x^n\}$ as in the statement of the lemma. Then, by definition of the TCP, we have, for each $n$,
\begin{align*} 0&=r_{\hU_*}(x^n)+\alpha\mu_{\hU_*}(x^n)\cdot D\hV_*(x)+\frac{\alpha}{2}trace(\sigma^2_{\hU_*}(x^n)\cdot D^2\hV_*(x^n))-(1-\alpha)\hV_*(x^n)\\
&=r_{\hU_*}(x^n)+\alpha\mu_{\hU_*}(x)\cdot D\hV_*(x^n)+
\alpha(\mu_{\hU_*}(x^n)-\mu_{\hU_*}(x))\cdot D\hV_*(x^n)\\
&\quad\quad+\frac{\alpha}{2}trace(\sigma^2_{\hU_*}(x^n)\cdot D^2\hV_*(x^n))-(1-\alpha)\hV_*(x^n),
\end{align*} where in the second line we added and subtracted $\alpha\mu_{\hU_*}(x)\cdot D\hV_*(x)$. By assumption $r_{\hU_*}(x^n),\sigma^2_{\hU_*}(x^n),D\hV_*(x^n),D^2\hV_*(x^n)\rightarrow
r_{\hU_*}(x),\sigma^2_{\hU_*}(x),D\hV_*(x),D^2\hV_*(x)$. Also, $(\mu_{\hU_*}(x^n)-\mu_{\hU_*}(x))\cdot D\hV_*(x)\rightarrow \eta_{\hU_*}(x)\cdot D\hV_*(x)=0$ so that
\begin{align*} (\mu_{\hU_*}(x^n)-\mu_{\hU_*}(x))\cdot D\hV_*(x^n)=
(\mu_{\hU_*}(x^n)-\mu_{\hU_*}(x))\cdot D\hV_*(x)+
(\mu_{\hU_*}(x^n)-\mu_{\hU_*}(x))\cdot (D\hV_*(x^n)-D\hV_*(x))\rightarrow 0.\end{align*} In sum, we have $$r_{\hU_*}(x)+\alpha \mathcal{L}_{\hU_*}\hV_*(x)-(1-\alpha)\hV_*(x)=\lim_{n\tinf}
r_{\hU_*}(x^n)+\alpha \mathcal{L}_{\hU_*}\hV_*(x^n)-(1-\alpha)\hV_*(x^n)=0,$$ as required. \eProof \vspace{0.5cm}

\paragraph{Proof of Theorem \ref{thm:diffusion}:}
For any stationary policy $U\in\mathbb{U}$ and any function $\Phi:\mathbb{Z}_+^d\to \mathbb{R}_+$, by standard arguments,
\begin{equation*}
\begin{split}
\mathbb{E}^{U}_x[\alpha^n \Phi(X_n)]
=&\Phi(x)+\mathbb{E}^{U}_x\left[\sum_{t=1}^{n}\alpha^{t-1}(\alpha \Phi(X_t)-\Phi(X_{t-1}))\right] \\
=&\Phi(x)+\mathbb{E}^{U}_x\left[\sum_{t=1}^{n}\alpha^{t-1}(\alpha P^U-I)\Phi(X_{t-1})
\right]\\
=&\Phi(x)+\mathbb{E}^{U}_x\left[\sum_{t=1}^{n}\alpha^{t-1}\left(\alpha (P^U-I)\Phi(X_{t-1})
-(1-\alpha)\Phi(X_{t-1})\right)\right].
\end{split}
\end{equation*}
As a result, if $\Phi$ and $U$ are such that, for all $x\in\mathbb{R}_+^d$, $\Ex_x^U[\alpha^n \Phi(X_n)]\rightarrow 0$ as $n\tinf$, then \begin{align*} 0
&=\Phi(x)+\mathbb{E}^{U}_x\left[\sum_{t=1}^{\infty}\alpha^{t-1}\left(\alpha (P^U-I)\Phi(X_{t-1})
-(1-\alpha)\Phi(X_{t-1})\right)\right]\\& =
\Phi(x)+\mathbb{E}^{U}_x\left[\sum_{t=0}^{\infty}\alpha^{t}\left(\alpha (P^U-I)\Phi(X_{t})
-(1-\alpha)\Phi(X_{t})\right)\right]. \end{align*}

Suppose that $\Phi$ is twice continuously differentiable,
 $f:\mathbb{Z}_+^d\to \mathbb{R}$ and $U\in\mathbb{U}$ are such that $\alpha \mathcal{L}_U \Phi(x)-(1-\alpha)\Phi(x)=-f(x)$ for all $x\in\mathbb{R}_+^d$, then
\begin{align*} 0&=\Phi(x)-\Ex_x^U\left[\sum_{t=0}^{\infty}\alpha^t f(X_t)\right]+\Ex_x^U\left[\sum_{t=0}^{\infty}\alpha^t \mathcal{A}_U[\Phi](X_t)\right],\end{align*}
where $\mathcal{A}_U[\Phi](x)=\alpha (P^U-I)\Phi(x)
-\alpha \mathcal{L}_U \Phi(x).$ Similarly, if $\alpha \mathcal{L}_U \Phi(x)-(1-\alpha)\Phi(x)\leq -f(x)$, then
\begin{align} 0&\leq \Phi(x)-\Ex_x^U\left[\sum_{t=0}^{\infty}\alpha^t f(X_t)\right]+\Ex_x^U\left[\sum_{t=0}^{\infty}\alpha^t \mathcal{A}_U[\Phi](X_t)\right].\label{eq:ito2}\end{align}
Taking $\Phi(x)=\hV_{\ast}^{\alpha}(x)$, $U=\hU_*$, and $f(x)=r_{\hU_*}(x)$, we then have
\begin{align*} 0&=\hV_{\ast}^{\alpha}(x)-\Ex_x^{\hU_*}\left[\sum_{t=0}^{\infty}\alpha^t r_{\hU_*}(X_t)\right]+\Ex_x^u\left[\sum_{t=0}^{\infty}\alpha^t \mathcal{A}_{\hU_*}[\hV_*^{\alpha}](X_t)\right]\\
& = \hV_{\ast}^{\alpha}(x)-V_{\hU_*}^{\alpha}(x) +\Ex_x^u\left[\sum_{t=0}^{\infty}\alpha^t \mathcal{A}_{\hU_*}[\hV_*^{\alpha}](X_t)\right].\end{align*}  Thus,
$$|\hV_{\ast}^{\alpha}(x)-V_{\hU_*}^{\alpha}(x)|\leq \Ex_x^{\hU_*}\left[\sum_{t=0}^{\infty}\alpha^t |\mathcal{A}_{\hU_*}[\hV_*^{\alpha}](X_t)|\right].$$ To obtain the bound under $U_*$, notice that  for all $x\in\mathbb{R}_+^d$ and any feasible action $u\in\mathcal{U}(x)$, then by the definition of $\hV_{*}^{\alpha}$,
$r_{u}(x)+\alpha \mathcal{L}_{u} \hV_{\ast}^{\alpha}(x)-(1-\alpha)\hV_{*}^{\alpha}(x)\leq 0$. Using \eqref{eq:ito2}, we have
$$0\leq \hV_{\ast}^{\alpha}(x)-\Ex_x^{U_*}\left[\sum_{t=0}^{\infty}\alpha^t r_{U_*}(X_t)\right]+\Ex_x^{U_*}\left[\sum_{t=0}^{\infty}\alpha^t \mathcal{A}_{U_*}[\hV_*^{\alpha}](X_t)\right].$$ Consequently (recalling $V_{\ast}^{\alpha}(x)=V_{U_*}^{\alpha}(x)$),
$$\hV_{\ast}^{\alpha}(x)\geq V_{\ast}^{\alpha}(x) -\Ex_x^{U_*}\left[\sum_{t=0}^{\infty}\alpha^t |\mathcal{A}_{U^*}[\hV_*^{\alpha}](X_t)|\right].$$

We conclude that for all $x\in\mathbb{R}_+^d$,
$$\hV_{*}^{\alpha}(x)-\Ex_x^{\hU_*}\left[\sum_{t=0}^{\infty}\alpha^t |\mathcal{A}_{\hU_*}[\hV_*^{\alpha}](X_t)|\right]\leq V_{\hU_*}^{\alpha}(x)\leq V_{\ast}^{\alpha}(x)\leq \hV_{*}^{\alpha}(x)+\Ex_x^{U_*}\left[\sum_{t=0}^{\infty}\alpha^t |\mathcal{A}_{U^*}[\hV_*^{\alpha}](X_t)|\right],$$
from which it follows both that \begin{equation}  |\hV_{*}^{\alpha}(x)-V_{\ast}^{\alpha}(x)|\leq \Ex_x^{U_*}\left[\sum_{t=0}^{\infty}\alpha^t |\mathcal{A}_{U^*}[\hV_*^{\alpha}](X_t)|\right]+\Ex_x^{\hU_*}\left[\sum_{t=0}^{\infty}\alpha^t |\mathcal{A}_{\hU_*}[\hV_*^{\alpha}](X_t)|\right],\label{eq:b1} \end{equation}
and
\begin{equation} |V_{\hU_*}^{\alpha}(x)-V_{\ast}^{\alpha}(x)|\leq \Ex_x^{U^*}\left[\sum_{t=0}^{\infty}\alpha^t |\mathcal{A}_{U^*}[\hV_*^{\alpha}](X_t)|\right]+\Ex_x^{\hU_*}\left[\sum_{t=0}^{\infty}\alpha^t |\mathcal{A}_{\hU_*}[\hV_*^{\alpha}](X_t)|\right].\label{eq:b2} \end{equation}

Thus, to conclude the proof, it only remains to bound the terms on the right-hand side of \eqref{eq:b1} and \eqref{eq:b2}. By Taylor's Theorem, for each pair $x, y\in\mathbb{R}_+^d$, there is a $\zeta_{x,y}$ connecting $x$ and $y$ such that
\begin{align*}
(P^u-I)f(x)=\sum_{y}P^u(x,y)\left(Df(x)'(y-x)+\frac{1}{2}(y-x)'D^2f(\zeta_{x,y})(y-x)\right).
\end{align*}

Then,
\begin{align*}
|\mathcal{A}_u[f](x)| \leq |(P^u-I-\mL_{u})f(x)|=&\frac{1}{2}\left|\sum_{y}P^{u}(x,y)(y-x)'(D^{2}f(\zeta_{x,y})-D^{2}f(x))(y-x)\right|\\
\leq &\frac{1}{2}\sum_{y}P^{u}(x,y)|y-x|^2[f]_{\beta,x\pm \frakj_u}^*|y-x|^{\beta}\leq \frakj_u^{2+\beta}[f]_{\beta,x\pm \frakj_u}^*.
\end{align*}

The first inequality follows because $\alpha\leq 1$. The second inequality holds because $\zeta_{x, y}\in[x, y]$ and then $|\zeta_{x, y}-x|\leq \frakj_{u}$.
Thus, for $u\in \{U^*,\hU_*\}$ we have
    \begin{equation*}
\begin{split}
\Ex_x^{u}\left[\sum_{t=0}^{\infty}\alpha^t |\mathcal{A}_{u}[\hV_{*}^{\alpha}](X_t)|\right] \leq \frakj_{u}^{2+\beta} \Ex^u_x\left[\sum_{t=0}^{\infty}\alpha^t [\hV_{\ast}^{\alpha}]_{\beta,X_t\pm \frakj_u}^*\right],
\end{split}
\end{equation*}
as required.

For Remark \ref{rem:unboundedjumps}, notice that the last steps can be generalized as follows:
\begin{align*}
|\mathcal{A}_u[f](x)| \leq |(P^u-I-\mL_{u})f(x)|=&\frac{1}{2}\left|\sum_{y}P^{u}(x,y)(y-x)'(D^{2}f(\zeta_{x,y})-D^{2}f(x))(y-x)\right|\\
\leq &\frac{1}{2}\sum_{y}P^{u}(x,y)|y-x|^2[f]_{\beta,x\pm \Delta_x}^*|y-x|^{\beta},
\end{align*}  where $\Delta_x$ is the random variable representing the jump, i.e., $\Delta_x=X_1-x$. Thus,
\begin{align*}
\Ex^u_x[\mathcal{A}_u[f](x)|] \leq \Ex^u_x[|X_1-x|^{2+\beta}[f]_{\beta,x\pm \Delta_x}^*],\end{align*}
and  \begin{equation}
\begin{split}
\Ex_x^{u}\left[\sum_{t=0}^{\infty}\alpha^t |\mathcal{A}_{u}[\hV_{*}^{\alpha}](X_t)|\right] \leq \Ex_x^u\left[\sum_{t=0}^{\infty}\alpha^t \Ex^u_{X_t}[|\Delta X_t|^{2+\beta}[f]_{\beta,X_t\pm \Delta_{X_t}}^*]\right].
\end{split}\label{eq:unboundedjumps}
\end{equation}

\eProof\vspace*{0.5cm}

\paragraph{Proof of Lemma \ref{lem:PDEbounds}.} The theory of gradient estimates for PDEs is somewhat involved. The existing bounds identify the primitives on which the constants depend but do not always make explicit the nature of this dependence. Here (in contrast, say, to \citeapp{Gurvich_AAP2014}), the only parameter being scaled is the discount factor and, fortunately, the constants in the PDE bounds do not depend on it. We rely in what follows on \citeapp{safonov1,safonov,safonov1995oblique}, where the dependence of the bounds on the domain $\Omega_M$ and on the reward function $r_u$ is relatively explicit.

It will help to have labels for the key conditions: \begin{align}
\label{cond:ellip}
& \exists \lambda>0 \mbox{ s.t. } \lambda|\xi|^2\leq \sum_{ij}\xi_{i}\xi_j {(\sigma^2_{u})}_{ij}(x)\leq \lambda^{-1}|\xi|^2, \mbox{ for all } \xi\in\mathbb{R}^d,~x\in\mathbb{R}_+^d,~u\in\mathbb{D},
 \tag{A0} \\
\label{cond:musig}
&\exists L \mbox{ s.t. } |\sigma^2_u|_{\Omega_M}^*,|\mu_u|_{\Omega_M}^*,[\mu_u]_{1,\Omega_M}^*,[\sigma^2_u]_{1,\Omega_M}^*\leq L,~~ u\in \mathbb{D},\tag{A1}\\
& \sup_{u\in \mathbb{D}}|r_u|_{\overline{\Omega_M}}^*+\sup_{u\in \mathbb{D}}[r_u]_{1,\overline{\Omega_M}}^*<\infty, \label{cond:goodreward} \tag{A2}
\\ \label{cond:OD}
&\exists \nu_0>0 \mbox{ s.t. } \eta(x)\cdot \theta(x) \geq \nu_0|\eta(x)|, ~x\in \partial \Omega_M, ~|D\eta|_{\partial\Omega_M}^*+[D\eta]_{1,\partial\Omega_M}^*\leq L,
\tag{A3}
\end{align}
where $\theta(x)$ is outward normal to $\partial{\Omega_M}$ at $x$. Notice that if a function is Lipschitz continuous (with constant $L$) and bounded by $L$ over a domain $\Omega$ then it is H\"{o}lder continuous for any exponent $\beta\in (0,1)$ with a suitable constant.

For completeness we repeat a standard definition. A bounded domain $\Omega$ in $\mathbb{R}_+^d$ and its boundary are of class $\mathcal{C}^{2,\beta}$ for $\beta\in [0,1]$ if for each point in $x_0$ there is a ball $B=B(x_0)$ and a one-to-one mapping $\Psi$ of $B$ onto $D\subset \mathbb{R}^d$ such that (i) $\Psi(B\cap \Omega)\subset \mathbb{R}_+^d$, (ii) $\Psi(B\cap \partial\Omega)\subset \partial \mathbb{R}_+^n$, and (iii) $\Psi\in \mathcal{C}^{2,\beta}(B)$, $\Psi^{-1}\in \mathcal{C}^{2,\beta}(D)$; see, e.g., \citeapp[page 94]{gilbarg2015elliptic}. When the estimate in Theorem \ref{thm:safonov} below is said to depend on $\partial \Omega_M$, the dependence is through the bounds on the $\mathcal{C}^{2,\beta}(\partial\Omega_M)$ norms of $\Psi$, i.e., on its first derivative, second derivative, and $\beta-$H\"{o}lder continuity constant. Imperative for us is that $\Psi$ for $\Omega_M$, as constructed, does not change with $\alpha$ or, indeed, with $M$.

The existence of a solution  $\hV\in \mathcal{C}^{2,\beta}(\Omega_M)$ for either the OD- or FOT-boundary TCP---under the assumptions on $\mu_u$, $\sigma^2_u$, and $r_u$ in Lemma \ref{lem:PDEbounds}---follows from general non-linear PDE results; see \citeapp[Theorem 7.11]{lieberman1986nonlinear}; see also the discussion in the supplementary remarks at the bottom of page 544 there. For the FOT-boundary, see the discussion at the top of that page. The H\"{o}lder constant $\beta$ in these existence results does not depend on our discount factor $\alpha$.

The following states that  any $V\in \mathcal{C}^{2,\beta}$ that satisfies the TCP in the interior of $\Omega_M$  must satisfy certain gradient properties. Recall that
$$\Omega_{--}^{\varrho(\cdot)}:=\{x\in \Omega_M: d(x,\mathcal{B}_i)> \varrho(x), \mbox{ for all $i$}\}.$$

\begin{theorem}[\citeapp{safonov1}, Theorem 1.1] Suppose that  \eqref{cond:ellip}-\eqref{cond:goodreward} hold. Then $V\in\mathcal{C}^{2,\beta}(\overline{\Omega_M})$ that satisfies $$0=\max_{u \in \mathbb{D}}\{r_u(x)+\alpha\mathcal{L}_uV(x)-(1-\alpha)V(x)\}, ~x\in \Omega_M,$$ also satisfies for all $x\in\Omega_{--}^{\varrho(\cdot)}$,  $$[D^2V]^*_{\beta,x\pm \frac{\varrho(x)}{2}}\leq \Gamma(L,\lambda,d,\beta,\partial\Omega_M)\left(\frac{|V|_{\mathcal{B}^{\varrho(x)}(x)}^*}{\varrho^{2+\beta}(x)}+
	\max_{u\in \mathbb{D}}[r_u]_{\beta,\mathcal{B}^{\varrho(x)}(x)}^*\right).$$
	With the OD boundary condition, this estimate holds also for $x$ in the larger set $\Omega_-^{\varrho(\cdot)}$.		
\end{theorem}
Theorem 1.1 of \citeapp{safonov1} is stated for the Dirichlet problem (i.e., with a boundary condition $V=\phi$ on $\partial\Omega_M$) but the reader should note that, once existence of a solution $V$ to the oblique-derivative (or first-order Tayloring) is established, this solution in particular solves a Dirichlet problem with the boundary function $\phi$ being given by the values of $V$; see also Theorem 5.1 in  \citeapp{safonov}. Safonov's Theorem 1.1 is stated for points that are suitably ``far'' from all boundaries. It is \cite[Theorem 8.2]{safonov} that allows us---for the OD boundary---to extend the bound from $\Omega_{--}^{\varrho(\cdot)}$ to the larger set $\Omega_{-}^{\varrho(\cdot)}$.

Finally, for the reader unfamiliar with the PDE notation used in \citeapp{safonov1}, notice that the interior norms are introduced in the Definition on page 599 there. Remark 1.1 explains how the Bellman equation belongs to the class covered by Theorem 1.1 and, in condition (F4) of Definition 1.1 there it is immediate that taking $\Omega\subseteq \Omega_M$, $K_1'=\sup_{u\in \mathbb{D}}[r_u]_{\beta,\Omega}^*$. Similarly, $L$ in the condition \eqref{cond:musig} replaces $K$ and $K_1$ in Safonov's statement and $\lambda$ replaces his $\nu$.

To treat points $x\in \overline{\Omega_M}\backslash \Omega_{-}^{\varrho(\cdot)}$, we use a global estimate relying again on \citeapp{safonov,safonov1}; see also \citeapp{lions1986linear} and \citeapp{lieberman2013oblique}. In applying Theorem 3.3 of \citeapp{safonov1995oblique}, notice that his function $g$ is identically $0$ in our case.

\begin{theorem}[\citeapp{safonov}, Theorem 8.3] Suppose that $V\in \mathcal{C}^{2,\beta}(\overline{\Omega_M})$ solves the TCP with either the OD- or the control-independent FOT-boundary conditions and that \eqref{cond:ellip}-\eqref{cond:OD} hold. Then, we have the global estimate  $$[D^2 V]_{\beta,\overline{\Omega_M}}^*\leq \Gamma(L,\lambda,d,\nu_0,\beta,\partial\Omega_M)\left(|V|_{\overline{\Omega_M}}^*+\max_{u \in \mathbb{D}}[r_u]_{\beta,\overline{\Omega_M}}^*\right).$$
	\label{thm:safonov} \end{theorem}

Theorem 8.3 of \citeapp{safonov} states that the constant ($\Gamma$ in our notation, $N$ in his) depends on the domain $\Omega_M$. This dependence is, however, only through the smoothness properties of $\partial \Omega_M$. In addition, the requirement that $|b_0|\geq v_0$ in the statement of Theorem 3.3 there is only used for the proof of the existence of a solution, not for the gradient bounds; see also Theorem 3.3 in \citeapp{safonov1995oblique}.
Finally, we notice that  Safonov's result applies only to the {\em control-independent} FOT boundary condition. While the general bounds in \citeapp{lieberman1986nonlinear} apply also to the general FOT boundary condition, they are less explicit in their dependence on, for example, the size of the state space $\Omega_M$.

These two theorems combined provide the supporting details for Lemma \ref{lem:PDEbounds}. Before turning to the proof of Theorem \ref{thm:explicit}, we state the analogue of Lemma \ref{lem:PDEbounds} for the control-independent FOT case. For the following we let $\eta(x)=\alpha \mu(x)$ (recall that with control-independent $\mu_u(x)\equiv \mu(x)$ for $x\in \partial \Omega_M$). With this notation, the only difference is the restriction of \eqref{eq:lipbound} to $\Omega_{--}^{\varrho(\cdot)}$.

\begin{lemma}	\label{lem:PDEboundsFOT}
Suppose that Assumption 1 holds. Let the conditions on $\Omega_M$, $\mathcal{U}(x)$, and $\eta$ be as in Lemma \ref{lem:PDEbounds} and suppose that Assumption \ref{asum:1} holds. Then, the TCP with the control-independent FOT boundary condition has a unique solution $\hV_*\in \mathcal{C}^{2,\beta}(\Omega_M)$ for some $\beta\in (0,1)$ (that does not depend on $\alpha$). Moreover, given a function $\varrho:\mathbb{Z}_+^d\to \mathbb{R}_{++}$, \eqref{eq:lipbound} holds with $\Omega_{--}^{\varrho(\cdot)}$ replacing  	$\Omega_{-}^{\varrho(\cdot)}$ and the global bound remains unchanged. \end{lemma}

Notice that Theorem \ref{thm:explicit} and Corollaries \ref{cor:alphaopt} and \ref{cor:orderopt} then hold with the obvious replacement of $\mathcal{T}_{CO}$ with
$$\mathcal{T}_{CO-}=\{t\geq 0: X_t\notin \Omega_{--}^{\varrho(\cdot)}\}.$$ \eProof \vspace{0.5cm}

\paragraph{Proof of Theorem \ref{thm:explicit}.} Set $\varrho(x)=(1+|x|)+\frac{1}{(1-\alpha)^{\delta}}+2\frakj$ with $\delta= m/(k+2)$. Then, per the assumptions of the theorem,
$$|\hV|_{x\pm \frac{\varrho(x)}{2}}^*\leq \Gamma\left(\frac{|x|^{k}}{1-\alpha}+\frac{1}{(1-\alpha)^{\delta k}}+\frac{1}{(1-\alpha)^m}\right),$$ and
$$[r_{\hU_*}]_{\beta,x\pm \frac{\varrho(x)}{2}}^*\leq \Gamma\left(1+|x|^{k-\beta}+\frac{1}{(1-\alpha)^{\delta(k-\beta)}}\right).$$
Consequently, per Lemma \ref{lem:PDEbounds}, we have for all $x\in\Omega_-^{\varrho(\cdot)}$,
\begin{align*}
[D^2\hV]^*_{\beta,x\pm \frac{\varrho(x)}{2}}& \leq \Gamma \left( \frac{|\hV|_{\mathcal{B}^{\varrho(x)}(x)}^*}{((1+|x|)+\frac{1}{(1-\alpha)^{\delta }})^{2+\beta}}+[r_{\hU_*}]_{\beta,x\pm \varrho(x)}^*\right)
\\
& \leq \Gamma\left( \frac{(1+|x|)^{k-(2+\beta)}}{1-\alpha}+\frac{1}{(1-\alpha)^{\delta(k-(2+\beta))}}+\frac{1}{(1-\alpha)^{m-(2+\beta)\delta }}+(1+|x|)^{k-\beta}+\frac{1}{(1-\alpha)^{\delta(k-\beta)}}\right).
\end{align*}

With $\delta= m/(k+2)$, we have $\delta(k-\beta)\vee (m-(2+\beta)\delta) = \frac{m(k-\beta)}{k+2}$ so that
$$[D^2\hV]^*_{\beta,x\pm \frac{\varrho(x)}{2}}\leq \Gamma\left(\frac{(1+|x|)^{k-(2+\beta)}}{1-\alpha}+(1+|x|)^{k-\beta}+\frac{1}{(1-\alpha)^\frac{m(k-\beta)}{k+2}}\right).$$
For $x:|x|\geq \frac{1}{\sqrt{1-\alpha}}$,
$\frac{(1+|x|)^{k-(2+\beta)}}{1-\alpha}\leq \Gamma (1+|x|)^{k-\beta}$ and for $x:|x|\leq \frac{1}{\sqrt{1-\alpha}}$,
$\frac{(1+|x|)^{k-(2+\beta)}}{1-\alpha}\leq \Gamma\frac{1}{(1-\alpha)^{\frac{1}{2}(k-\beta)}}$. Thus,
$\frac{(1+|x|)^{k-(2+\beta)}}{1-\alpha}\leq \Gamma \left((1+|x|)^{k-\beta}+\frac{1}{(1-\alpha)^{\frac{1}{2}(k-\beta)}}\right)$
and we arrive at
$$[D^2\hV]^*_{\beta,x\pm  \frac{\varrho(x)}{2}}\leq \Gamma\left((1+|x|)^{k-\beta}+\frac{1}{(1-\alpha)^\frac{m(k-\beta)}{k+2}}
+\frac{1}{(1-\alpha)^{\frac{1}{2}(k-\beta)}}\right), ~ x\in\Omega_-^{\varrho(\cdot)}.$$

Because $|x+\bar{\frakj}|^{k-1}\leq \Gamma(|x|^k+\frakj^k)$,
\begin{align*}
[D^2\hV]_{\beta, x\pm \bar{\frakj}}^*\leq
\sup_{y:|y-x|\leq \bar{\frakj}}[D^2\hV]_{\beta,y\pm \frac{\varrho(y)}{2}}^*
\leq \Gamma\left((1+|x|)^{k-\beta}+\frac{1}{(1-\alpha)^\frac{m(k-\beta)}{k+2}}
+\frac{1}{(1-\alpha)^{\frac{1}{2}(k-\beta)}}\right),
\end{align*}
for a redefined constant $\Gamma$ (that depends on $\bar{\frakj}$).
Notice that by the definition of $\Omega_-^{\varrho(\cdot)}$, $y$ in $x\pm \bar{\frakj}$ is such that $d(y,\mathcal{B}_i)\leq \varrho(y)$ for at most one $i$, so that we can use the bounds.

We conclude that for any stationary control, $$\Ex_x^{U}\left[\sum_{t\notin\mathcal{T}_{CO}}\alpha^t[D^2\hV]_{\beta,X_t\pm \bar{\frakj}}^*\right]\leq \Gamma\left(V_U^{\alpha}[f_{k-\beta}](x)+\frac{1}{(1-\alpha)^{\frac{m(k-\beta)}{k+2}+1}}+\frac{1}{(1-\alpha)^{\frac{1}{2}(k+2-\beta)}}\right),
$$ where, notice, the constant $\Gamma$, depends on $\bar{\frakj}$. \eProof \vspace*{0.5cm}

\paragraph{Proof of Corollary \ref{cor:alphaopt}.} Because the Markov chain has jumps bounded (uniformly in the control) by $\bar{\frakj}$, we claim that for all $x\geq (1-\alpha)^{-1}$ and all stationary control $U$,
\be
V_U^{\alpha}[f_{k-\beta}](x)\leq \Gamma(1-\alpha)^{\beta}V_U^{\alpha}[f_k](x), \mbox{ and } V_U^{\alpha}[f_k](x)\geq \gamma \frac{1}{(1-\alpha)^{k+1}}.
\label{eq:interim}\ee
This would imply that for all such $x$, because $m\leq k+1$ by assumption,
$(1-\alpha)^{\beta}V_U^{\alpha}[f_k](x)\geq  \gamma\frac{1}{(1-\alpha)^{\frac{m(k-\beta)}{k+2}+1}}$.
This concludes the proof of the Corollary under the assumptions that
$\gamma V_{U_*}^{\alpha}[f_k](x)\leq V_{U_*}^{\alpha}(x)$ and
$\gamma V_{\hU_*}^{\alpha}[f_k](x)\leq  V_{\hU_*}^{\alpha}(x)$
for all such $x$.

It remains to prove \eqref{eq:interim}.
Because the drift (in absolute value) is bounded from above by $\bar{\frakj}$ we have $|X_t|\geq |x|/2$ for all $n\leq |x|/(2\bar{\frakj})$, so that
$$V_U^{\alpha}[f_k](x)=\Ex_x^u\left[\sum_{t=0}^{\infty}\alpha^t|X_t|^k\right]\geq \sum_{t=0}^{\lfloor \frac{|x|}{2\bar{\frakj}}\rfloor }\alpha^t (|x|-\bar{\frakj}n)^k\geq \gamma|x|^{k}\frac{1-\alpha^{\frac{|x|}{2\bar{\frakj}}}}{1-\alpha}.$$
There exists $\epsilon<1$ so that for all $x:|x|\geq (1-\alpha)^{-1}$, and all $\alpha\geq 0.5$, $\alpha^{\frac{|x|}{2\bar{\frakj}}}\leq \epsilon<1$.
This result is based on the simple fact that $\alpha^{\frac{1}{1-\alpha}}\rightarrow e^{-1}$ as $\alpha\uparrow 1$. Thus, we have
$$V_U^{\alpha}[f_k](x)\geq \gamma\frac{|x|^k}{1-\alpha},$$ for a redefined constant $\gamma$ and all $x:|x|\geq \frac{1}{1-\alpha}$.
Next note that because $|X_t|^{k-\beta}\leq (1-\alpha)^{\beta}|X_t|^k+\frac{1}{(1-\alpha)^{k-\beta}}\leq (1-\alpha)^{\beta}\left(|X_t|^k+\frac{1}{(1-\alpha)^{k}}\right)$,
\begin{align*}
V^{\alpha}_U[f_{k-\beta}](x)&\leq \Ex^U_x\left[\sum_{t=0}^{\infty}\alpha^t|X_t|^{k-\beta}\right]\leq
(1-\alpha)^{\beta}\left(\Ex^U_x\left[\sum_{t=0}^{\infty}\alpha^t|X_t|^{k}\right]+\frac{1}{(1-\alpha)^{k+1}}\right)\\
\\&=(1-\alpha)^{\beta}\left(V_U^{\alpha}[f_k](x)+\frac{1}{(1-\alpha)^{k+1}}\right) \leq
\Gamma(1-\alpha)^{\beta}V_U^{\alpha}[f_k](x).  \end{align*} The last inequality follows because we have already shown that for all $x:|x|\geq \frac{1}{1-\alpha}$, $V_U^{\alpha}[f_k](x)\geq \gamma|x|^k/(1-\alpha)\geq \frac{1}{(1-\alpha)^{k+1}}$.\eProof \vspace*{0.5cm}

\paragraph{Proof of Corollary \ref{cor:orderopt}.} Let $\tau=\inf\{t\geq 0:|X_t|=0\}$. Because jumps are bounded by $\bar{\frakj}$ uniformly in the control, there exists $\mathfrak{b}$ such that $\tau\geq \mathfrak{b}|x|$, a.s. Also, for all $t\leq \tau$, $|X_t|\geq |x|-\bar{\frakj}t$. Because $\alpha<1$ we have
$$V_U^{\alpha}[f_{k-\beta}](x)=\Ex_x^u\left[\sum_{t=0}^{\infty}\alpha^t |X_t|^{k-\beta}\right]\geq \alpha^{\mathfrak{b}|x|}\sum_{t=0}^{\mathfrak{b}|x|} (|x|-\bar{\frakj}t)^{k-\beta}\geq \gamma\alpha^{\mathfrak{b}|x|}|x|^{k+1-\beta}.$$

There exists $\nu>0$ such that $\gamma\alpha^{\mathfrak{b}|x|}\geq \nu>0$ for all $\alpha\in (0,1)$ and $x:|x|\leq \frac{1}{1-\alpha}$ (this result follows again from the properties of $(1-\epsilon)^{\frac{1}{\epsilon}}\to e^{-1}$ as $\epsilon\downarrow 0$). We thus have for all such $x$ that
$V_U^{\alpha}[f_{k-\beta}](x)\geq \nu|x|^{k+1-\beta}.$
Because $\chi(m,k)=\max\{\frac{m(k-\beta)}{k+2}+1,\frac{1}{2}(k+2-\beta)\}\leq k+1-\beta$
we have, for all $x:\frac{1}{(1-\alpha)^{\zeta(m,k)}}\leq x\leq \frac{1}{1-\alpha}$,
\begin{align} \label{eq:temp1} V_U^{\alpha}[f_{k-\beta}](x)\geq \nu \frac{1}{(1-\alpha)^{\chi(m,k)}}.\end{align}
As in the proof of Corollary \ref{cor:alphaopt}, for all $x:|x|\geq \frac{1}{1-\alpha}$,
\begin{align}  V_U^{\alpha}[f_{k-\beta}](x)\geq \gamma\frac{|x|^{k-\beta}}{1-\alpha}\geq \gamma\frac{1}{(1-\alpha)^{k+1-\beta}}\geq \gamma \frac{1}{(1-\alpha)^{\chi(m,k)}},\end{align}
where we used again the fact $\zeta(m,k)\leq k+1-\beta$. Overall, we have   $$V_U^{\alpha}[f_{k-\beta}](x)\geq \frac{\gamma}{(1-\alpha)^{\chi(m,k)}}, \mbox{ for all } x:|x|\geq \frac{1}{(1-\alpha)^{\zeta(m,k)}}.$$ The corollary now follows directly from equation \eqref{eq:expbound} in Theorem \ref{thm:explicit}. \eProof \vspace*{0.5cm}

\begin{lemma} \label{lem:1dimbound} Let $\mu:\mathbb{R}_+\to \mathbb{R}$ satisfy $|\mu(x)|\leq \Gamma(1+x)$ for all $x\geq 0$ and $\sigma^2:\mathbb{R}_+\to \mathbb{R}_+$ be continuously differentiable and satisfy $0<\lambda \leq \sigma^2\leq \frac{1}{\lambda}$ for some $\lambda>0$. Consider a solution $\hV$ to the ordinary differential equation \begin{align*}
	0&= f(x)+\mu(x)V'(x)+\frac{1}{2}\sigma^2(x) V''(x)-(1-\alpha) V(x),\quad x> 0, \\
	0&=V'(0).
	\end{align*}
	If, for all $x\geq 0$, $\mu(x)\leq 0$ and $f(x)\leq \frac{A}{(1-\alpha)^{\frac{k}{2}}}+Bx^k$, then
	$$\hV(x)\leq \Gamma\left(\frac{x^k}{1-\alpha}+\frac{1}{(1-\alpha)^{\frac{k+2}{2}}}\right).$$
\end{lemma}

The simple proof relies on the connection between the differential equation and a reflected diffusion process. Given such a connection, it is plausible that one can build on the extensive literature on reflected diffusion processes to obtain preliminary bounds on $\hV$.

\bProof Given a Brownian motion $W$, let $\widetilde{X}$ be the unique strong solution to the stochastic differential equation
\begin{align*}
&d\widetilde{X}(t)=\mu(\widetilde{X}(t))dt+\sigma(\widetilde{X}(t))dW(t)+d\widetilde{L}(t), \\&
\widetilde{X}(t)d\widetilde{L}(t)=0,  \\&
\widetilde{X}(t)\geq  0,\quad \mbox{$\widetilde{L}(0)=0$ and $\widetilde{L}$ is nondecreasing.}
\end{align*} See \citeapp{zhang1994strong} for the existence and uniqueness of $\widetilde{X}$. Applying Ito's formula to $V(X(t))$ (see, e.g., \citeapp[Chapter 8]{klebaner2012introduction}) it follows that $V(x)$ has the representation
\begin{equation*}
\hV(x)=\mathbb{E}_x\left[\int_0^{\infty}e^{-(1-\alpha) s}f(\widetilde{X}(s))ds\right],
\end{equation*}

It is a standard result that the increasing process $\widetilde{L}(t)$ is given by the Skorohod map $\widetilde{L}(t)=\sup_{0\leq s\leq t}[x+\int_0^s \mu(\widetilde{X}(z))dz+\int_0^s \sigma(\widetilde{X}(z))dW(z)]^+$ and, in particular,
\begin{equation*}
\begin{split}
0\leq \widetilde{X}(t)\leq & 2\max_{0\leq s\leq t}\left[x+\int_0^s\mu(\widetilde{X}(z))dz+\int_0^s\sigma (\widetilde{X}(z))dW(z)\right]^+\\
\leq &2\max_{0\leq s\leq t}\left|x+\int_0^s\sigma (\widetilde{X}(u))dW(u)\right|,
\end{split}
\end{equation*} where the last inequality follows from our assumption that $\mu\leq 0$. For $f$ as in the assumptions of the lemma we have a constant $C_k,D_k$ such that
\begin{align*}
f(\widetilde{X}(t))&\leq \frac{A}{(1-\alpha)^{\frac{k}{2}}}+ C_k\left(x^k+\Ex_x\left[\max_{0\leq s\leq t}\left|\int_0^s \sigma(\widetilde{X}(z))dW(z)\right|^k\right]\right)\\
&\leq \frac{A}{(1-\alpha)^{\frac{k}{2}}}+D_k\left(x^k+t^{\frac{k}{2}}\right),
\end{align*} where the second-to-last inequality follows from the Burkholder-Davis-Gundy inequality applied to the zero-mean martingale $\int_0^t \sigma(\widetilde{X}(s))dW(s)$ and using the assumed bound on $\sigma$. Integrating with respect to $e^{-(1-\alpha)s}ds$ we then have
\begin{equation*}
\begin{split}
|\hV(x)|\leq& \int_{0}^{\infty}e^{-(1-\alpha) s}\left(\frac{A}{(1-\alpha)^{\frac{k}{2}}}+D_k(x^k+s^{\frac{k}{2}})\right)ds\\
\leq &  \frac{A}{(1-\alpha)^{\frac{k+2}{2}}}+D_k\frac{x^k}{1-\alpha}+\frac{D_k}{1-\alpha}\int_{0}^{\infty}(1-\alpha)e^{-(1-\alpha) s}s^{\frac{k}{2}}ds
\\
\leq &  \frac{A}{(1-\alpha)^{\frac{k+2}{2}}}+D_k\frac{x^k}{1-\alpha}+E_k\frac{1}{(1-\alpha)^{\frac{k+2}{2}}}.
\end{split}
\end{equation*} The integral in the second line is just the $k/2$ moment of an $exp(1-\alpha)$ random variable and hence the last inequality. Taking $M=\max\{A+E_k,D_k\}$ concludes the proof. \eProof\vspace*{0.5cm}

We conclude the appendix with an example that shows the relationship discussed in \S \ref{sec:HT} in the context of a control problem rather than the performance analysis (fixed control) in that section. We consider here the controlled version of Example \ref{example:SSQ} but we take a different approach. Rather than directly ``attacking'' the Taylored equation, we build on priori knowledge about queues, specifically the fact that the optimal service rate will indeed place the queue in heavy-traffic and that, as in \S \ref{sec:HT}, the natural relationship between time scaling (as reflected in the discount factor) and utilization is $\sqrt{1-\alpha}\approx (1-\rho)$. Scaling space $(1-\rho)$ is the same as shrinking it by $\sqrt{1-\alpha}$.

\begin{example} \label{example:pdescaling}
	
	In the setting of Example \ref{example:SSQ0}, suppose the state space is restricted to $[0, M]$ and the transition probability at $M$ is $P^u_{M, M-1}=1.$ We take $\mathcal{U}(x)\equiv \mathbb{D}=[0,1]\bigcap\mathbb{Q}$.
	
	We will allow the truncated $M$ to scale with $\alpha$, specifically that $$M(\alpha)=\Upsilon/\sqrt{1-\alpha},$$  for $\Upsilon>0$. We will scale the space so that we study the (controlled chain) $Y=(\sqrt{1-\alpha}X_t,t\geq 0)$ instead of $X=(X_t,t\geq 0)$. The scaled chain's states space is then $[0,\Upsilon]$.
	
	From a feasible policy $U$, we can define $\delta( y)=\frac{1}{\sqrt{1-\alpha}}\left(U\left(\frac{y}{\sqrt{1-\alpha}}\right)-\frac{1}{2}\right),$ in which case the control set becomes
$$\widehat{\mathcal{U}}\equiv \left[-\frac{1}{2\sqrt{1-\alpha}},\frac{1}{2\sqrt{1-\alpha}}\right]\bigcap\mathbb{Q};$$ we are assuming here that $\sqrt{1-\alpha}$ is a rational number. The drift of $Y$ can be written in terms of $\delta$ as
$$\widehat{\mu}_{\delta}(y)=\sqrt{1-\alpha}\Ex_x^{U(x)}[X_1-x]=\sqrt{1-\alpha}\mu_u(x)=\sqrt{1-\alpha}(1-2u)=-2(1-\alpha)\delta(y).$$

The diffusion coefficient satisfies
$$\widehat{\sigma}^2_{\delta}(y)=(1-\alpha)\Ex_x^{U(x)}[(X_1-x)^2]=(1-\alpha)\sigma^2_u(x)\equiv (1-\alpha).$$

As in Example \ref{example:SSQ}, we take the reward function $r_u(x)=-x^4-\frac{c_s}{1-u}$. Scaling space and controls, this function translates to
$-\frac{y^4}{(1-\alpha)^2}-\frac{c_s}{\frac{1}{2}-\sqrt{1-\alpha}\delta}.$ Let us take the reward function

$$\widehat{r}_{\delta}(y)=-(1-\alpha)y^4-(1-\alpha)^3\frac{c_s}{\frac{1}{2}-\sqrt{1-\alpha}\delta}, ~y\in [0,\Upsilon],~ \delta \in \widehat{\mathcal{U}}.$$ The reward function $\tilde{r}_{\delta}(y):=\frac{\widehat{r}_{\delta}(y)}{1-\alpha},$ does not, then, grow with $\alpha$ for fixed $\delta$.
	
	The Taylored equation for the controlled problem for $(Y_t,t\geq 0)$ is then given by
	\begin{equation}\label{eqn:HJB-scaled}
		\begin{split}
			0= &\max_{\delta\in \widehat{\mathcal{U}}}\left\{\widehat{r}_{\delta}(y)-2(1-\alpha)\delta(y)V'(y)+\frac{1}{2}(1-\alpha)V''(y)\right\}- (1-\alpha)V(y), \quad\mbox{for}\quad y\in(0, \Upsilon),\\
			0= & V'(0), \\
			0= & V'(\Upsilon),
		\end{split}
	\end{equation}
	which, upon dividing by $(1-\alpha)$, becomes
			\begin{equation}\label{eqn:HJB-scaled}
	\begin{split}
	0= &\max_{\delta\in \widehat{\mathcal{U}}}\left\{\tilde{r}_{\delta}(y)-2\delta V'(y)+\frac{1}{2}V''(y)\right\}- V(y), \quad\mbox{for}\quad y\in(0, \Upsilon),\\
		0= & V'(0), \\
		0= & V'(\Upsilon),
		\end{split}
		\end{equation}
	
Notice that all of the reward function, drift, and diffusion coefficient are bounded and Lipschitz continuous, so that the existence of a solution is guaranteed by means of the same results underlying Lemma \ref{lem:PDEbounds}.
It is easy to argue that this solution has $V'(y)\geq \gamma$ for some $\gamma>0$ and all $y$ and, in turn, that given a constant $\Theta$, $\delta_*(y)\in [-\Theta,\Theta]$ for all $y$.
Thus, we may restrict $\delta$ to this range.
Let $\widehat{V}_*$ be the solution to the equation with $\delta$ now restricted; then, all of the reward, drift and diffusion coefficient are bounded and Lipschitz continuous uniformly in $\alpha$ over the bounded domain.
In turn, by the same results we used before, there exists a constant $C$ (depending on $\Upsilon$, $c_s$ and $\Theta$) such that
	\begin{equation}\label{eq:scaledbound}
[D^2\widehat{V}_*]_{\theta, \Omega_{M(\alpha)}}\leq C,
	\end{equation}
	 and if we define \be \bar{V}^{\alpha}_*(x)=\frac{\hV_*(\sqrt{1-\alpha}x)}{(1-\alpha)^3}\label{eq:barVdefin}\ee and $\bar{U}_*(x)=\sqrt{1-\alpha}\delta_*(\sqrt{1-\alpha}x)+\frac{1}{2}$, then $(\bar{U}_*, \bar{V}^{\alpha}_*)$ solves the TCP for the original chain:
	\begin{equation}\label{eqn:HJB-unscaled}
		\begin{split}
			0=&\max_{u\in \mathcal{U}(x)}\left\{r_u(x)+\alpha\left((1-2\mu_u(x))\widehat{V}'(x)+\frac{1}{2}\widehat{V}''(x)\right)\right\}-(1-\alpha) \widehat{V}(x),\quad x\in(0, M),\\
			0=& V'(0),  \\
			0=& V'(M).
		\end{split}
	\end{equation}
	
	From \eqref{eq:barVdefin} it then follows that
	\begin{equation*}
		[D^2\bar{V}^{\alpha}_*]_{\theta,\Omega_{M(\alpha)}}^*\leq \frac{C}{(1-\alpha)^{2-\frac{\theta}{2}}}.
	\end{equation*} Because $V_{\ast}^{\alpha}(x)\geq \frac{1}{(1-\alpha)^3}$ for all $x\geq \frac{1}{\sqrt{1-\alpha}}$ (see the proof of Corollary \ref{cor:orderopt}), we conclude from Theorem \ref{thm:diffusion} that
	\begin{equation*}
		|V_{\ast}^{\alpha}(x)-\bar{V}^{\alpha}_*(x)|\leq \sum_{t=0}^{\infty}\alpha^t [D^2\bar{V}^{\alpha}_*]_{\theta,\Omega_{M(\alpha)}}^*\leq \frac{C}{(1-\alpha)^{3-\frac{\theta}{2}}}\leq \sqrt{1-\alpha}^{\theta}V_{\ast}^{\alpha}(x),
	\end{equation*} for all $x\geq \frac{1}{\sqrt{1-\alpha}}$.
	\hfill \bsq
\end{example}

\bibliographystyleapp{ormsv080}
\bibliographyapp{central}
 \end{document}